\documentclass[moor,sglanonrev]{informs4}
\RequirePackage{tgtermes}
\RequirePackage{newtxtext}
\RequirePackage{newtxmath}
\RequirePackage{bm}
\RequirePackage{endnotes}

\OneAndAHalfSpacedXII % Current default line spacing
\usepackage{algorithm}
\usepackage{algpseudocode}
\usepackage{tikz}
\usepackage[sort&compress]{natbib}
 \bibpunct[, ]{[}{]}{,}{n}{}{,}%
 \def\bibfont{\small}%
\newcommand{\tabitem}{~~\llap{\textbullet}~~}

\EquationsNumberedThrough    % Default: (1), (2), ...
\TheoremsNumberedThrough     % Preferred (Theorem 1, Lemma 1, Theorem 2)
\ECRepeatTheorems  %  

\MANUSCRIPTNO{MOOR-0001-2024.00}

\begin{document}
%%%%%%%%%%%%%%%%

% Outcomment only when entries are known. Otherwise leave as is and
%   default values will be used.
%\setcounter{page}{1}
%\VOLUME{00}%
%\NO{0}%
%\MONTH{Xxxxx}% (month or a similar seasonal id)
%\YEAR{0000}% e.g., 2005
%\FIRSTPAGE{000}%
%\LASTPAGE{000}%
%\SHORTYEAR{00}% shortened year (two-digit)
%\ISSUE{0000} %
%\LONGFIRSTPAGE{0001} %
%\DOI{10.1287/xxxx.0000.0000}%

% Author's names for the running heads
% Sample depending on the number of authors;
% \RUNAUTHOR{Jones}
% \RUNAUTHOR{Jones and Wilson}
% \RUNAUTHOR{Jones, Miller, and Wilson}
% \RUNAUTHOR{Jones et al.} % for four or more authors
% Enter authors following the given pattern:
%\RUNAUTHOR{}
\RUNAUTHOR{Cao et al.}

% Title or shortened title suitable for running heads. Sample:
% \RUNTITLE{Predictive Maintenance in Manufacturing}
% Enter the (shortened) title:
\RUNTITLE{Optimal Partition for a Multi-Type Queueing System}

% Full title. Sample:
% \TITLE{Optimal Resource Allocation in Humanitarian Logistics: A Stochastic Programming Approach}
% Enter the full title:
\TITLE{Optimal Partition for a Multi-Type Queueing System}

% Block of authors and their affiliations starts here:
% NOTE: Authors with same affiliation, if the order of authors allows,
%   should be entered in ONE field, separated by a comma.
%   \EMAIL field can be repeated if more than one author
\ARTICLEAUTHORS{%
%\AUTHOR{John Doe,\textsuperscript{a} Jane Smith,\textsuperscript{b}}
%\AFF{\textsuperscript{a}Department of Industrial Engineering, University of XYZ, \EMAIL{john.doe@xyz.edu; \textsuperscript{b}Department of Computer Science, University of ABC, \EMAIL{jane.smith@abc.edu}} 
\AUTHOR{Shengyu Cao}
\AFF{Rotman School of Management, University of Toronto, Canada, M5S3E6, \EMAIL{shengyu.cao@rotman.utoronto.ca}}

\AUTHOR{Simai He}
\AFF{Antai College of Economics and Management, Shanghai Jiaotong University, China, 200030, \EMAIL{simaihe@sjtu.edu.cn}}

\AUTHOR{Zizhuo Wang}
\AFF{School of Data Science, The Chinese University of Hong Kong, Shenzhen, China, 518172, \EMAIL{wangzizhuo@cuhk.edu.cn}}

\AUTHOR{Yifan Feng}
\AFF{Department of Analytics and Operations, NUS Business School, National University of Singapore, Singapore, 119245, \EMAIL{yifan.feng@nus.edu.sg}}

% Enter all authors
} % end of the block

\ABSTRACT{%
% Enter your abstract
We study an optimal server partition and customer assignment problem for an uncapacitated FCFS queueing system with heterogeneous types of customers. Each type of customers is associated with a Poisson arrival, a certain service time distribution, and a unit waiting cost.
The goal is to minimize the expected total waiting cost by partitioning the server into sub-queues, each with a smaller service capacity, and routing customer types probabilistically.
First, we show that by properly partitioning the queue, it is possible to reduce the expected waiting costs by an arbitrarily large ratio. Then, we show that for any given server partition, the optimal customer assignment admits a certain geometric structure, enabling an efficient algorithm to find the optimal assignment. Such an optimal structure also applies when minimizing the expected sojourn time. Finally, we consider the joint partition-assignment optimization problem. The customer assignment under the optimal server partition admits a stronger structure. Specifically, if the first two moments of the service time distributions satisfy certain properties, it is optimal to \textit{deterministically} assign customer types with consecutive service rates to the same sub-queue. This structure allows for more efficient algorithms. Overall, the common rule of thumb to partition customers into continuous segments ranked by service rates could be suboptimal, and our work is the first to comprehensively study the queue partition problem based on customer types.
}%

%Supplemental Material:
%Data Ethics & Reproducibility Note:

% Sample
%\KEYWORDS{Stochastic programming, Decision support,Uncertainty, Disaster response, Optimization}

% Fill in data. If unknown, outcomment the field
%\HISTORY{Received: Month DD, YYYY; Accepted: Month DD, YYYY; Published Online: Month DD, YYYY}

\maketitle
%%%%%%%%%%%%%%%%%%%%%%%%%%%%%%%%%%%%%%%%%%%%%%%%%%%%%%%%%%%%%%%%%%%%%%

% Text of your paper here
\section{Introduction}
\label{sec:orga4d462c}
Queues are prevalent in daily activities. The analysis of queue pooling and partitioning has important ramifications in improving the performance of queueing systems. Consequently, whether pooling the queue will reduce the average waiting time has been extensively studied in the literature in the contexts of emergency departments \citep[see, e.g.,][]{SaghafianHVDK12}, call centers \citep[see, e.g.,][]{JouiniDN08}, computer server farms \citep[see, e.g.,][]{Harchol-BalterSY09}, among others. A classical conclusion is that when customer service time distributions are homogeneous, a pooled queue is always more appealing \citep{SmithW81}. Otherwise, dedicated queues can outperform in some cases \citep{Whitt99}. Indeed, separating arrivals to different queues according to different service times (or service time distributions) is ubiquitous in practice, including express lines in supermarkets, fast track in emergency departments, and size interval task assignment (SITA) policy employed in computer server farms. Notwithstanding its prevalent use, very few results rigorously discuss how to make an optimal partition of queues based on the service time (distribution) of different types of customers. Such insufficient study of queue partitioning, owing partly to the difficulty of analysis, motivates this work.

In this paper, we consider a service system with several types of customers, each with exponential arrival time and a certain service time distribution. We consider the scenario where a single server can be partitioned to multiple servers with different serving capacities, and each type of customers can be assigned to each partition with a certain probability. In our work, we are concerned with the expected waiting time of all customers. We aim to answer the following questions in our paper.
\begin{enumerate}
	\item Can server partition lead to a shorter expected waiting time?
	\item What is the optimal way to assign each type of customers with \textit{given} server partition?
	\item What is the optimal way to determine the server partition and assignment rule jointly?
\end{enumerate}

We build a queueing model with two sets of decisions to answer these questions. One is allocating the server resources (we call the partition decision), and the other is assigning each type of customers (we call the assignment decision).
First, we show that by properly partitioning the queue and making proper assignments, it is possible to reduce customers' expected total waiting time, and the benefit of partition could be arbitrarily large. Then, we consider two settings of this problem, depending on whether the partition is given (thus only assignment decision is to be made) or not.

When the partition of the servers is given, the optimal assignment is related to a convex partition of points in a two dimensional space defined by the first and second order moments of service time for each customer type. In particular, when there are only two queues, the optimal assignment rule is equivalent to a linear segmentation on the two dimensional space. In the joint partition-assignment setting, we further show that under certain conditions, the optimal solution structure is to partition customers into sub-intervals based on their service rates and assign customers within the same sub-interval to the same server, where customers of the same type will not be assigned to different servers.

We also consider a simplified case where the service time follows an exponential distribution. When the partition of the servers is given, we show that it is sometimes optimal to pool the fastest customers with the slowest in the same queue. This counter-intuitive finding is in line with \citet{feng2005optimal} and emphasizes that the common rule of thumb which partitions customers based on their service rate, ranked from slowest to fastest, may be far from optimum. Specifically, when there are two servers with a given serving capacity, the optimal assignment contains at most three continuous segments ranked by the customers' service rate, where customers in the first and third segments are assigned to one server while customers in the second segment are assigned to the other.
In the joint partition-assignment optimization problem, we show that it is always optimal to partition customers into sub-intervals based on their service rates, where customers within the same sub-interval are assigned to the same server. Moreover, customers of the same type will be assigned to the same server under the optimal assignment.

For better illustration, we summarize our results and main insights in two tables. In Table \ref{table:summary of results for k=2}, we summarize our main results for $ k = 2 $, i.e., there are two queues. We then present our main results in full generality in Table \ref{table:summary of results for general k}, where we consider a general $ k > 2 $.

\renewcommand{\arraystretch}{1.6}
\setlength{\tabcolsep}{12pt}
\begin{table}[ht!]
	\caption{\textbf{Summary of main results for two-queues $ (k=2) $.}} \label{table:summary of results for k=2}
	\begin{center}
		\begin{tabular}{|l || c |  c  | c |}
			\hline 
			%& Direction 
			& Assignment Only  & Partition and Assignment\\ %[0.5ex] 
			\hline
			
			\hline
			\begin{tabular}{@{}c@{}}
				\textit{General} \\
				service time \\
				distribution
			\end{tabular} & %$ \leq $ &
			\renewcommand{\arraystretch}{1.2}
			\begin{tabular}{@{}l@{}}
				\noalign{\medskip}
				\multicolumn{1}{c}{(Theorem \ref{thm:fix-alpha})}   \\
				\hline
				\textit{Optimal assignment structure:}\\
				%\tabitem Configuration set based \\
				\tabitem Ranked by a linear combination  \\
				$ \quad $ of the first two moments of  \\
				$ \quad  $   service time distribution\\
				\hline 
				\textit{Enumeration complexity:} $ O(n^2) $ \\
				\noalign{\medskip}
			\end{tabular}
			& 
			\renewcommand{\arraystretch}{1.2}
			\begin{tabular}{@{}l@{}}
				\multicolumn{1}{c}{(Theorem \ref{thm:opt-alpha})$ ^\diamondsuit$ } \\
				\hline
				\textit{Optimal assignment structure:}\\
				\tabitem Consecutive-service-rate based 
				\\
				\tabitem At most 2 segments \\
				\tabitem No fractional solutions \\
				\hline
				\textit{Enumeration complexity:} $ O(n)$
			\end{tabular} \\
			
			\hline
			\begin{tabular}{@{}c@{}}
				\textit{Exponential} \\
				service time \\
				distribution
			\end{tabular}
			%& $ \leq $	
			& \renewcommand{\arraystretch}{1.2}
			\begin{tabular}{@{}l@{}}  
				\noalign{\medskip}
				\multicolumn{1}{c}{(Theorem \ref{thm:fix-alpha-exp})}   \\
				\hline
				\textit{Optimal assignment structure:}\\
				\tabitem Service rate based 
				\\
				\tabitem At most 3 segments \\
				\tabitem At most 1 fractional solution \\
				\hline
				\textit{Enumeration complexity:} $ O(n^2) $ \\
				\noalign{\medskip}
			\end{tabular}
			& 
			\renewcommand{\arraystretch}{1.2}
			\begin{tabular}{@{}l@{}}
				\multicolumn{1}{c}{(Corollary \ref{cor:opt-alpha-exp})}    \\
				\hline
				\textit{Optimal assignment structure:}\\
				\tabitem Consecutive-service-rate based 
				\\
				\tabitem At most 2 segments \\
				\tabitem No fractional solutions \\
				\hline
				\textit{Enumeration complexity:} $ O(n) $ 
			\end{tabular}
			
			\\
			\hline
		\end{tabular}
	\end{center}\vspace{0.2cm}
	\footnotesize
	%\hspace{0.5cm}
	
	$\diamondsuit$ We need an additional condition on the first two moments of the service distribution described in Assumption \ref{as:reg-2seg} (e.g., exponential distribution) to obtain this result.
\end{table}

\renewcommand{\arraystretch}{1.6}
\setlength{\tabcolsep}{12pt}
\begin{table}[ht!]
	\caption{\textbf{Summary of main results for multi-queues $ (k > 2) $.}} \label{table:summary of results for general k}
	\begin{center}
		\begin{tabular}{|l || c |  c  | c |}
			\hline 
			%& Direction 
			& Assignment Only  & Partition and Assignment\\ %[0.5ex] 
			\hline
			
			\hline
			\begin{tabular}{@{}c@{}}
				\textit{General} \\
				service time \\
				distribution
			\end{tabular}
			& %$ \leq $ &
			\renewcommand{\arraystretch}{1.2}
			\begin{tabular}{@{}l@{}}
				\noalign{\medskip}
				\multicolumn{1}{c}{(Theorem \ref{thm:multi-fixed-alpha})}   \\
				\hline
				\textit{Optimal assignment structure:}\\
				%\tabitem Configuration set based\\
				\tabitem Based on ``configuration sets", \\
				$ \quad $ constructed from the first two  \\
				$ \quad $ moments of service time \\
				\hline 
				\textit{\small Enumeration complexity:} $ O(n^{k(k-1)}) $ \\
				\noalign{\medskip}
			\end{tabular}
			& 
			\renewcommand{\arraystretch}{1.2}
			\begin{tabular}{@{}l@{}}
				\multicolumn{1}{c}{(Theorem \ref{thm:multi-opt-alpha})$ ^\diamondsuit$ } \\
				\hline
				\textit{Optimal assignment structure:}\\
				\tabitem Consecutive-service-rate based 
				\\
				\tabitem At most $ k $ blocks \\
				\tabitem No fractional solutions \\
				\hline
				\textit{\small Enumeration complexity:} $ O(n^{k-1})$
			\end{tabular} \\
			
			\hline
			\begin{tabular}{@{}c@{}}
				\textit{Exponential} \\
				service time \\
				distribution
			\end{tabular}
			%& $ \leq $	
			& \renewcommand{\arraystretch}{1.2}
			\begin{tabular}{@{}l@{}}  
				\noalign{\medskip}
				\multicolumn{1}{c}{(Theorem \ref{thm:multi-exp})}   \\
				\hline
				\textit{Optimal assignment structure:}\\
				\tabitem Service rate based 
				\\
				\tabitem At most $ (2k-1) $ blocks \\
				%\tabitem At most $ k $ fractional solutions \\
				\hline
				\textit{\small Enumeration complexity:} $ O(n^{2k-2}) $ \\
				\noalign{\medskip}
			\end{tabular}
			& 
			\renewcommand{\arraystretch}{1.2}
			\begin{tabular}{@{}l@{}}
				\multicolumn{1}{c}{(Theorem \ref{thm:multi-exp})}    \\
				\hline
				\textit{Optimal assignment structure:}\\
				\tabitem Consecutive-service-rate based 
				\\
				\tabitem At most $ k $ blocks \\
				\tabitem No fractional solutions \\
				\hline
				\textit{\small Enumeration complexity:} $ O(n^{k-1}) $ 
			\end{tabular}
			
			\\
			\hline
		\end{tabular}
	\end{center}\vspace{0.2cm}
	\footnotesize
	%\hspace{0.5cm}
	
	$\diamondsuit$ We need an additional condition on the first two moments of the service distribution described in Assumption \ref{as:reg-2seg} (e.g., exponential distribution) to obtain this result.
\end{table}

We also study a few extensions.
First, when different waiting costs for different customer types are considered, we show in  Theorem \ref{thm:multi-opt-alpha-cost} that the minimization of expected waiting costs can be equivalently converted to a queue partition and assignment problem where unit waiting costs are identical for each customer type. In particular, when both partition and assignment can be optimized, we show that under certain conditions, an optimal assignment will partition customers into sub-intervals based on the product of their costs per unit waiting time and service rates, where customers within the same sub-interval are assigned to the same server. Such a result builds a bridge between the queue partition problem and the priority sequencing problem in which the famous \(c\mu\) rule is shown to be optimal.

Finally, we extend our results to another performance measure, the expected sojourn time. We show in Proposition \ref{pro:improvement-s} that the gap between not partitioning versus the optimal partition and assignment could also be arbitrarily large. We also show that the optimal assignment with minimal expected sojourn time is also connected to the partition of points in a two-dimensional space in Theorem \ref{thm:multi-fixed-alpha-sojourn}.

To summarize, our work provides a comprehensive analysis of the queue partition problem based on customer type. We identify the computational complexity and the solution structure, and propose algorithms for each case of the problem. In particular, the common rule of thumb to partition customers into two continuous segments ranked by their service rates may be suboptimal. To the best of our knowledge, our results are the first to identify the optimal structure for the queue partition problem under general service distribution, which may help improve the efficiency of queueing systems in practice.

\vspace{0.4 cm }

\noindent \textbf{Roadmap.}  The rest of the paper is organized as follows. In Section \ref{sec:literature}, we review the related literature. We formulate the queue partition problem in Section \ref{sec:model}. In Section \ref{sec:general}, we analyze the queue partition problem for both identical and non-identical waiting cost settings and relate it to the \(c\mu\) rule. We reveal the optimal structure under general service time distribution and also highlight the special case of exponential distribution. We extend our results with sojourn time as the performance metric in Section \ref{sec:sojourn}. Finally, we conclude in Section \ref{sec:conclusion}.

\section{Literature Review}
\label{sec:literature}
Our work is closely related to the literature that studies queue partitions under various settings.
Among these studies, two partition decisions, including arrival partition (also known as assignment or routing decision) and service partition (also known as capacity allocation decision), are considered. One stream considers the arrival partition, where the partition is done according to the job size. In particular, they assume that upon arrival, the service time is known. Based on the known service time, assignment decisions are made. Under such an assumption, a class of policy called the size interval task assignment (SITA) policy is proposed. A seminal result in this stream is by \citet{Harchol-BalterSY09}, who state that when traffic is heavy and job-size distributions have high variability, the SITA policy outperforms the policy that assigns incoming jobs to the queue with the least total job size remaining, which is called the Least-Work-Left (LWL) policy. \citet{feng2005optimal} later generalize it to a nested size interval strategy and proves its optimality. Such partition based on the actual job size is reasonable for deterministic computer workloads. However, in fields of emergency departments, call centers, and cloud computing system with general tasks (e.g., machine learning or optimization tasks), service time is often unknown before the service is completed. Instead, it is common that we only have the type information, i.e., only the distribution of service time is known upon customers' arrival. In this paper, we assume that the type information is known but the actual service time is unknown, and study the partition problem under this setting. Since deterministic distribution for each type is equivalent to known sizes with a discrete-support job size distribution, our work generalizes nested size interval strategies to general service distributions for each arrival type.

In the queueing literature, there are two types of models characterizing how the service can be pooled or partitioned. In the first type of model, the decision maker determines the number of servers assigned to each sub-queue \citep[see, e.g.,][]{Whitt99,HungP07,HuB09}, while in the second type of model, the decision maker decides the service rate assigned to each sub-queue \citep[see, e.g.,][]{IyerJ04,HassinSY15,YuBG15}. When pooling service rates, two single-server queues are merged into a new single-server queue, where service rates are aggregated \citep[see, e.g.,][]{Kleinrock76,MandelbaumR98,AndradottirAD17}. {The partition on service rates rather than servers is common in the application of service facility, business process management, conveyor and storage systems, cloud computing and web server environments \citep{AllonF08,dieker2017optimal,hassin2019delay}. Our work adopts the second model and assumes that the service rates can be pooled or partitioned. Nevertheless, the findings may still provide valuable perspectives for server partitioning systems. For example, \citet{MandelbaumR98} point out that when traffic is heavy, the performance of pooling service rates coincides with that of pooling servers.}

Recently, researchers have analyzed the effect of pooling under more sophisticated settings. For instance, \citet{ArgonZ09} discuss how imperfect classification influences the decision process of partitioning. \citet{SunarTZ21} analyze the case where customers are delay-sensitive and discuss the benefit of pooling. \citet{CaoHHL20} argue that with a proper routing policy, idle time of servers in dedicate systems can be reduced. In another work, \citet{HuB09} study how to partition the queueing systems during rush hour, where a plethora of customer arrivals occur in a short time window and few or even no customers appear thereafter. They assume instant arrival of customers and use single-queue formula to analyze the best partition decision. They prove that separating each customer type is optimal, and give the optimal allocation of servers to each customer type. In contrast with this work, we consider the optimal partition in a stationary environment. Also, we assume the server capacity is divisible. Our goal is to analyze whether a subset of customer types should be pooled and how many server resources should be assigned to this group of customer types. {\citet{hassin2019delay} explore the capacity partitioning within a GI/M/$\infty$ system. They aim to minimize the expected sojourn times by allocating a predefined service capacity among an infinite number of sub-queues. In their model, allocating substantial service capacity to a few sub-queues ensures rapid service for customers in these queues, but if these faster queues are all occupied, new customers are relegated to the remaining queues with significantly slower service rates. Optimizing this allocation involves balancing these two factors. Contrasting with their model, we assume customer assignment is based on their type information and each queue follows a first-come-first-served (FCFS) principle.}
The work that is closest to ours is \citet{Whitt99}. In this work, the author points out that if the queues are pooled, then the economies of scale will increase the service resource utilization and minimize the idle time of the server. However, when the variation of the service time distributions is very large, separating fast customers from the others may save them from being blocked and offset the lower utilization disadvantage, thus increasing the overall efficiency. \citet{Whitt99} considers when and how to partition a pool of identical indivisible servers into sub-groups to minimize the overall number of servers required to meet the underlying requirement of system delay. However, the analysis is based on a heuristic method and does not analyze the structure of the optimal partition. We complement this work and establish the optimal structure of customer assignment by rigorously quantifying the trade-off between servers being under-utilized and fast jobs getting blocked by slow jobs.

Another line of work that is related to ours is the study of stochastic scheduling problems. This literature aims to minimize the average waiting cost by optimizing the allocation sequence of a single resource to \(n\) parallel queues. A classical result is that a \(c\mu\) decision rule is optimal. More specifically, the \(c\mu\) rule orders the queues based on the product of the unit waiting cost and the service rate, and selects the queue with the largest product under various conditions. For exponential arrival processes and generally distributed service times, \citet{CoxS61} prove that the \(c\mu\) rule is optimal. Later, \citet{BarasMM85} prove its optimality for general arrival processes and geometrically distributed service times. Following works have extended the optimality of \(c\mu\) rule in more general settings \citep[see, e.g.,][]{HirayamaKN89,Mieghem95,MandelbaumS04}. Nevertheless, as mentioned by \citet{HuB09}, the priority sequencing problem can be impractical in some settings. Partitioning, as an alternative, can be used to boost the system's efficiency. Our work shows that the \(c\mu\) decision rule is also optimal for the queue partition problem under certain conditions, bridging the queue partition problem and the priority sequencing problem.

\section{Model}
\label{sec:model}
We consider a queueing system with \(n\) types of customers and \(k\) queues. Type \(i\) customers arrive at the system according to a Poisson process with an arrival rate \(\lambda_i\). The service time for type \(i\) customer under a unit server is \(S_i\) with \(\mathbb{E}[S_i] = \frac{1}{\mu_i}\) and \(\mathbb{E}[S_i^2] = v_i\). We assume the total serving capacity is \(1\) (this assumption is without loss of generality) and the serving capacity is divisible, meaning that it can be divided into several queues, each with serving capacity \(\alpha_j\) such that \(\sum_{j=1}^k \alpha_j = 1\). When a server with capacity \(\alpha_j\) serves the \(i\)th type of customers, the service time distribution becomes \(\frac{S_i}{\alpha_j}\). In the following, we consider the problem of dividing the service capacity and assigning customers to each sub-queue with their type information only. We assume the unit waiting cost for each customer is identical, and our objective is to minimize the expected total waiting time of all customers in the queueing system. In Sections \ref{sec:general} and \ref{sec:sojourn}, we will consider the case in which the waiting cost for each customer is heterogeneous and the case in which the objective is to minimize the expected total sojourn time respectively.

We consider two types of problems. In the first type of problem, we assume that the partition \(\boldsymbol{\alpha} = (\alpha_1, \alpha_2, \dots, \alpha_k)\) is given, and we are interested in finding the optimal assignment of each type of customers. In the second type of problem, we jointly decide the partition \(\boldsymbol{\alpha}\) and the assignment of each type of customers.

For each aforementioned problem, we consider a stochastic assignment policy where customers of the same type can be assigned to different queues, each with certain probabilities. Specifically, we use \(X\in[0, 1]^{n\times k}\) to indicate the assignment of each type of customer, where \(X_{ij}\) denotes the probability of the \(i\)th type of customers assigned to the \(j\)th queue. In other words, the \(i\)th row of \(X\) represents how the \(i\)th type of customers are assigned to each queue, and the \(j\)th column of \(X\) represents the proportion of each type of customers assigned to the \(j\)th queue.

We define \(R(\boldsymbol{\alpha}) = \{j: \alpha_j > 0\}\) to represent queues with positive serving capacity. Given partition \(\boldsymbol{\alpha}\) and assignment \(X\), we can provide a closed form expression for expected waiting time using the Pollaczek-Khinchine formula \citep[pg.~209]{Haigh13}:

$$
f(X, \boldsymbol{\alpha}) := \sum_{j\in R(\boldsymbol{\alpha})}
\left(\frac{\sum_{i=1}^n \lambda_i X_{ij}v_i}{\alpha_j^2 - \alpha_j\sum_{i=1}^n \frac{\lambda_i X_{ij}}{\mu_i}}
\cdot \frac{\sum_{i=1}^n X_{ij} \lambda_i}{2\sum_{i=1}^n \lambda_i}\right).
$$

Note that for \(f(X, \boldsymbol{\alpha})\) to be meaningful, we must have
$$
(X, \boldsymbol{\alpha}) \in \mathcal{F} = \left\{
(X, \boldsymbol{\alpha}):
\alpha_j > \sum_{i=1}^n \frac{\lambda_i X_{ij}}{\mu_i}, \text{ for } j \in R(\boldsymbol{\alpha})
\right\},
$$
or equivalently,
$$
X \in \mathcal{F}_{\boldsymbol{\alpha}} = \left\{
X:
\alpha_j > \sum_{i=1}^n \frac{\lambda_i X_{ij}}{\mu_i}, \text{ for } j \in R(\boldsymbol{\alpha})
\right\}.
$$

Classified by the type of problem (assignment only versus partition and assignment), we are interested in solving the following two problems:
\begin{itemize}
	\item The \emph{multi-queue assignment problem} for given \(\boldsymbol{\alpha}\):
	\begin{eqnarray}\label{eq:k-qa-w}
	(k\textup{-QAP})\quad\inf_{X} &&\quad f(X, \boldsymbol{\alpha}) \nonumber \\
	\text{s.t.} &&\quad X \in \mathcal{F}_{\boldsymbol{\alpha}}; \nonumber \\
	&&\quad \sum_{j=1}^k X_{ij} = 1, \quad i = 1, \dots, n; \nonumber \\
	&&\quad X_{ij}\in[0,1], \quad i = 1, \dots, n,\; j = 1, \dots, k.
	\end{eqnarray}
	\item The \emph{multi-queue partition (and assignment) problem}:
	\begin{eqnarray}\label{eq:k-qp-w}
	(k\textup{-QPP})\quad\inf_{X, \boldsymbol{\alpha}} &&\quad f(X, \boldsymbol{\alpha}) \nonumber \\
	\text{s.t.} &&\quad (X, \boldsymbol{\alpha}) \in \mathcal{F}; \nonumber \\
	&&\quad \sum_{j=1}^k \alpha_j = 1; \qquad  \alpha_j\ge 0, \quad j = 1, \dots, k; \nonumber \\
	&&\quad \sum_{j=1}^k X_{ij} = 1, \quad i = 1, \dots, n; \nonumber \\
	&&\quad X_{ij}\in[0,1], \quad i = 1, \dots, n,\; j = 1, \dots, k.
	\end{eqnarray}
\end{itemize}

Before we proceed, we show that, in general, proper partitioning and assignment can reduce the expected waiting time by an arbitrarily large amount. The precise result is stated in the proposition below.

\vspace{0.5 cm}
\begin{proposition}
	\label{pro:improvement}
	For any \(\epsilon > 0\), there exist input parameters \(\{\lambda_i\}_{i=1}^n, \{\mu_i\}_{i=1}^n, \{v_i\}_{i=1}^n\) and \(X, \boldsymbol{\alpha}\) such that \(f(X, \boldsymbol{\alpha}) < \epsilon f(\bar{X}, \bar{\boldsymbol{\alpha}})\), where the pooling decision is
	$$
	\bar{X} = \begin{pmatrix}
	0 &\cdots & 0 &1 \\
	\vdots &\vdots &\vdots &\vdots \\
	0 &\cdots & 0 &1 \\
	\end{pmatrix}, \quad
	\bar{\boldsymbol{\alpha}} = \begin{pmatrix}
	0 \\ \vdots \\ 0 \\ 1
	\end{pmatrix}.
	$$
\end{proposition}

To illustrate, consider \(k = 2\) queues and \(n = 2\) types of customers whose service times both follow exponential distributions. In particular, \(\lambda_1 = t, \lambda_2 = 1\) and \(\frac{\lambda_1}{\mu_1} = \frac{1}{2}, \frac{\lambda_2}{\mu_2} = \frac{t}{2}\) where \(0 < t \le \frac{1}{4}\). For partition \(\boldsymbol{\alpha} = (1 - t, t)\) and assignment
$$X = \begin{pmatrix}
1 &0\\
0 &1\\
\end{pmatrix},$$
we have \(f(X, \boldsymbol{\alpha}) = \left(\frac{1/4}{(1-t)(1/2-t)}+\frac{t^2/4}{t^2/2}\right)\frac{1}{1+t}\). We can also calculate \(f(\bar{X}, \bar{\boldsymbol{\alpha}}) = \frac{1/t+t^2}{2(1-t)}\). Therefore,
$$
\lim_{t\to0}\frac{f(X, \boldsymbol{\alpha})}{f(\bar{X}, \bar{\boldsymbol{\alpha}})} = \lim_{t\to0}\frac{1}{1/(2t)} = 0.
$$
The intuition is the negative externality of queues: When a plethora of customers who can be served quickly are blocked by some slow minority, a huge reduction to the average waiting time can be obtained by allocating some dedicated resources for the vast majority. Therefore, partitioning the queue can significantly improve the overall efficiency in this case.

We note that \citet{Whitt99} also shows that by properly partitioning the queue, it is possible to reduce the total average waiting time. (The model in \citealt{Whitt99} considers the split of servers, while in this paper, we consider the split of service capacity. Therefore, our model can be viewed as allowing a more general partition than the one in \citealt{Whitt99}.) In Proposition \ref{pro:improvement}, we further show that the improvement can be arbitrarily large under certain instances. This also implies that there does not exist a constant bound for the expected waiting time between the optimal split and a pooled server.

Both in the literature (e.g., \citealp{Whitt99}) and in practice, a rule of thumb is to partition customers into two continuous segments ranked by their service rates where all customers with high service rates are assigned to one queue, and all the rest of the customers are assigned to the other queue. However, in the following example, we show that the performance of this principle can be poor.

\vspace{0.3 cm}
\begin{example}
	\label{eg:3better}
	Consider three types of customers with \(\lambda_1 = 13, \lambda_2 = 0.008, \lambda_3 = 0.006\), \(\mu_1 = 20, \mu_2 = 0.3, \mu_3 = 0.2\), \(\mathbb{E}[S_1^2] = 0.01, \mathbb{E}[S_2^2] = 60, \mathbb{E}[S_3^2] = 25.4\). For $ \boldsymbol{\alpha} = (0.95, 0.05) $, if we partition the customers based on the expected service time, the best assignment is \(x_1 = 1, x_2 = 0.65, x_3 = 0\) with an expected waiting time of 1.03. However, for assignment \(x_1 = 1, x_2 = 0, x_3 = 1\), the expected waiting time is 0.68, which is 34\% less. Intuitively, the first and the third types of customers have lower service time variance. By grouping them, we can benefit from pooling while still having a small service time variance within that group. In contrast, the second type of customers has a large service time variance. Grouping this type of customers with any other type (either type one or type three) could block the other type of customers, resulting in a larger average waiting time. Therefore, it is better to group type one and type three customers in this case. \(\hfill\square\)
\end{example}
\vspace{0.3 cm}

In this example, we can see that not only the service rate but also the variation of service time has an impact on the optimal assignment. Consequently, the partition of two continuous segments fails to be optimal. In subsequent sections, we will analyze the structure of the above-described problems and propose efficient algorithms to solve the partition and assignment problem.

\section{Analysis}
\label{sec:general}

In this section, we analyze the optimal structures of the queue assignment problem and the queue partition problem when each service time distribution follows a general distribution, e.g., Erlang distribution \citep{NagH17}, Pareto distribution \citep{XuXP06}, or log-normal distribution \citep{SakovZ00}.

\subsection{The Assignment Problem for Two-Queues}
\label{ssec:qap-g}
In this subsection, we first consider the queue assignment problem in which the serving capacity can only be divided into two queues. For notational simplicity, we use \(\alpha\) to represent the serving capacity of the first queue. Consequently, the serving capacity of the second queue is \(1 - \alpha\). We also use \(\boldsymbol{x}\in[0, 1]^n\) to denote the assignment of each type of customer, where \(x_i\) denotes the probability of the \(i\)th type of customers assigned to the first queue. With slight abuse of notation, we let
\begin{align*}
f(\boldsymbol{x}, \alpha) :=
\begin{cases}
\frac{\sum\limits_{i=1}^n \lambda_i x_i v_i}{\alpha^2 - \alpha\sum\limits_{i=1}^n \frac{\lambda_i x_i}{\mu_i}} \cdot \frac{\sum\limits_{i=1}^n x_i \lambda_i}{2\sum\limits_{i=1}^n \lambda_i} + \frac{\sum\limits_{i=1}^n \lambda_i(1 - x_i)v_i}{(1 - \alpha)^2 - (1 - \alpha)\sum\limits_{i=1}^n \frac{\lambda_i(1 - x_i)}{\mu_i}} \cdot \frac{\sum\limits_{i=1}^n (1 - x_i)\lambda_i}{2\sum\limits_{i=1}^n \lambda_i} &\quad \text{ if } 0 < \alpha < 1 \\[18pt]
\frac{1}{2}\left(\sum_{i=1}^n \lambda_i v_i\right)/\left(1 - \sum_{i=1}^n \lambda_i/\mu_i\right) &\quad \text{ if } \alpha = 0, \boldsymbol{x} = \boldsymbol{0} \text{ or } \alpha = 1, \boldsymbol{x} = \boldsymbol{1},
\end{cases}
\end{align*}
where \(\boldsymbol{0}\) (\(\boldsymbol{1}\), respectively) is an all-\(0\) (all-\(1\), respectively) vector. 

A few words about the function \(f(\boldsymbol{x}, \alpha)\). First, for it to be meaningful, we must have
$$
(\boldsymbol{x}, \alpha) \in \mathcal{F} = \left\{
(\boldsymbol{x}, \alpha):
\sum_{i=1}^n \frac{\lambda_i x_i}{\mu_i} < \alpha < \sum_{i=1}^n \frac{\lambda_i x_i}{\mu_i} + 1 - \sum_{i=1}^n \frac{\lambda_i}{\mu_i}\right\} \bigcup \left\{(\boldsymbol{0}, 0), (\boldsymbol{1}, 1)\right\},
$$
or equivalently, \(\boldsymbol{x}\in\mathcal{F}_\alpha\) where
$$
\mathcal{F}_\alpha = 
\left\{
\boldsymbol{x}:
\sum_{i=1}^n \frac{\lambda_i x_i}{\mu_i} < \alpha < \sum_{i=1}^n \frac{\lambda_i x_i}{\mu_i} + 1 - \sum_{i=1}^n \frac{\lambda_i}{\mu_i}
\right\},
$$
if \(0 < \alpha < 1\) and \(\mathcal{F}_0 = \left\{\boldsymbol{0}\right\}, \mathcal{F}_1 = \left\{\boldsymbol{1}\right\}\). Second, since \(f(\boldsymbol{x}, \alpha) = f(\boldsymbol{1} - \boldsymbol{x}, 1 - \alpha)\), we may assume \(\alpha \ge \frac{1}{2}\) without loss of generality. We also assume \(\mu_1 > \mu_2 > \cdots > \mu_n\). Finally, Lemma \ref{lemma:convex} below reveals the convexity property of \(f(\boldsymbol{x}, \alpha)\) Its proof is by directly analyzing the second-order derivative of \(f(\boldsymbol{x}, \alpha)\) and thus is omitted. 
\begin{lemma}
	\label{lemma:convex}
	For \((\boldsymbol{x}, \alpha) \in \mathcal{F}\) and \(x_i \ge 0\), \(i = 1, \dots, n\), \(f(\boldsymbol{x}, \alpha)\) is convex in \(\alpha\) and is convex in each \(x_i\), for \(i = 1, \dots, n\).
\end{lemma}

We analyze the queue assignment problem given \(\alpha\):
\begin{eqnarray}\label{eq:qa-w}
(\textup{QAP})\quad\inf_{\boldsymbol{x}} &&\quad f(\boldsymbol{x}, \alpha) \nonumber \\
\text{s.t.} &&\quad\boldsymbol{x} \in \mathcal{F}_\alpha; \quad x_i\in [0, 1], \; i = 1, \dots, n.
\end{eqnarray}

Despite Lemma \ref{lemma:convex}, we note that \(f(\boldsymbol{x}, \alpha)\) is not jointly convex in \(\boldsymbol{x}\) and \(\alpha\) (it is not jointly convex in \(\boldsymbol{x}\) either). Therefore, directly minimizing \(f(\boldsymbol{x}, \alpha)\) may be challenging. To further characterize the structure of this problem, we introduce the following notation. For a set of service time \(\mathcal{I} = \{S_i\}_{i=1}^n\), define a piecewise linear function \(v_{\mathcal{I}}: [\mathbb{E}[S_1], \mathbb{E}[S_n]] \to \mathbb{R}_+\) such that for \(i = 1, \dots, n - 1\) and \(\mathbb{E}[S_i] \le u \le \mathbb{E}[S_{i+1}]\),
\begin{equation}\label{eq:v-u-relation}
v_{\mathcal{I}}(u) =
\frac{\mathbb{E}[S^2_{i+1}] - \mathbb{E}[S^2_i]}{\mathbb{E}[S_{i+1}] - \mathbb{E}[S_i]} (u - \mathbb{E}[S_i]) + \mathbb{E}[S_i^2].
\end{equation}
Intuitively, the graph of \(v_{\mathcal{I}}\) contains both the average service rate information (how fast can a customer type be served) and the service time variation information (how uncertain a customer type's service time is).

Given \(\mathcal{I} = \{S_i\}_{i=1}^n\), for a linear function \(l(u)\), consider the induced partition \(U_{\mathcal{I}, l} \cup E_{\mathcal{I}, l} \cup L_{\mathcal{I}, l} = \{1, 2, \dots, n\}\) where
\begin{align*}
U_{\mathcal{I}, l} &= \left\{i \in \{1, 2, \dots, n\}: v_{\mathcal{I}}(\mathbb{E}[S_i]) > l(\mathbb{E}[S_i])\right\}, \\
E_{\mathcal{I}, l} &= \left\{i \in \{1, 2, \dots, n\}: v_{\mathcal{I}}(\mathbb{E}[S_i]) = l(\mathbb{E}[S_i])\right\}, \\
L_{\mathcal{I}, l} &= \left\{i \in \{1, 2, \dots, n\}: v_{\mathcal{I}}(\mathbb{E}[S_i]) < l(\mathbb{E}[S_i])\right\}.
\end{align*}

\begin{figure}[!ht]
	\centering
	\begin{minipage}{.45\textwidth}
		\centering
		\includegraphics[height=5.6cm]{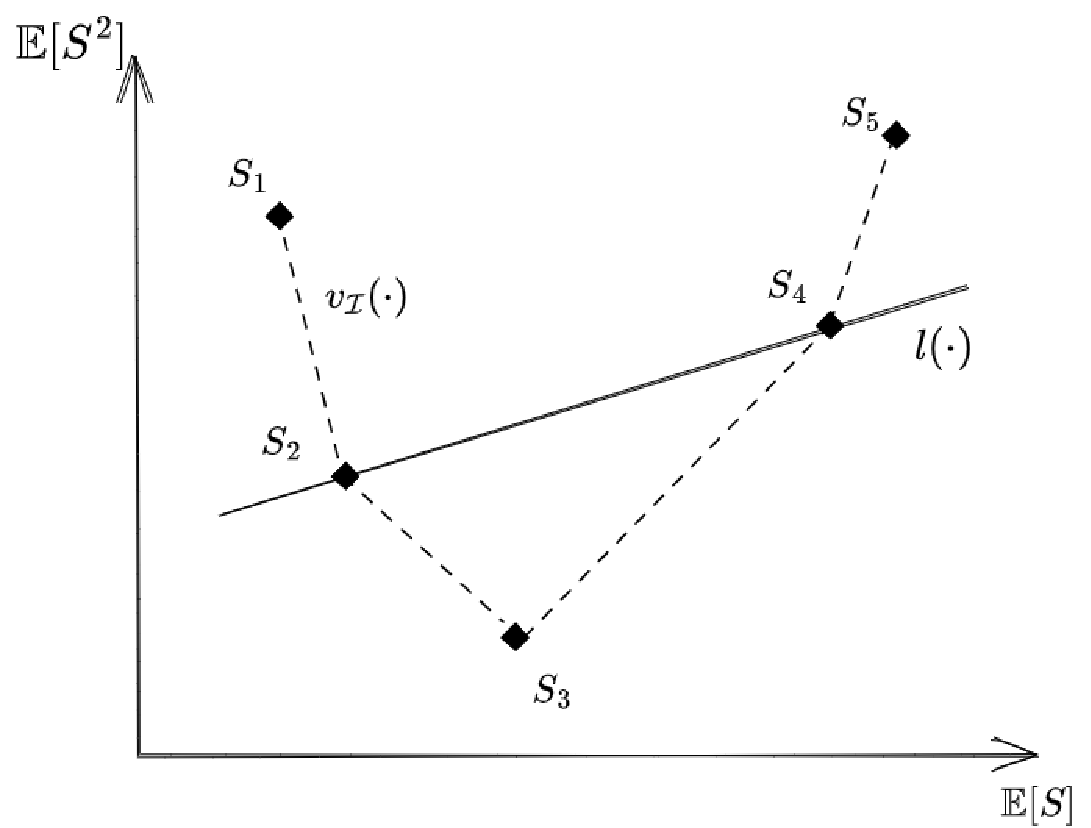}
		\caption{There are 5 types of customers, and the induced partition by $l(\cdot)$ is $U_{\mathcal{I}, l} = \{1, 5\}, E_{\mathcal{I}, l} = \{2, 4\}, L_{\mathcal{I}, l} = \{3\}$.}
		\label{fig:point-partition}
	\end{minipage}
	\hfill
	\begin{minipage}{.45\textwidth}
		\centering
		\includegraphics[height=5.6cm]{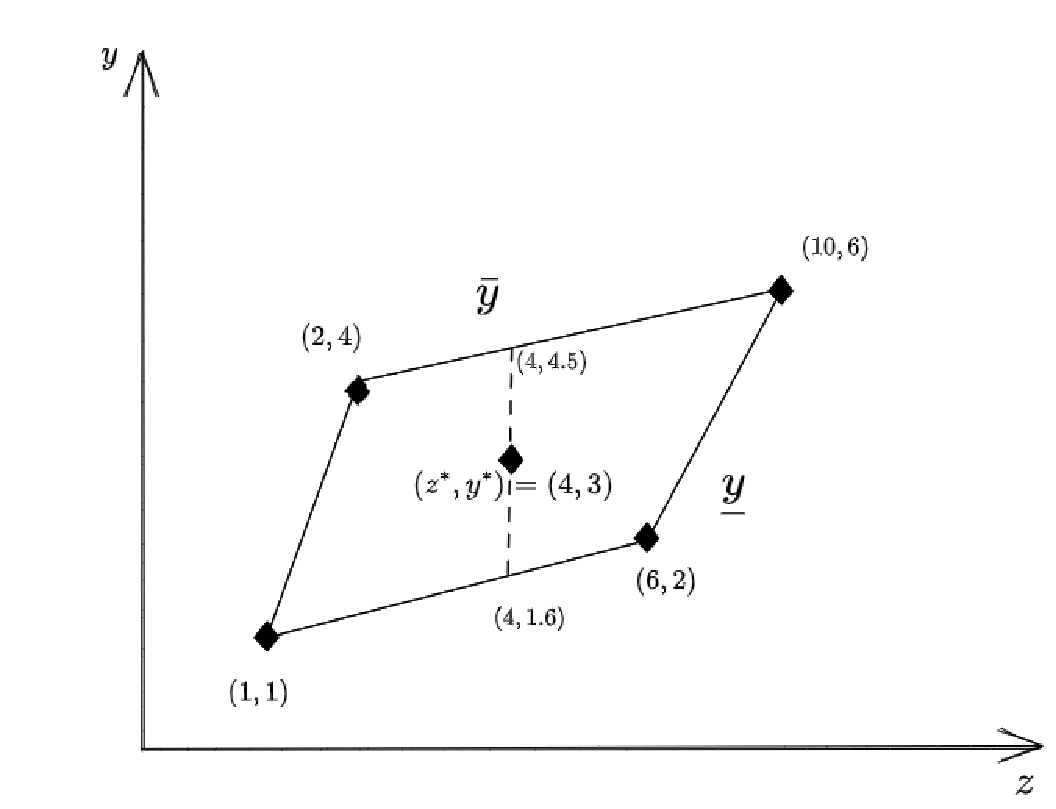}
		\caption{Recover optimal assignment from $(z^{ * }, y^{ * }) = (4, 3)$ with $E = \{1, 2\}$ where $\mu_1 > \mu_2$.}
		\label{fig:bivariate}
	\end{minipage}
\end{figure}

To illustrate, consider a QAP problem with \(n = 5\). As is shown in Figure \ref{fig:point-partition}, \(v_{\mathcal{I}}(\cdot)\) is a piecewise linear function with four pieces. For the linear function \(l(\cdot)\) shown in Figure \ref{fig:point-partition}, \(U_{\mathcal{I}, l} = \{1, 5\}\) is the set of customer types of which the service time has larger second moment compared with linear rule \(l(\cdot)\). \(L_{\mathcal{I}, l} = \{3\}\) is the set of customer types with smaller second moment compared with linear rule \(l(\cdot)\). \(E_{\mathcal{I}, l} = \{2, 4\}\) is the set of customer types whose second moment of the service time distribution is exactly prescribed by \(l(\cdot)\).

We define the configuration set of \(\mathcal{I}\) as
$$ \mathcal{C}(\mathcal{I}) := \left\{(U, E, L): (U, E, L) \text{ is the induced partition of some linear function } l, |E|\ge 2 \right\}. $$

Intuitively, the function \(l(u)\) provides a linear rule to classify service time distributions into three groups such that \(|E_{\mathcal{I}, l}| \ge 2\). Within \(U_{\mathcal{I}, l}\), service time has a larger variation scaled by its expectation, while \(L_{\mathcal{I}, l}\) contains those with a smaller variation. The configuration set \(\mathcal{C}(\mathcal{I})\) contains all possible classifications by such linear rules. Because every two points determine a line in a two-dimensional space, it is readily shown that \(|\mathcal{C}(\mathcal{I})| \le \binom{n}{2}\). The following theorem reveals the relation between the optimal solution of QAP and \(\mathcal{C}(\mathcal{I})\).

\vskip 0.3 cm
\begin{theorem}
	\label{thm:fix-alpha}
	Suppose \(\mathcal{I}\) with \(\mu_1 > \cdots > \mu_n, v_1, \dots, v_n\),\(\lambda_1, \dots, \lambda_n\), and $\alpha$ are given. For any optimal solution \(\boldsymbol{x}^*\) to the assignment problem, there exist \((U, E, L) \in \mathcal{C}(\mathcal{I})\) such that either
	\begin{align*}
	\begin{cases} x_i^{ * } = 1 & i\in U \\ x_i^{ * } = 0 & i\in L \\ x_i^{ * } \in [0, 1] & i\in E \end{cases} \qquad \text{ or } \qquad \begin{cases} x_i^{ * } = 0 & i\in U \\ x_i^{ * } = 1 & i\in L \\ x_i^{ * } \in [0, 1] & i\in E \end{cases}.
	\end{align*}
\end{theorem}
\vskip 0.3 cm

Theorem \ref{thm:fix-alpha} characterizes the structure of the optimal partition by connecting the partition of the customers to the partition of points in a two-dimensional space by a line. To illustrate, consider the example in Figure \ref{fig:point-partition}. The partition induced by function \(l(\cdot)\) corresponds to two candidates of optimal solution, \(x_1 = x_5 = 0, x_3 = 1\) and \(x_1 = x_5 = 1, x_3 = 0\). Theorem \ref{thm:fix-alpha} states that one of all these \(2|\mathcal{C}(\mathcal{I})|\) candidates must be optimal for the corresponding QAP problem.

Consequently, we can find the optimal assignment by enumerating all possible configurations, whose size is at most \(\mathcal{O}(n^2)\). Then, within each configuration \((U, E, L)\), we search for the best assignment for customers within set \(E\). If \(E\) contains exactly two points (which is usually likely to be the case), then solving the best assignment for customers within set \(E\) amounts to solving a bi-variate optimization problem. Next, we show that even if \(E\) contains more than two points, searching for the best assignment for customers within set \(E\) can still be converted to a bi-variate optimization problem, which can be solved easily. In particular, given \(x_j, j \in U \cup L\), for \(i \in E\), we assume \(\mathbb{E}[S_i^2] = k_1 \mathbb{E}[S_i] + k_2\). Then we can reduce decision variables \(x_i, i \in E\) to \(y = \sum_{i \in E}\frac{\lambda_i}{\mu_i} x_i\) and \(z = \sum_{i \in E}\lambda_i x_i\). We show that optimization for \(x_i, i \in E\) can be reduced to an optimization problem over \(y\) and \(z\). Under the definition of \(y\) and \(z\), we have \(\sum_{i \in E}\lambda_i v_i x_i = k_1 \sum_{i \in E}\frac{\lambda_i}{\mu_i}x_i + k_2 \sum_{i \in E}\lambda_i x_i = k_1 y + k_2 z\). Let \(E = \{i_1, i_2, \dots, i_{|E|}\}\) where \(i_1 < i_2 < \cdots < i_{|E|}\). The objective function can be rewritten as
\begin{align}
&\frac{\sum_{i \in L \cup U} \lambda_i x_i v_i + k_1 y + k_2 z}{\alpha^2 - \alpha\left(\sum_{i\in L \cup U} \frac{\lambda_i x_i}{\mu_i} + y\right)} \cdot \frac{\sum_{i \in L\cup U} x_i \lambda_i + z}{2\sum_{i=1}^n \lambda_i} \nonumber \\ +
&\frac{\sum_{i=1}^n \lambda_i v_i - \sum_{i \in L \cup U} \lambda_i x_i v_i - k_1 y - k_2 z}{(1 - \alpha)^2 - (1 - \alpha)\left(\sum_{i=1}^n \frac{\lambda_i}{\mu_i} - \sum_{i\in L \cup U} \frac{\lambda_i x_i}{\mu_i} - y\right)} \cdot \frac{\sum_{i=1}^n \lambda_i - \sum_{i \in L\cup U} x_i \lambda_i - z}{2\sum_{i=1}^n \lambda_i}. \label{eq:bivariate}
\end{align}
Fix \(z = \sum_{i \in E}\lambda_i x_i\), and let \(s\) be the largest integer such that \(\sum_{k=1}^{s-1} \lambda_{i_k} \le z\). Then the assignment \(\underline{x}_{i_1} = \cdots = \underline{x}_{i_{s-1}} = 1, 0 \le \underline{x}_{i_s} = \frac{z - \sum_{k=1}^{s-1}\lambda_{i_k}}{\lambda_{i_s}} < 1, \underline{x}_{i_{s+1}} = \cdots = \underline{x}_{i_{|E|}} = 0\) corresponds to the minimal value of \(y\) given \(z\), denoted by \(\underline{y}(z)\). By enumerating \(s\) from \(1\) to \(|E|\), we can see that \(\underline{y}(z)\) is a convex piecewise linear function with \(|E|\) pieces, and thus can be represented efficiently by \(|E|+1\) pairs of \((y, z)\). Similarly, let \(t\) be the smallest integer such that \(\sum_{k=t+1}^{|E|}\lambda_{i_k} \le z\). The assignment \(\bar{x}_{i_1} = \cdots = \bar{x}_{i_{t-1}} = 0, 0 \le \bar{x}_{i_t} = \frac{z - \sum_{k=t+1}^{|E|}\lambda_{i_k}}{\lambda_{i_t}} < 1, \bar{x}_{i_{t+1}} = \cdots = \bar{x}_{i_{|E|}} = 1\) corresponds to the maximal value of \(y\) given \(z\), denoted by \(\bar{y}(z)\), which is a concave piecewise linear function with \(|E|\) pieces. Let \((y^{ * }, z^{ * })\) be the optimal solution of the bi-variate optimization problem. Then \(x_i = \frac{y^{ * } - \underline{y}(z^{ * })}{\bar{y}(z^{ * }) - \underline{y}(z^{ * })}\bar{x}_i(z^{ * }) + \frac{\bar{y}(z^{ * }) - y^{ * }}{\bar{y}(z^{ * }) - \underline{y}(z^{ * })} \underline{x}_i(z^{ * }), i \in E\) is an optimal solution to the original problem. The detail of the algorithm is given in Algorithm \ref{alg:bivariate}. We also provide an illustration in the following example.
\begin{algorithm}[!htp]
	\caption{Algorithm for Searching for Optimal Fractional Assignment}
	\label{alg:bivariate}
	\begin{algorithmic}[1]
		\Require $\lambda_i, \mu_i, v_i, i = 1, \dots, n$ and $x_i, i \in L\cup U$, $\alpha$
		\Ensure $x_i, i \in E$
		
		\State Initialize $x_i = 0, i \in E$, start from the minimal $\mu_i$, $x_i \gets 1$, compute upper region of $\bar{y}(z)$ and corresponding assignment $\bar{x}_i(z), i\in E$, parameterized by $|E|+1$ pairs of $(z, \bar{y}(z))$
		\State Initialize $x_i = 0, i \in E$, start from the maximal $\mu_i$, $x_i \gets 1$, compute lower region of $\underline{y}(z)$ and corresponding assignment $\underline{x}_i(z), i\in E$, parameterized by $|E|+1$ pairs of $(z, \underline{y}(z))$
		\State Minimize \eqref{eq:bivariate} with $\underline{y}(z) \le y \le \bar{y}(z)$ with optimal solution $y^{ * }, z^{ * }$ \label{lst:line:bivariate}
		\State $x_i \gets \frac{y^{ * } - \underline{y}(z^{ * })}{\bar{y}(z^{ * }) - \underline{y}(z^{ * })}\bar{x}_i(z^{ * }) + \frac{\bar{y}(z^{ * }) - y^{ * }}{\bar{y}(z^{ * }) - \underline{y}(z^{ * })} \underline{x}_i(z^{ * }), i \in E$
	\end{algorithmic}
\end{algorithm}

\vspace{0.4 cm}
\begin{example}
	In this example, we show how to recover assignment decisions from the values of \(y\) and \(z\). Consider \(E = \{1, 2\}\) with \(\mu_1 > \mu_2\). \(\bar{y}(z)\) passes through \((1, 1), (2, 4)\) and \((10, 6)\), while \(\underline{y}(z)\) passes through \((1, 1), (6, 2)\) and \((10, 6)\). The optimal solution of the bivariate problem is \(y^{ * } = 3\) and \(z^{ * } = 4\). As is shown in Figure \ref{fig:bivariate}, \((4, 1.6)\) corresponds to the assignment \(x_1 = 0.6, x_2 = 0\) and \((4, 4.5)\) corresponds to the assignment \(x_1 = 0.25, x_2 = 1\). Thereby, the optimal assignment is \(x_1 = \frac{25}{58}\) and \(x_2 = \frac{14}{29}\). (In this case, since \(E\) only contains two points, one can also directly optimize \(x_1\) and \(x_2\). This example is mainly for showing the conversion from \(x\) space to \((y, z)\) space.) \(\hfill\square\)
\end{example}
\vspace{0.4 cm}

Finally, by enumerating all candidate configurations and solving corresponding optimization subproblems, we are able to provide an efficient algorithm to solve the assignment problem as stated in Algorithm \ref{alg:fixed-alpha}.

\begin{algorithm}[!htp]
	\caption{Algorithm for Solving QAP}
	\label{alg:fixed-alpha}
	\begin{algorithmic}[1]
		\Require $\lambda_i, \mu_i, v_i, i = 1, \dots, n$ and $\mu_1 > \mu_2 > \cdots > \mu_n, \alpha$
		\Ensure $\boldsymbol{x}^*\in\mathbb{R}^n$
		
		\State Initialize $f^*\gets+\infty$
		\For{each $(U, E, L) \in \mathcal{C}(\mathcal{I})$}
		\State $x_i \gets \begin{cases}
		0 &\text{if }i \in U\\
		1 &\text{if }i \in L
		\end{cases}$
		\State Solve for optimal $x_{i}, i \in E$ using Algorithm \ref{alg:bivariate}
		\State Compute objective function value $f$ ($f \gets +\infty$ if infeasible)
		\If {$f < f^*$}
		\State $f^* \gets f, \boldsymbol{x}^* \gets \boldsymbol{x}$
		\EndIf
		\State $x_i \gets \begin{cases}
		1 &\text{if }i \in U\\
		0 &\text{if }i \in L
		\end{cases}$
		\State Solve for optimal $x_{i}, i \in E$ using Algorithm \ref{alg:bivariate}
		\State Compute objective function value $f$ ($f \gets +\infty$ if infeasible)
		\If {$f < f^*$}
		\State $f^* \gets f, \boldsymbol{x}^* \gets \boldsymbol{x}$
		\EndIf
		\EndFor
	\end{algorithmic}
\end{algorithm}

\subsection{The Assignment Problem for Two-Queues under Exponential Service Distribution}
\label{ssec:qap}
In this subsection, we consider a special case of exponential service time distribution. Under this setting, we are able to simplify the formulation and reveal a special block-type structure.

We follow the notations in Section \ref{ssec:qap-g} where \(v_i = 2/\mu_i^2\) for \(i = 1, 2, \dots, n\) (corresponding to exponential distribution) and let
$$
\\[5pt]
f(\boldsymbol{x}, \alpha) :=
\begin{cases}
\frac{\sum\limits_{i=1}^n \frac{\lambda_i x_i}{\mu_i^2}}{\alpha^2 - \alpha\sum\limits_{i=1}^n \frac{\lambda_i x_i}{\mu_i}} \cdot \frac{\sum\limits_{i=1}^n x_i \lambda_i}{\sum\limits_{i=1}^n \lambda_i} + \frac{\sum\limits_{i=1}^n \frac{\lambda_i(1 - x_i)}{\mu_i^2}}{(1 - \alpha)^2 - (1 - \alpha)\sum\limits_{i=1}^n \frac{\lambda_i(1 - x_i)}{\mu_i}} \cdot \frac{\sum\limits_{i=1}^n (1 - x_i)\lambda_i}{\sum\limits_{i=1}^n \lambda_i}  & \text{ if } 0 < \alpha < 1 \\[20pt]
\left(\sum_{i=1}^n \frac{\lambda_i}{\mu_i^2}\right)/\left(1 - \sum_{i=1}^n \frac{\lambda_i}{\mu_i}\right)  & \text{ if } \alpha = 0, \boldsymbol{x} = \boldsymbol{0} \text{ or } \alpha = 1, \boldsymbol{x} = \boldsymbol{1}.\\[5pt]
\end{cases}
$$

We assume \(\mu_1 > \mu_2 > \cdots > \mu_n\), since if \(\mu_i = \mu_j\), we can redefine a type of customers with arrival rate \(\lambda_i + \lambda_j\) and service rate \(\mu_i\), to replace type \(i\) and \(j\) customers. In the following, we analyze the queue assignment problem \eqref{eq:qa-w} under this model.

Similar as in Section \ref{ssec:qap-g}, we can verify that although \(f(\boldsymbol{x}, \alpha)\) is convex in each \(x_i\), it is not jointly convex in \(\boldsymbol{x}\). Therefore, directly minimizing \(f(\boldsymbol{x}, \alpha)\) may be challenging. However, it turns out that the optimal solution to \eqref{eq:qa-w} has a special block-type structure. We describe it in the following theorem.

\vspace{0.3 cm}
\begin{theorem}
	\label{thm:fix-alpha-exp}
	Suppose \(\lambda_1, \dots, \lambda_n, \mu_1 > \cdots > \mu_n\) and \(\alpha \ge \frac{1}{2}\) are given. For any optimal solution \(\boldsymbol{x}^*\) to the QAP, there exist \(0 \le l < h \le n\) such that:
	
	\begin{itemize}
		\item When \(i > h\) or \(i < l\), \(x^*_i = 1\);
		\item When \(l < i < h\), \(x^*_i = 0\);
		\item When $ i \in \{l,h\} $, $ x^*_i \in [0,1] $. Moreover, the set \(\mathcal{M}^* = \left\{i: 0 < x^*_i < 1\right\}\) is either empty or singleton.
	\end{itemize}
\end{theorem}
\vspace{0.3 cm}

Theorem \ref{thm:fix-alpha-exp} reveals a few important structures of the optimal solution to \eqref{eq:qa-w}. First, the optimal assignment can have at most three continuous segments, ranked by the service rate. More specifically, customer types with sufficiently high and low service rates are assigned to one sub-queue, and those with medium service rates assigned to the other. This at-most-three-segment structure can sometimes be degenerate with fewer segments: Note that \(l\) could be \(0\), in which case the first segment does not exist. Similarly, \(h\) could be \(n\), in which case the third segment does not exist. When \(l = 0\) and \(h = n\), only the second segment exists.

In addition, the optimal assignment is \textit{almost} deterministic: At most one element (either the first or the last) of the segment can be fractional. If we refer to allowing fractional assignment of all customer types as \textit{full flexibility}, and allowing only \textit{one} customer type's assignment to be fractional as \textit{restricted flexibility}, then Theorem \ref{thm:fix-alpha-exp} shows that restricted flexibility can bring as much benefit as full flexibility. 

\begin{figure}[ht]
	\centering
	\includegraphics[height=8cm]{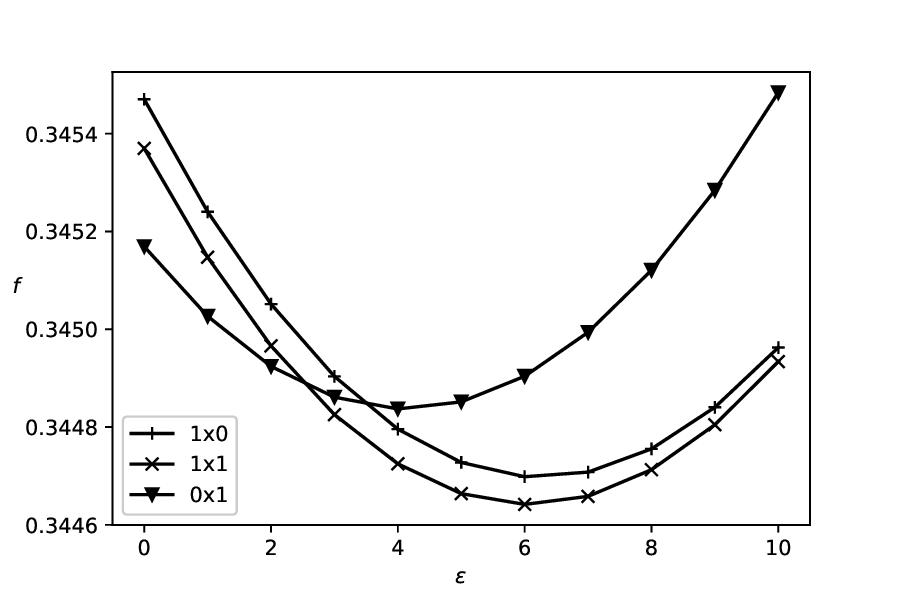}
	\caption{Objective function value under the assignment in the form of \((0, x_2, 1), (1, x_2, 0)\) and \((1, x_2, 1)\). \label{fig:example}}
\end{figure}

The optimal structure that we discover is a mild generalization of the common rule of thumb to partition customers into two continuous segments ranked by their service rates. In fact, this generalization is necessary: It is possible that the customers with the largest and the lowest service rates are assigned to the same queue, while the customers with the intermediate service rate are assigned to the other queue. We illustrate such a case in Example \ref{example:three segment} below.
\vspace{0.4 cm}
\begin{example}\label{example:three segment}
	Consider three types of customers with \(\lambda_1 = 0.4, \lambda_2 = 8, \lambda_3 = 0.2, \mu_1 = 16, \mu_2 = 12, \mu_3 = 10\). For partition \(\alpha = 0.8\), since \(1 - \alpha < \frac{\lambda_2}{\mu_2}\), we can see that \(x_2 = 0\) is not feasible. Therefore, by Theorem \ref{thm:fix-alpha-exp} we know that the optimal solution must be in the form of \((x_1, 1, 1)\), \((1, 1, x_3)\), \((0, x_2, 1)\), \((1, x_2, 0)\) or \((1, x_2, 1)\). According to the convexity of \(f(\boldsymbol{x}, \alpha)\) in \(\boldsymbol{x}\), \(\frac{\partial}{\partial x_1}f\left((x_1, 1, 1), 0.8\right) \ge \frac{\partial}{\partial x_1}f\left((x_1, 1, 1), 0.8\right)\biggr\rvert_{x_1=0} = 0.18 > 0\), thereby the minimum of \((x_1, 1, 1)\) type of solution is obtained at \(f\left((0, 1, 1), 0.8\right) = 0.61\).
	Similarly, \(\frac{\partial}{\partial x_3}f\left((1, 1, x_3), 0.8\right) \ge \frac{\partial}{\partial x_3}f\left((1, 1, x_3), 0.8\right)\biggr\rvert_{x_3=0} = 0.15 > 0\), thereby the minimum of \((1, 1, x_3)\) type of solution is obtained at \(f\left((1, 1, 0), 0.8\right) = 0.65\). Since \(f\left((1, 0.8, 1), 0.8\right) = 0.37\), we can see that the optimal assignment cannot be in the form of \((x_1, 1, 1)\) or \((1, 1, x_3)\). Therefore, the optimal solution must be of form \((0, x_2, 1)\), \((1, x_2, 0)\) or \((1, x_2, 1)\). Figure \ref{fig:example} shows how the values of \(f\left((0, 0.87 + 0.001\times\epsilon, 1), 0.8\right)\), \(f\left((1, 0.86 + 0.001\times\epsilon, 0), 0.8\right)\) and \(f\left((1, 0.83+0.001\times\epsilon, 1), 0.8\right)\) change with \(\epsilon \in [0, 10]\).
	By Lemma \ref{lemma:convex}, \(f\) is convex in \(x_2\) given \(x_1\) and \(x_3\). Therefore, the result in Figure \ref{fig:example} indicates that the optimal value is in the form of \((1, x_2, 1)\). By Algorithm \ref{alg:fixed-alpha-exp}, we can calculate that the optimal solution is \(\boldsymbol{x}^* = (1, 0.836, 1)\). Under this assignment, all of type \(1\) and type \(3\), and part of type \(2\) customers are assigned to the first queue, and the rest of type \(2\) customers are assigned to the second queue.\(\hfill\square\)
\end{example}

{It is worth noting that this optimal structure revealed in Theorem \ref{thm:fix-alpha-exp} is closely related to the nested size interval strategy proposed by \citet{feng2005optimal}. In particular, in \citet{feng2005optimal}, the authors consider a queue assignment problem where service time follows a general distribution and is known upon arrival. They show that under the optimal assignment, the service rate range of customers assigned to the lower capacity is nested within those assigned to the higher capacity queues. In other words, there exist $\mu_l \le \mu_h$ such that if the service time is among the range of $\frac{1}{\mu_h}$ and $\frac{1}{\mu_l}$ then customers are assigned to the lower capacity queue. Otherwise, if the service time is outside the range of $\frac{1}{\mu_h}$ and $\frac{1}{\mu_l}$, then customers are assigned to the higher capacity queue. In fact, there is indeed a connection between our analysis and theirs. Based on our results for the general distribution case (Theorem \ref{thm:fix-alpha}), the optimal partition is determined by the relative relation of $\mathbb{E}[S_i]$ and $\mathbb{E}[S_i^2]$ for different $i$'s. When we compare exponential service time with deterministic service time, the only difference is that $\mathbb{E}[S_i^2]$ will be two-fold for exponential service time compared to deterministic service time. Therefore, the nested size interval structure mentioned in Theorem \ref{thm:fix-alpha-exp} can be viewed as an extension of the results in \citet{feng2005optimal}. However, our results differ from theirs in three key aspects. First, we extend the results in \citet{feng2005optimal} and show that the optimal assignment has a restricted flexibility property with at most one type of customers splitting randomly to two queues. Second, as is shown in Theorem \ref{thm:fix-alpha-convex} in the appendix, we can prove that the nested size interval property always holds as long as the service distributions satisfy $\frac{\mathbb{E}[S_n^2] - \mathbb{E}[S_{n-1}^2]}{\mathbb{E}[S_{n}] - \mathbb{E}[S_{n-1}]} \ge \cdots \ge \frac{\mathbb{E}[S_2^2] - \mathbb{E}[S_1^2]}{\mathbb{E}[S_2] - \mathbb{E}[S_1]} \ge \frac{\mathbb{E}[S_1^2]}{\mathbb{E}[S_1]}$ for $\mathbb{E}[S_n] \ge \cdots \ge \mathbb{E}[S_1]$, which is a generalization of exponential distribution and cannot be reduced to the deterministic case. Third, as we will show in Section \ref{sec:sojourn}, the nested size interval strategy is also optimal in the sojourn time optimization setting, which extends the fixed load partitioning results of Proposition 10 in \citet{feng2005optimal} to general heterogeneous system. Moreover, in this paper, we also analyze the effect of joint partition-assignment optimization on the optimal assignment decision.}

Finally, we can leverage our insights to develop an efficient algorithm to solve the QAP. Specifically, by Proposition \ref{lemma:convex}, we know \(f(\boldsymbol{x}, \alpha)\) is convex in \(x_i\) for all \(i\). Thereby, an optimal solution to the QAP can be found by first enumerating all possible \(0\) or \(1\) elements in the solution (by Theorem \ref{thm:fix-alpha-exp} there are \(\mathcal{O}(n^2)\) of possible candidates), and then for each candidate solving a single-variable convex optimization for the fractional element (if it exists). The detail of the algorithm is given in Algorithm \ref{alg:fixed-alpha-exp}.

\begin{algorithm}[!htp]
	\caption{Algorithm for Solving QAP under Exponential Distribution}
	\label{alg:fixed-alpha-exp}
	\begin{algorithmic}[1]
		\Require $\alpha \ge \frac{1}{2}, \lambda_i, \mu_i, i = 1, \dots, n$ and $\mu_1 > \mu_2 > \cdots > \mu_n$
		\Ensure $\boldsymbol{x}^*\in\mathbb{R}^n$
		
		\State Initialize $f^*\gets+\infty$
		\For{$i = 0, 1, 2, \dots,  n$}
		\For{$j = i + 1, \dots, n + 1$}
		\State $x_k \gets \begin{cases}
		0 &\text{if }i < k < j\\
		1 &\text{if }k > j \text{ or } k < i\\
		\end{cases}$
		\State $x_i \gets 0$, solve for optimal $x_j$
		\State Compute objective function value $f$ ($f \gets +\infty$ if infeasible)
		\If {$f < f^*$}
		\State $f^* \gets f, \boldsymbol{x}^* \gets \boldsymbol{x}$
		\EndIf
		\State For $x_j = 0, 1$, solve for optimal $x_i$, repeat line 6 to line 9
		\EndFor
		\EndFor
	\end{algorithmic}
\end{algorithm}

\subsection{The Joint Queueing Partition and Assignment Problem for Two-Queues}
In this subsection, we consider the (joint) queueing partition (and assignment) problem
\begin{eqnarray}\label{eq:qp-w}
\textup{(QPP)}\quad\inf_{\boldsymbol{x}, \alpha} &&\quad f(\boldsymbol{x}, \alpha) \nonumber \\
\text{s.t.} &&\quad (\boldsymbol{x}, \alpha) \in \mathcal{F}; \nonumber \\
&&\quad 0\le\alpha\le1; \quad x_i\in [0, 1], \; i = 1, \dots, n.
\end{eqnarray}

We can incorporate the optimization of \(\alpha\) in line \ref{lst:line:bivariate} of Algorithm \ref{alg:bivariate} to solve a tri-variate optimization problem. In this way, we are able to solve the QPP problem \eqref{eq:qp-w} and the details are given in Algorithm \ref{alg:opt-alpha}.

\begin{algorithm}[!htp]
	\caption{Algorithm for Solving QPP}
	\label{alg:opt-alpha}
	\begin{algorithmic}[1]
		\Require $\lambda_i, \mu_i, v_i, i = 1, \dots, n$ and $\mu_1 > \mu_2 > \cdots > \mu_n$
		\Ensure $\boldsymbol{x}^*\in\mathbb{R}^n, \alpha^{ * } \in \mathbb{R}$
		
		\State Initialize $f^*\gets+\infty$
		\For{each $(U, E, L) \in \mathcal{C}(\mathcal{I})$}
		\State $x_i \gets \begin{cases}
		0 &\text{if }i \in U\\
		1 &\text{if }i \in L
		\end{cases}$
		\State Solve for optimal $x_{i}, i \in E$ and $\alpha$
		\State Compute objective function value $f$ ($f \gets +\infty$ if infeasible)
		\If {$f < f^*$}
		\State $f^* \gets f, (\boldsymbol{x}^*, \alpha^{ * }) \gets (\boldsymbol{x}, \alpha)$
		\EndIf
		\State $x_i \gets \begin{cases}
		1 &\text{if }i \in U\\
		0 &\text{if }i \in L
		\end{cases}$
		\State Solve for optimal $x_{i}, i \in E$
		\State Compute objective function value $f$ ($f \gets +\infty$ if infeasible)
		\If {$f < f^*$}
		\State $f^* \gets f, (\boldsymbol{x}^*, \alpha^{ * }) \gets (\boldsymbol{x}, \alpha)$
		\EndIf
		\EndFor
	\end{algorithmic}
\end{algorithm}
 
Recall that the geometric property introduced in Theorem \ref{thm:fix-alpha} tells us that it is not necessarily optimal to partition customers into two continuous segments ranked by their service rates. Nevertheless, its prevalent use in practice motivates us to analyze the QPP problem further. Indeed, we find that partitioning customers into two segments is optimal in a variety of cases. In particular, we consider the following assumption on the service time.

\begin{assumption}
	\label{as:reg-2seg}
	For a set of service time \(\mathcal{I} = \{S_i\}_{i=1}^n\), we have \(\frac{\mathbb{E}[S_i^2]}{\mathbb{E}[S_i]} \ge \frac{\mathbb{E}[S_j^2]}{\mathbb{E}[S_j]} \text{ for any } \mathbb{E}[S_i] \ge \mathbb{E}[S_j]\).
\end{assumption}

Intuitively, Assumption \ref{as:reg-2seg} regulates the service time distributions such that service time distribution with larger mean also has larger adjusted variation (second moment divided by the mean). Under this assumption, customers with low service rate also suffer from high variability of service time. Consequently, there is no incentive to pool customers with high or low service rates while keeping those with moderate service rates in a dedicate queue. Note that this assumption is satisfied in many settings including exponential distributions, half-normal distributions and Rayleigh distributions. We show that under this assumption, somewhat surprisingly, the optimal assignment in this problem assigns each type of customers to a unique queue, which we refer to as a deterministic assignment policy. Furthermore, the optimal partition is to separate customers with large service rates from those with small service rates.

By analyzing the properties of the optimal solution to the QPP, we have the following result.

\vskip 0.3 cm
\begin{theorem}
	\label{thm:opt-alpha}
	Given \(\lambda_i, i = 1, 2, \dots, n\) and \(\mathcal{I}\) with \(\mu_1 > \cdots > \mu_n, v_1, \dots, v_n\) satisfying Assumption \ref{as:reg-2seg}. For any optimal solution \((\boldsymbol{x}^*, \alpha^*)\) to the queue partition problem, there exists \(i^*\) such that either 
	\begin{align*}
	x^*_i = \begin{cases}
	1 &\text{if }i \le i^*\\
	0 &\text{if }i > i^* 
	\end{cases} \qquad \text{ or } \qquad  x^*_i = \begin{cases}
	0 &\text{if }i \le i^*\\
	1 &\text{if }i > i^* 
	\end{cases}.
	\end{align*}
\end{theorem}
\vskip 0.3 cm

Theorem \ref{thm:opt-alpha} has several implications. First, it indicates that the optimal assignment in the QPP assigns each type of customers to a unique queue. Recall that we refer to allowing assignment of all customer types to be fractional as full flexibility, and refer to allowing no customer type's assignment to be fractional as no flexibility. Theorem \ref{thm:opt-alpha} indicates that the additional flexibility of allowing fractional assignment does not have extra value in this queue partition problem when Assumption \ref{as:reg-2seg} holds. This is different from the queue assignment problem in which allowing fractional assignment may further reduce the expected waiting time compared to the assignment where each type of customers can only be assigned to a unique queue. In particular, if we restrict our attention to the deterministic assignment policy, i.e. \(x_i \in \{0, 1\}\) for \(i = 1, \dots, n\), then finding the optimal assignment is computationally harder than \eqref{eq:qa-w} (we show in Theorem \ref{thm:np-hard} in the appendix that problem \eqref{eq:qa-w} with integer constraints is NP-hard). In comparison, the optimal fractional assignment and the optimal deterministic assignment coincide in the queue partition problem when Assumption \ref{as:reg-2seg} holds. Second, the optimal partition in the QPP has \textit{at most two} continuous segments, in which all customers with high service rates are assigned to one queue, and all the rest of the customers are assigned to the other queue. Note that it is possible that one segment is empty (equivalently, \(i^*=0\) or \(n\)), in which case the optimal partition is to choose \(\alpha = 1\) and assign all customers to one queue, and there is no benefit of splitting the queue.

In particular, since the exponential distribution family satisfies Assumption \ref{as:reg-2seg}, it is readily shown the following results.

\vspace{0.4 cm}
{\begin{corollary}
	\label{cor:opt-alpha-exp}
	Given \(\lambda_i, i = 1, 2, \dots, n\) and \(\mu_1 > \mu_2 > \cdots > \mu_n\). For any optimal solution \((\boldsymbol{x}^*, \alpha^*)\) to the QPP under exponential service time distribution, there exists \(i^*\) such that either
	\begin{align*}
	x^*_i = \begin{cases}
	1 &\text{if }i \le i^*\\
	0 &\text{if }i > i^* 
	\end{cases} \qquad \text{ or } \qquad x^*_i = \begin{cases}
	0 &\text{if }i \le i^*\\
	1 &\text{if }i > i^* 
	\end{cases}.
	\end{align*}
\end{corollary}}
\vspace{0.4 cm}

With Theorem \ref{thm:opt-alpha}, we can design a more efficient algorithm
to solve the queue partition problem. Note that given \(\boldsymbol{x}\), \(f(\boldsymbol{x}, \alpha)\) is convex in \(\alpha\). Therefore, the algorithm enumerates over all candidates of \(\boldsymbol{x}\) (at most \(n\) of them) and for each \(\boldsymbol{x}\) solves a single-variable convex optimization for the optimal \(\alpha\). The detailed algorithm is given in Algorithm \ref{alg:opt-alpha-spec}.
\begin{algorithm}[!htp]
	\caption{Algorithm for Solving QPP under Assumption \ref{as:reg-2seg}}
	\label{alg:opt-alpha-spec}
	\begin{algorithmic}[1]
		\Require $\lambda_i, i = 1, \dots, n$, $\mu_1 > \mu_2 > \cdots > \mu_n$ and $v_i, i = 1, \dots, n$
		\Ensure $\boldsymbol{x}^*\in\mathbb{R}^n, \alpha^*\in\mathbb{R}$
		
		\State Initialize $f^*\gets+\infty$
		\For{$i = 1, 2, \dots,  n + 1$}
		\State $x_k = \begin{cases}
		0 &\text{if }k < i\\
		1 &\text{if }k \ge i\\
		\end{cases}$
		\State $f \gets \min_\alpha f(\boldsymbol{x}, \alpha)$, $\bar{\alpha} \gets \arg\min_\alpha f(\boldsymbol{x}, \alpha)$
		\If {$f < f^*$}
		\State $f^* \gets f, \boldsymbol{x}^* \gets \boldsymbol{x}, \alpha^* \gets \bar{\alpha}$
		\EndIf
		\EndFor
	\end{algorithmic}
\end{algorithm}

\subsection{From Two-Queues to Multiple Queues}
In this subsection, we treat \(k\) as a given constant and analyze the problems introduced in Section \ref{sec:model}. Before proceeding, we first generalize the definition of the configuration set. Let \(\textbf{int}\,P\) and \(\textbf{bd}\,P\) denote the interior and boundary of polygon \(P\) respectively. Let \(U = \{(\mathbb{E}[S_i], \mathbb{E}[S_i^2]): i = 1, 2, \dots, n\}\). Given \(k \ge 2\), for any \(k\) convex polygons \(\mathcal{P} = (P_1, P_2, \dots, P_k)\) such that
\begin{equation}\label{eq:polygon}
\begin{aligned}
U \subseteq \bigcup_{j=1}^k P_j, \quad
\textbf{int}\,P_{j_1} \bigcap \textbf{int}\,P_{j_2} = \emptyset, \,\, \forall j_1 \neq j_2,
\text{ and if } P_{j_1} \bigcap P_{j_2} \neq \emptyset, \text{ then } \left|P_{j_1} \bigcap P_{j_2} \bigcap U\right| \ge 2,
\end{aligned}
\end{equation}
we define sets
\begin{equation} \label{eq:config-k}
\begin{aligned}
T_{\mathcal{P},j} = & \left\{i \in \{1, 2, \dots, n\}: (\mathbb{E}[S_i], \mathbb{E}[S_i^2]) \in \textbf{int}\, P_j \right\},\qquad \forall j = 1, 2, \dots, k,\\
E_{\mathcal{P},j} = & \left\{i \in \{1, 2, \dots, n\}: (\mathbb{E}[S_i], \mathbb{E}[S_i^2]) \in \textbf{bd}\, P_j \right\},\qquad \forall j = 1, 2, \dots, k.
\end{aligned}
\end{equation}

\begin{figure}[ht]
	\centering
	\includegraphics[height=9cm]{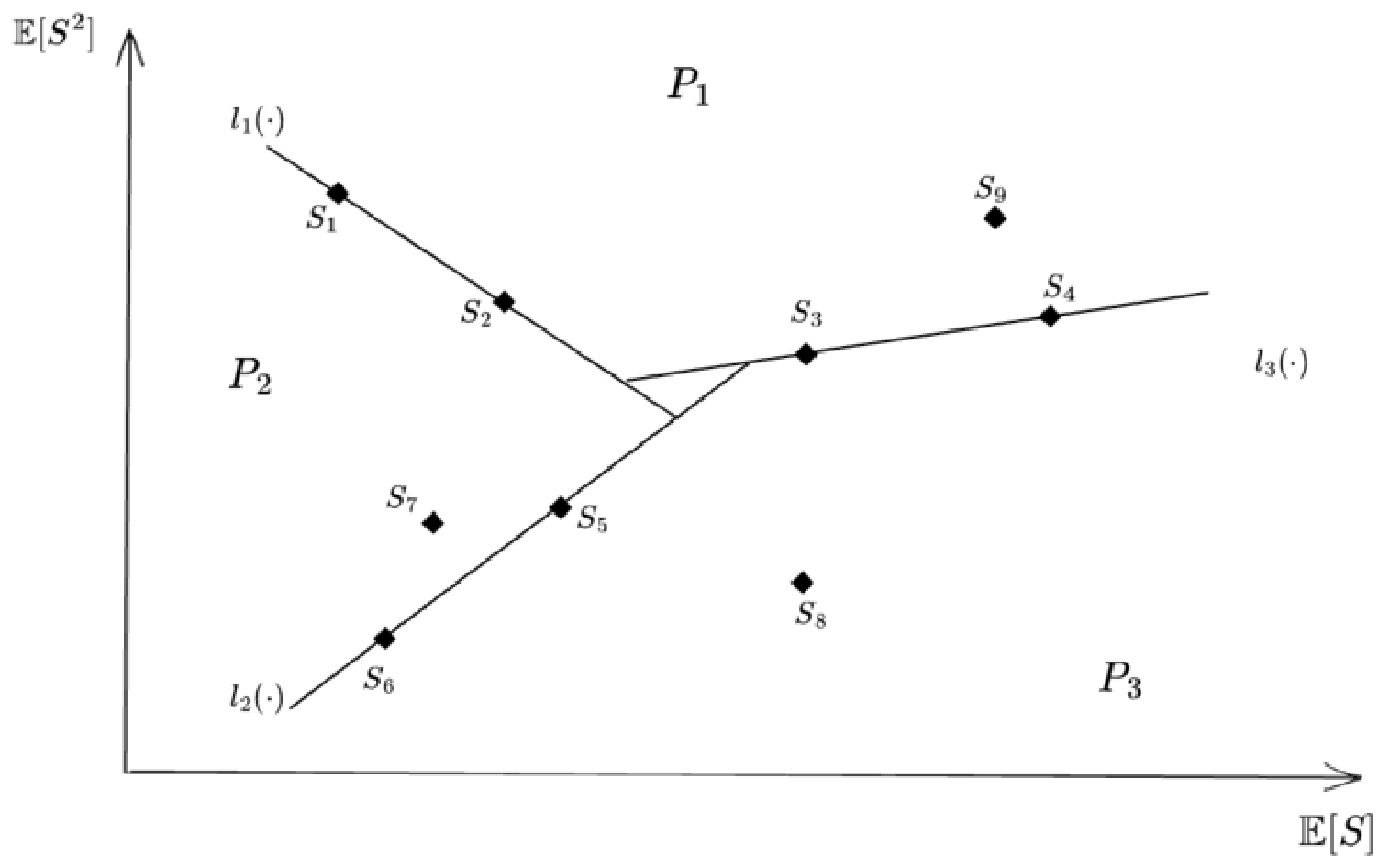}
	\caption{There are 9 types of customers and 3 queues. The induced sets are \(T_{\mathcal{P}, 1} = \{9\}, E_{\mathcal{P}, 1} = \{1, 2, 3, 4\}, T_{\mathcal{P}, 2} = \{7\}, E_{\mathcal{P}, 2} = \{1, 2, 5, 6\}, T_{\mathcal{P}, 3} = \{8\}, E_{\mathcal{P}, 3} = \{3, 4, 5, 6\}\). \label{fig:ppartitionk3}}
\end{figure}

To illustrate, consider a \(k\)-QAP problem with \(n = 9\) and \(k = 3\). As is shown in Figure \ref{fig:ppartitionk3}, \(P_1\) is the upper-right polygon, \(P_2\) is the left polygon and \(P_3\) is the lower-right polygon\footnote{{As is shown in Figure \ref{fig:ppartitionk3}, $P_1, P_2, P_3$ do not necessarily partition the two-dimensional space, since there exists an empty triangle in the middle. However, we must ensure that there is no $(\mathbb{E}[S_i], \mathbb{E}[S_i^2])$ in the interior of the middle triangle.}}. The induced sets are \(T_{\mathcal{P}, 1} = \{9\}, E_{\mathcal{P}, 1} = \{1, 2, 3, 4\}, T_{\mathcal{P}, 2} = \{7\}, E_{\mathcal{P}, 2} = \{1, 2, 5, 6\}, T_{\mathcal{P}, 3} = \{8\}, E_{\mathcal{P}, 3} = \{3, 4, 5, 6\}\).

Intuitively, polygons in \(\mathcal{P}\) classify service time distributions into \(k\) clusters according to the mean and second-order moment. In the graph of \((\mathbb{E}[S], \mathbb{E}[S^2])\), each cluster is separated from the other \(k-1\) clusters by at most \(k - 1\) different lines. In particular, when \(k = 2\), it reduces to the classification by one line.

In the following, we define the set of partitions satisfying \eqref{eq:polygon} as \(\mathscr{P}\). We define the \(k\)th order configuration set of \(\mathcal{I}\)
\begin{equation}\label{eq:def-config-k}
\mathcal{C}(\mathcal{I}, k) = \left\{(T_1, \dots, T_k, E_1, \dots, E_k): (T_1, \dots, T_k, E_1, \dots, E_k) \text{ satisfies \eqref{eq:config-k}} \text{ for some } \mathcal{P} \in \mathscr{P} \right\}.
\end{equation}

Intuitively, \(\mathcal{C}(\mathcal{I}, k)\) defines a collection of classifications of customer types based on the mean and second-order moment of the service time distribution. In the following, we characterize the size of the set \(\mathcal{C}(\mathcal{I}, k)\). First, from the separating hyperplane theorem \citep[pg.~46]{BoydV04} we know that for each \((T_1, \dots, T_k, E_1, \dots, E_k) \in \mathcal{C}(\mathcal{I}, k)\), there exist a set \(\mathcal{L}\) of \(\frac{(k - 1)k}{2}\) linear functions \(\{l_1(\cdot), l_2(\cdot), \dots, l_{\frac{(k-1)k}{2}}(\cdot)\}\) and \(k\) convex polygons \(\mathcal{P} = (P_1, \dots, P_k)\) such that \(\left|\{(u, v): v = l_i(u)\} \bigcap U\right|\ge2\) for \(i = 1, \dots, \frac{k(k-1)}{2}\), \(T_1, \dots, T_k, E_1, \dots, E_k\) satisfy \eqref{eq:config-k}, and \(P_{j_1} \bigcap P_{j_2} \subseteq \displaystyle\bigcup_{i=1}^{\frac{(k-1)k}{2}}\left\{(u, v): v = l_i(u)\right\}\) for any \(j_1 \neq j_2\). In particular, each linear function is the separating line between \(\textbf{int}\,(T_{j_1}\cup E_{j_1})\) and \(\textbf{int}\,(T_{j_2}\cup E_{j_2})\) for some \(j_1 \neq j_2\). It indicates that for each \(P_j\), there exists a set \(\mathcal{M}\) of indicators \(e_1, e_2, \dots, e_{\frac{k(k-1)}{2}} \in \{-1 , 0, 1\}\) such that \(P_j = \{(u, v): e_k \cdot (v - l_{j_k}(u)) \ge 0, k = 1, 2, \dots, \frac{k(k-1)}{2}\}\). In other words, each element of \(\mathcal{C}(\mathcal{I}, k)\) can be determined by \(\mathcal{L}\) and \(\mathcal{M}\). To illustrate, consider the example in Figure \ref{fig:ppartitionk3}. We have \(P_1 = \{(u, v): v \ge l_1(u), v \ge l_3(u)\}\), \(P_2 = \{(u, v): v \le l_1(u), v \ge l_2(u)\}\), and \(P_3 = \{(u, v): v \le l_2(u), v \le l_3(u)\}\). Note that the number of possible choices of each element in \(\mathcal{L}\) is bounded by \(\binom{n}{2}\). By enumerating all possible \(\mathcal{L}\) and \(\mathcal{M}\), we have  \(|\mathcal{C}(\mathcal{I}, k)| = \mathcal{O}(n^{k(k-1)})\).

Based on \(\mathcal{C}(\mathcal{I}, k)\), we are ready to state the optimal structure of the multi-queue assignment problem in the following theorem.

\vspace{0.4 cm}
\begin{theorem}
	\label{thm:multi-fixed-alpha}
	For any given \(k\), suppose \(\lambda_1, \dots, \lambda_n, \mu_1 > \cdots > \mu_n, v_1, \dots, v_n\) and \(\alpha\) are given. Then for any optimal solution \(X\) to the \(k\)-QAP problem, there exist a permutation \(\boldsymbol{\pi} = (\pi_1, \dots, \pi_k)\) of \(\{1, 2, \dots, k\}\) and \((T_1, \dots, T_k, E_1, \dots, E_{k}) \in \mathcal{C}(\mathcal{I}, k)\) such that
	$$\begin{cases}
	X_{i \pi_j} = 1 & i \in T_j\\
	X_{i \pi_j} \in [0, 1] & i \in E_{j}\\
	X_{i \pi_j} = 0 & \text{otherwise}.
	\end{cases}$$
\end{theorem}
\vspace{0.4 cm}

Theorem \ref{thm:multi-fixed-alpha} reveals a special property of the optimal solution. Specifically, the optimal partition of \(k\) queues is connected to the partition of points in a two-dimensional space. In the graph of \((\mathbb{E}[S], \mathbb{E}[S^2])\), \(k\) convex polygons correspond to \(k\) queues such that if \((\mathbb{E}[S_i], \mathbb{E}[S_i^2])\) is in the interior of \(j\)th polygon then customers of type \(i\) will be assigned to the \(j\)th queue. To illustrate, consider the example in Figure \ref{fig:ppartitionk3}. One candidate of the optimal assignment induced by functions \(l_1(\cdot)\) and \(l_2(\cdot)\) is
$$
X = \begin{pmatrix}
-& -& 0\\
-& -& 0\\
-& 0& -\\
-& 0& -\\
0& -& -\\
0& -& -\\
0& 1& 0\\
0& 0& 1\\
1& 0& 0
\end{pmatrix},
$$
where ``\(-\)'' denotes a value in \([0, 1]\). Theorem \ref{thm:multi-fixed-alpha} states that one of all these candidates must be optimal for the corresponding \(k\)-QAP problem. In particular, when \(k = 2\), Theorem \ref{thm:multi-fixed-alpha} reduces to the structure in Theorem \ref{thm:fix-alpha}.

With this property, we can develop an algorithm to solve the \(k\)-QAP problem. We first enumerate all possible \(0\)-\(1\) elements based on \(\mathcal{C}(\mathcal{I}, k)\) and then solve the optimization problem for the fractional elements (if exist). The enumeration contains \(\mathcal{O}(n^{(k-1)k})\) possibilities. The detailed algorithm is given in Algorithm \ref{alg:multi-fixed-alpha}.

\begin{algorithm}[!htp]
	\caption{Algorithm for Solving $k$-QAP}
	\label{alg:multi-fixed-alpha}
	\begin{algorithmic}[1]
		\Require $\alpha_j, j = 1, \dots, k$, $\lambda_i, \mu_i, i = 1, \dots, n$ and $\mu_1 > \mu_2 > \cdots > \mu_n$
		\Ensure $X^*\in\mathbb{R}^n$
		
		\State Initialize $f^*\gets+\infty$
		\State Enumerate all possible integer elements derived from $\mathcal{C}(\mathcal{I}, k)$. Let $\mathcal{X}$ be a set of tuple $(X, M)$ where $X$ is the candidate assignment, and $M$ is the fractional indices set
		\For{$(X, M) \in \mathcal{X}$}
		\State Solve for optimal $X_{ij}$ for all $(i, j) \in M$ \label{lst:line:multi-variate}
		\State Compute objective function value $f$ ($f \gets +\infty$ if infeasible)
		\If {$f < f^*$}
		\State $f^* \gets f, X^* \gets X$
		\EndIf
		\EndFor
	\end{algorithmic}
\end{algorithm}

Next, we show that in the special case of exponential service time distribution, a special block-type structure can be utilized to develop a more efficient algorithm. By letting \(v_i = 2/\mu_i^2\) for \(i = 1, 2, \dots, n\), we state the result in the following theorem.

\begin{theorem}
	\label{thm:multi-exp}
	Let $ k $, $\lambda_1, \dots, \lambda_n$ and $ \mu_1 > \cdots > \mu_n$ be given. Suppose \(\alpha_1, \dots, \alpha_k > 0\) are given. Any optimal solution \(X^*\) to the \(k\)-QAP consists of at most $ (2k-1) $ ``blocks" of non-zero entries. More precisely, there exist \(0 = l_0 \le l_1 \le l_2 \le \cdots \le l_{2k-2} \le l_{2k-1} = n + 1\) and \(o_1, o_2, \dots, o_{2k-1} \in \{1, 2, \dots, k\}\) such that:
	\begin{itemize}
		\item  For every \(1 \le h \le 2k - 1\), we have \(X^*_{i o_h} = 1\) for all \(l_{h-1} < i < l_{h}\);
		\item For every $ j \in \{1, 2, \ldots, k\} $, we have 
		\begin{itemize}
			\item \(X^*_{i j} = 0\) for all $ i \in \{1,2,\ldots,n\} \setminus \cup_{h: o_h = j} \{i: l_{h-1} \leq j \leq l_h\} $;
			\item \(X^*_{i j} \in [0,1]\) for all $ i \in \cup_{h: o_h = j} \{l_{h-1}, l_h\} $.
			% \item The set \(\mathcal{M}_j^* = \left\{i: 0 < X^*_{ij} < 1\right\}\) is either empty or singleton.
		\end{itemize} 
	\end{itemize}
\end{theorem}

\begin{figure}[ht]
	\centering
	\includegraphics[height=4cm]{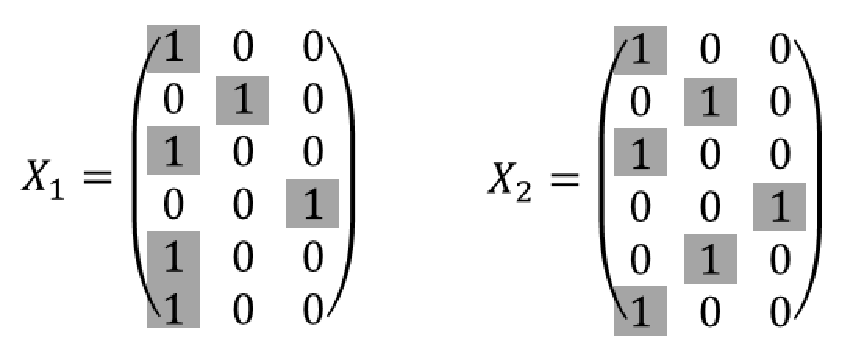}
	\caption{There are 5 blocks in \(X_1\) and 6 blocks in \(X_2\), where blocks are marked as the shaded area. \label{fig:block}}
\end{figure}

We now provide some explanations about Theorem \ref{thm:multi-exp}, which generalizes Theorem \ref{thm:fix-alpha-exp} by letting $ k = 2 $. The first part of Theorem \ref{thm:multi-exp} reveals a special column-wise property of the optimal \(k\)-QAP solution. Specifically, let \(X^* \in \mathbb{R}^{n\times k}\) be an optimal solution to the \(k\)-QAP. For each column, if we call consecutive $ 1 $-valued entries as a block, then \(X^*\) contains at most \(2k - 1\) blocks. To illustrate, consider a \(k\)-QAP problem with \(n = 5\) and \(k = 3\). As is shown in Figure \ref{fig:block}, \(X_1\) is a valid candidate for the optimal solution because the number of blocks is \(5 \le 2k - 1 = 5\), while \(X_2\) is not because the number of blocks is \(6 > 2k - 1 = 5\). 

With this property, we can now develop a more efficient algorithm to solve the \(k\)-QAP problem. Specifically, an optimal solution to the \(k\)-QAP problem can be obtained by first enumerating all possible \(0\)-\(1\) elements and then solving for the fractional elements (if they exist). By Theorem \ref{thm:multi-exp}, the enumeration can be obtained by first enumerating values of \(0 \le l_1 \le \cdots \le l_{2k-2} \le n + 1\) with \(\mathcal{O}(n^{2k-2})\) possibilities, then enumerating values of \(o_1, \dots, o_{2k-1}\) with \(\mathcal{O}(1)\) possibilities (recall \(k\) is treated as a constant in our analysis), and finally deciding the values of \(X_{l_1 o_1}, X_{l_{2k-2} o_{2k-1}}\) and \(X_{l_{i} o_{i+1}}, X_{l_{i+1} o_{i+1}}\) for \(i = 1, 2, \dots, 2k-3\) (this requires a grid search on these values, the complexity of such a step depends on $k$ and is independent of $n$). Thereby, the enumeration complexity is \(\mathcal{O}(n^{2k-2})\). The detailed algorithm is given in Algorithm \ref{alg:multi-fixed-alpha-exp}.

\begin{algorithm}[!htp]
	\caption{Algorithm for Solving $k$-QAP under Exponential Distribution}
	\label{alg:multi-fixed-alpha-exp}
	\begin{algorithmic}[1]
		\Require $\alpha_1 \ge \alpha_2 \ge \cdots \ge \alpha_k, \lambda_i, \mu_i, i = 1, \dots, n$ and $\mu_1 > \mu_2 > \cdots > \mu_n$
		\Ensure $X^*\in\mathbb{R}^n$
		
		\State Initialize $f^*\gets+\infty$
		\State Enumerate all possible integer elements. Let $\mathcal{X}$ be a set of tuple $(X, M)$ where $X$ is the candidate assignment, and $M = \{(l_1, o_1), (l_1, o_2), (l_2, o_2), \dots, (l_{2k-3}, o_{2k-2}), (l_{2k-2}, o_{2k-2}), (l_{2k-2}, o_{2k-1}) \}$ is the fractional indices set
		\For{$(X, M) \in \mathcal{X}$}
		\State Solve for optimal $X_{ij}$ for all $(i, j) \in M$
		\State Compute objective function value $f$ ($f \gets +\infty$ if infeasible)
		\If {$f < f^*$}
		\State $f^* \gets f, X^* \gets X$
		\EndIf
		\EndFor
	\end{algorithmic}
\end{algorithm}

Lastly, we provide a multi-queue analog for the queue partition problem, which is analogous to Theorem \ref{thm:opt-alpha}.

\begin{theorem}
	\label{thm:multi-opt-alpha}
	Given \(\lambda_i, i = 1, 2, \dots, n\) and \(\mathcal{I}\) with \(\mu_1 > \cdots > \mu_n, v_1, \dots, v_n\) satisfying Assumption \ref{as:reg-2seg}. If the feasible set of the \(k\)-QPP is non-empty, then there exists an optimal solution \((X^*, \boldsymbol{\alpha}^*)\) to the \(k\)-QPP and \(1 = i_0 \le i_1 \le i_2 \le \cdots \le i_{k-1} \le i_k = n + 1\) such that 
	\begin{align*}
	X^*_{i j} = \begin{cases}
	1 &\text{if }i_{j-1} \le i < i_{j} \\
	0 &\text{otherwise},
	\end{cases}
	\end{align*}
	for \(j = 1, 2, \dots, k\).
\end{theorem}

Theorem \ref{thm:multi-opt-alpha} shows that under Assumption \ref{as:reg-2seg}, the optimal partition must be to partition the customer types according to their service rates. More specifically, it must be optimal to assign customers with consecutive service rates to the same queue. Such a result can give rise to an efficient algorithm for the \(k\)-QPP. Particularly, note that given \(X\), we can verify that \(f(X, \boldsymbol{\alpha})\) is jointly convex in \(\boldsymbol{\alpha}\). Therefore, we can enumerate all candidates of \(X\)
(by Theorem \ref{thm:multi-opt-alpha}, there are at most \(\mathcal{O}(n^{k-1})\) of them) and for each \(X\) solve a convex optimization for the optimal \(\boldsymbol{\alpha}\). The detailed algorithm is given in Algorithm \ref{alg:multi-opt-alpha}.

\begin{algorithm}[!htp]
	\caption{Algorithm for Solving $k$-QPP under Assumption \ref{as:reg-2seg}}
	\label{alg:multi-opt-alpha}
	\begin{algorithmic}[1]
		\Require $\lambda_i, i = 1, \dots, n$ and $\mu_1 > \mu_2 > \cdots > \mu_n$
		\Ensure $X^*\in\mathbb{R}^{n\times k}, \boldsymbol{\alpha}^*\in\mathbb{R}^k$
		
		\State Initialize $f^*\gets+\infty$
		\State $i_0\gets 0, i_{k+1} \gets n+1$
		\For{$1 \le i_1 \le i_2 \le \cdots \le i_k \le n + 1$}
		\State $X_{ij} = \begin{cases}
		1 &\text{if }i_{j-1} \le i < i_{j} \\
		0 &\text{Otherwise}
		\end{cases}$
		\State $f \gets \min_{\boldsymbol{\alpha}} f(X, \boldsymbol{\alpha})$, $\boldsymbol{\bar{\alpha}} \gets \arg\min_{\boldsymbol{\alpha}} f(X, \boldsymbol{\alpha})$
		\If {$f < f^*$}
		\State $f^* \gets f, X^* \gets X, \boldsymbol{\alpha}^* \gets \boldsymbol{\bar{\alpha}}$
		\EndIf
		\EndFor
	\end{algorithmic}
\end{algorithm}

\subsection{Choice of Customer Types}
{
In previous subsections, we assume that the customer types are konwn in advance. However, in reality, customer types are often unknown and might need to be determined by predictive or clustering algorithms. In this subsection, we show how to determine the number of customer types and illustrate the effect of choosing a coarser or finer type grouping for the customers.

In practical applications, customer segmentation is typically derived from various customer features. These features enable machine leanring algorithms to cluster customers in the finest level. With a sufficient sample size of historical data, we can use this customer segmentation to develop queueing partition rules. However, in cases where historical data are limited, a coarser customer grouping can be beneficial. To see this, we consider a queueing system where customer arrival follows a Poisson process with rate 1. Each customer's service time follows an exponential distribution, where service rate is uniformly randomly chosen from $1, 3, 10, 11, 13$. In other words, there are $n=5$ types of customers in the system. Since the focus is to determine the optimal number of customer groups, for simplicity, we assume each customer's type is known upon arrival based on customer features from an exogenous machine learning algorithm. When corresponding service time distribution is known, we can formulate the queueing partition problem with $n=5$ customer types and find the optimal partition and assignment decision. However, the corresponding service time distribution is usually unknown. Instead, we assume there are 50 customer service time and customer type data, based on which we can infer the moment information of the service time distributions. Therefore, merging different customer types can be beneficial in that the pooling of more data may lead to a better moment estimation. Specifically, we start with 2 customer groups, iteratively computing the corresponding optimal waiting time. In particular, when we try to group customers into $m$ groups, we use the K-Nearest Neighbors algorithm based on their first and second order sample moment. For each group, we then estimate the first two moments of service time distribution, and use this value to compute the optimal partition $\alpha$ and assignment $\boldsymbol{x}$ for the $m$ customer type system. Since the first two moments are based on sample estimation, the corresponding waiting time serves as an estimation. The true expected waiting time is also computed as comparison. We repeat this procedure 1000 times. By averaging across these 1000 iterations, we illustrate how expected waiting time changes as $m$ increases in Figure \ref{fig:type-choice}.
\begin{figure}[ht]
	\centering
	\includegraphics[width=0.8\linewidth]{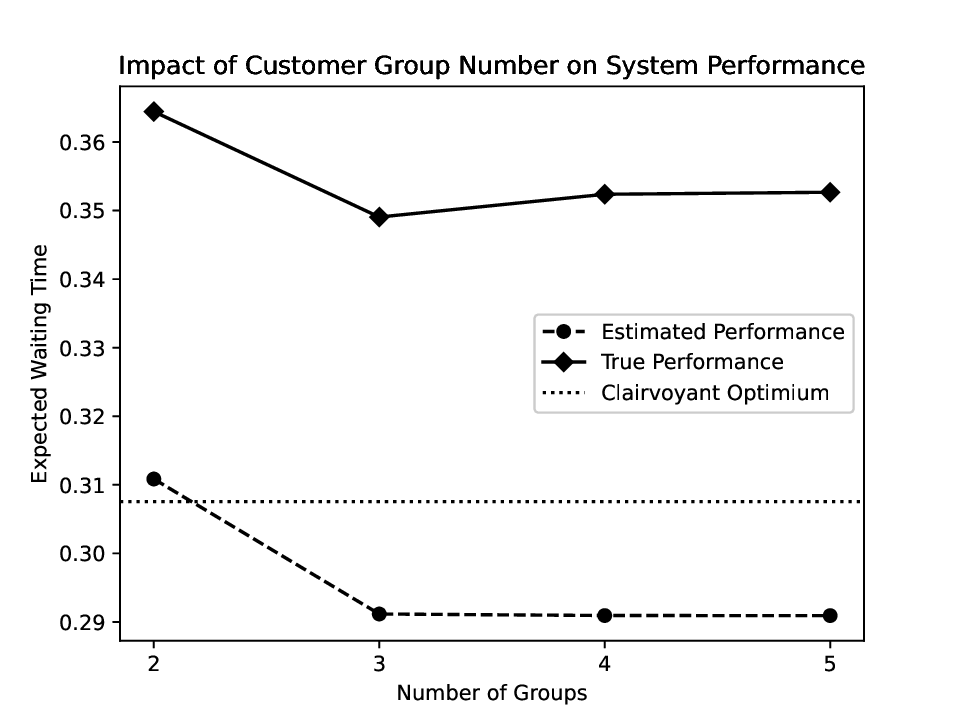}
	\caption{Clustering customers with Exponential service time into $2 \le m \le 5$ groups. \label{fig:type-choice}}
\end{figure}

As is shown in the figure, the true expected waiting time first decreases as the number of customer groups increases and then increases, and the best customer group number in this case is 3. We also notice that Figure \ref{fig:type-choice} depicts a diminishing marginal benefit in the estimated performance correlated with progressively refined customer segmentation. When the customer group number is small, a finer grouping for customer types may allow finer partition of queues, which has the potential of decreasing the expected waiting time. In contrast, when there are already a sufficient number of customer groups, more customer groups will exacerbate the severity of overfitting and model misspecification. Consequently, it is essential to strike a balance between the depth of segmentation and considerations of accuracy. In practice, the optimal number of customer groups can be determined by, e.g., cross validation.

% We illustrate this phenomenon in Figure \ref{fig:type-choice-vis} which visualizes the clustering results of each customer. Recall that the optimal assignment is to allocate customers with service rate no smaller than 5 to one queue and the rest customers to the other. When group number is smaller than 5, we can benefit from adding customer group by separating customer with service rate 5 from those with service rate 1. On the other hand, when group number is larger than 6, since history data is limited, the problem of model misspecification is more severe, leading to a worse cluster result. Overall, the marginal benefit of finer grouping decreases as we increase the number of customer types. Consequently, in practical applications, an overly granular grouping is generally superfluous and unnecessary, given the escalating risk of model misspecification.
% \begin{figure}[ht]
% 	\centering
% 	\includegraphics[width=\linewidth]{figure/type-choice-vis.eps}
% 	\caption{Case study: clustering customers with exponential service time into $m = 2, 4, 6, 8$ groups. \label{fig:type-choice-vis}}
% \end{figure}
}

\subsection{Non-identical Waiting Cost}
\label{sec:c-mu}
In previous subsections, we analyze how to minimize the expected waiting time of customers. In practice, however, the costs of waiting for different types of customers could be different. In this subsection, we extend the discussion to a cost minimization problem. In the stochastic scheduling setting, a well-known \(c\mu\) rule is proposed in the literature to minimize the total expected cost. Particularly, the \(c\mu\) rule assigns priorities to different customer types by calculating the product of unit waiting cost \(c\) and the service rate \(\mu\). \citet{CoxS61} show that \(c\mu\) rule is optimal for the \(M/G/1\) queue in a preemptive system. However, as discussed by \citet{HuB09}, the priority sequencing problem can be impractical in some settings when preemptive is not feasible. As an alternative, partitioning can be used to improve the efficiency of the system. In this subsection, we consider how to minimize the overall expected waiting cost when different types of customers have different costs per unit waiting time. We show that a modified version of the \(c\mu\) rule is optimal in the context of queueing partition under certain assumptions.

We follow the notations in Section \ref{sec:model}. We further let \(c_i\) denote the cost per unit waiting time of customer type \(i\) for \(i = 1, 2, \dots, n\). Then by the Pollaczek–Khinchine formula, given partition \(\boldsymbol{\alpha}\) and assignment \(X\), the expected waiting cost of the \(i\)th customer type can be calculated as
$$
\sum_{j\in R(\boldsymbol{\alpha})}
\left(\frac{\sum_{i=1}^n \lambda_i X_{ij}v_i/2}{\alpha_j^2 - \alpha_j\sum_{i=1}^n \frac{\lambda_i X_{ij}}{\mu_i}}
\right)c_i.
$$
Thus, the total expected waiting cost is
$$
\hat{f}(X, \boldsymbol{\alpha}) := \sum_{j\in R(\boldsymbol{\alpha})}
\left(\frac{\sum_{i=1}^n \lambda_i X_{ij}v_i}{\alpha_j^2 - \alpha_j\sum_{i=1}^n \frac{\lambda_i X_{ij}}{\mu_i}}
\cdot \frac{\sum_{i=1}^n X_{ij} \lambda_i c_i}{2\sum_{i=1}^n \lambda_i}\right).
$$
Replacing \(f(X, \boldsymbol{\alpha})\) with \(\hat{f}(X, \boldsymbol{\alpha})\) in \eqref{eq:k-qa-w} (\eqref{eq:k-qp-w}, resp.), we have the assignment (partition, resp.) problem with non-identical unit waiting cost, which we will denote as (\(k\)-QAP-\(c\)) (\(k\)-QPP-\(c\), resp.). In particular, the model in Section \ref{sec:model} can be viewed as a special case where \(c_1, c_2, \dots, c_n\) are identical.

Now we show that \(k\)-QAP-\(c\) and \(k\)-QPP-\(c\) can be transformed to \(k\)-QAP and \(k\)-QPP with a proper transformation of variables. We define \(\tilde{\lambda}_i  = c_i\lambda_i, \tilde{v}_i = v_i/c_i\) and \(\tilde{\mu}_i = c_i\mu_u\) for \(i = 1, 2, \dots, n\). Then, the \(k\)-QAP-\(c\) problem with parameters \(\{(\lambda_i, \mu_i, v_i)\}_{i=1}^n\) is equivalent to the \(k\)-QAP problem with parameters \(\{(\tilde{\lambda}_i, \tilde{\mu}_i, \tilde{v}_i)\}_{i=1}^n\). Therefore, the results obtained in previous subsections can be directly extended to the non-identical cost setting. Particularly, we have the following version of the modified \(c\mu\) rule in the context of the queueing partition problem.

\begin{theorem}
	\label{thm:multi-opt-alpha-cost}
	Given \(\lambda_i, c_i, i = 1, 2, \dots, n\) and \(\mathcal{I}\) with \(\mu_i, v_i, i = 1, 2, \dots, n\) such that \(c_1\mu_1 > \cdots > c_n \mu_n\) and \(\mathcal{I}\) satisfies Assumption \ref{as:reg-2seg}. If the feasible set of the \(k\)-QAP-\(c\) problem is non-empty, then there exist an optimal solution \((X^*, \boldsymbol{\alpha}^*)\) and \(1 = i_0 \le i_1 \le i_2 \le \cdots \le i_{k-1} \le i_k = n + 1\) such that
	\begin{align*}
	X^*_{i j} = \begin{cases}
	1 &\text{if }i_{j-1} \le i < i_{j} \\
	0 &\text{otherwise},
	\end{cases}
	\end{align*}
	for \(j = 1, 2, \dots, k\).
\end{theorem}

Theorem \ref{thm:multi-opt-alpha-cost} shows that under Assumption \ref{as:reg-2seg}, the optimal partition must be to partition the customer types according to the product of per unit waiting cost and the service rate, i.e., \(c\mu\). In other words, it is optimal to assign customers with consecutive \(c\mu\) into the same queue. Such a result can give rise to an efficient algorithm for the queue partition problem under different waiting costs. Particularly, note that given \(X\), we can verify that \(\hat{f}(X, \boldsymbol{\alpha})\) is jointly convex in \(\boldsymbol{\alpha}\). Therefore, we can enumerate all candidates of \(X\)
(by Theorem \ref{thm:multi-opt-alpha-cost}, there are at most \(\mathcal{O}(n^{k-1})\) of them) and for each \(X\) solve a convex optimization for the optimal \(\boldsymbol{\alpha}\). The detailed algorithm is given in Algorithm~\ref{alg:multi-opt-alpha-cost}.

\begin{algorithm}[!htp]
	\caption{Algorithm for Solving Queue Partition Problem under Non-Identical Waiting Cost}
	\label{alg:multi-opt-alpha-cost}
	\begin{algorithmic}[1]
		\Require $\lambda_i, \mu_i, v_i, c_i, i = 1, \dots, n$ such that $c_1\mu_1 > c_2\mu_2 > \cdots > c_n\mu_n$
		\Ensure $X^*\in\mathbb{R}^{n\times k}, \boldsymbol{\alpha}^*\in\mathbb{R}^k$
		
		\State Initialize $f^*\gets+\infty$
		\State $i_0\gets 0, i_{k+1} \gets n+1$
		\For{$1 \le i_1 \le i_2 \le \cdots \le i_k \le n + 1$}
		\State $X_{ij} = \begin{cases}
		1 &\text{if }i_{j-1} \le i < i_{j} \\
		0 &\text{Otherwise}
		\end{cases}$
		\State $f \gets \min_{\boldsymbol{\alpha}} f(X, \boldsymbol{\alpha})$, $\boldsymbol{\bar{\alpha}} \gets \arg\min_{\boldsymbol{\alpha}} f(X, \boldsymbol{\alpha})$
		\If {$f < f^*$}
		\State $f^* \gets f, X^* \gets X, \boldsymbol{\alpha}^* \gets \boldsymbol{\bar{\alpha}}$
		\EndIf
		\EndFor
	\end{algorithmic}
\end{algorithm}

\section{Extension: Analysis of Sojourn Time}
\label{sec:sojourn}
In previous sections, we investigate the optimal partition and assignment of queues by using waiting time as the performance metric. However, in practice, the decision maker may also care about other performance metrics of the system. One metric that is of special interest is the sojourn time, which is the waiting time plus the service time. Since in our setting we partition the queue by splitting the service capacity, the service times can be dramatically increased if a partition is made. In this section, we consider the optimal partition and assignment problem with the objective of minimizing the sojourn time. We show that the optimal structure discussed earlier also applies to the sojourn time minimization setting.

Specifically, following the notations in Section \ref{sec:model}, by the Pollaczek–Khinchine formula, given partition \(\boldsymbol{\alpha}\) and assignment \(X\), the expected sojourn time can be calculated as
$$
\tilde{f}(X, \boldsymbol{\alpha}) := \frac{1}{2\sum_{i=1}^n \lambda_i}\sum_{j\in R(\boldsymbol{\alpha})}
\left(\frac{\sum_{i=1}^n \lambda_i X_{ij}v_i \cdot \sum_{i=1}^n X_{ij} \lambda_i}{\alpha_j^2 - \alpha_j\sum_{i=1}^n \frac{\lambda_i X_{ij}}{\mu_i}} +
\frac{2\sum_{i=1}^n \frac{\lambda_i X_{ij}}{\mu_i}}{\alpha_j}\right).
$$
Here, the first term corresponds to the expected waiting time, which is the same as in \(f(X, \boldsymbol{\alpha})\). The second term is the expected service time. In the following, we first show that in general, by making proper partition and assignment, one can reduce the expected sojourn time. Moreover, the improvement can be arbitrarily large. We have the following proposition.

\begin{proposition}
	\label{pro:improvement-s}
	For any \(\epsilon > 0\), there exist input parameters \(\{\lambda_i\}_{i=1}^n, \{\mu_i\}_{i=1}^n, \{v_i\}_{i=1}^n\) and \(X, \boldsymbol{\alpha}\) such that sojourn time \(\tilde{f}(X, \boldsymbol{\alpha}) < \epsilon \tilde{f}(\bar{X}, \bar{\boldsymbol{\alpha}})\), where the pooling decision is
	$$
	\bar{X} = \begin{pmatrix}
	0 &\cdots & 0 &1 \\
	\vdots &\vdots &\vdots &\vdots \\
	0 &\cdots & 0 &1 \\
	\end{pmatrix}, \quad
	\bar{\boldsymbol{\alpha}} = \begin{pmatrix}
	0 \\ \vdots \\ 0 \\ 1
	\end{pmatrix}.
	$$
\end{proposition}

To illustrate, consider the same example as is shown in Section \ref{sec:model}. We show that even if we measure the system by the expected sojourn time, doing a proper partition can improve the system efficiency significantly. Specifically, recall that there are \(k = 2\) queues and \(n = 2\) types of customers with \(\lambda_1 = t, \lambda_2 = 1\). \(S_1\) and \(S_2\) both follow exponential distributions with \(\mu_1 = 2t, \mu_2 = \frac{2}{t}\). For \(0 < t \le \frac{1}{4}\), under partition \(\boldsymbol{\alpha} = (1 - t, t)\) and assignment \(X = \begin{pmatrix}
1 &0\\
0 &1\\
\end{pmatrix}\), we have \(\tilde{f}(X, \boldsymbol{\alpha}) = \left(\frac{1/4}{(1-t)(1/2-t)}+\frac{t^2/4}{t^2/2}+\frac{1/2}{1-t}+\frac{1}{2}\right)\frac{1}{t+1}\). We can also calculate \(\tilde{f}(\bar{X}, \bar{\boldsymbol{\alpha}}) = \frac{1}{t+1}\left(\frac{(\frac{1}{4t}+\frac{t^2}{4})(t+1)}{1/2-t/2}+\frac{1}{2}+\frac{t}{2}\right)\). Therefore, \(\lim_{t\to0}\frac{\tilde{f}(X, \boldsymbol{\alpha})}{\tilde{f}(\bar{X}, \bar{\boldsymbol{\alpha}})} = \lim_{t\to0}\frac{2}{1/(2t)} = 0\). The intuition is similar to before. When a plethora of customers who can be served quickly are blocked by some slow minority, we can allocate some dedicated resources for the vast majority and boost the average sojourn time dramatically.

Moreover, we show that partitioning customers into two continuous segments ranked by their service rates can also behave poorly under the expected sojourn time measure.

\vspace{0.4 cm}
\begin{example}
	Consider the same instance as in Example \ref{eg:3better}. For \(\alpha = 0.95\), if we partition the customer based on the expected service time, the best assignment is \(x_1 = 1, x_2 = 0.66, x_3 = 0\) with expected sojourn time 1.14. However, for assignment \(x_1 = 1, x_2 = 0, x_3 = 1\), the expected sojourn time is 0.77, decreased by 32\%. Similar to the waiting time case, the first and the third types of customers have smaller service time variance. By grouping them, we can benefit from pooling while having a small service time variance within that group.\(\hfill\square\)
\end{example}
\vspace{0.4 cm}

In the following, for an instance \(\mathcal{I}\), we adopt the same notation \(\mathcal{C}(\mathcal{I}, k)\) as defined in \eqref{eq:def-config-k}. In the following theorem, we state that the optimal structure of sojourn time minimization is the same as that of waiting time minimization.

\begin{theorem}
	\label{thm:multi-fixed-alpha-sojourn}
	For any given \(k\), suppose \(\lambda_1, \dots, \lambda_n, \mu_1 > \cdots > \mu_n, v_1, \dots, v_n\) and \(\alpha\) are given. Then there exists an optimal solution \(X\) to the sojourn time minimization problem such that
	$$\begin{cases}
	X_{i \pi_j} = 1 & i \in T_j\\
	X_{i \pi_j} \in [0, 1] & i \in E_{j}\\
	X_{i \pi_j} = 0 & \text{otherwise},
	\end{cases}$$
	holds for a permutation \(\boldsymbol{\pi} = (\pi_1, \dots, \pi_k)\) of \(\{1, 2, \dots, k\}\) and \((T_1, \dots, T_k, E_1, \dots, E_{k}) \in \mathcal{C}(\mathcal{I}, k)\).
\end{theorem}

With this property, we can use Algorithm \ref{alg:multi-fixed-alpha} to find the optimal solution of the sojourn time minimization problem by replacing the computation of the objective function. Moreover, we can find the optimal partition by incorporating the optimization of \(\boldsymbol{\alpha}\) in line \ref{lst:line:multi-variate} of Algorithm \ref{alg:multi-fixed-alpha}.

{
It is worth noting that although the optimal assignment under sojourn time minimization possesses the same structural property as that under waiting time minimization, the types of customers grouped in the same queue may be different under these two settings. Consider the instance in Example \ref{eg:3better}. When we vary the second moment of the third-type customers' service time, the customer groups can be different when objective varies. We illustrate this phenomenon in Figure \ref{fig:w_vs_s}, where the three bars represent the assignment of customer type 1, 2, 3 respectively. In particular, when $\mathbb{E}[S_3^2] = 79$, if we aim to minimize the expected waiting time, then it is optimal to group the first two types of customers. In contrast, if we focus on the expected sojourn time minimization, grouping the first type of customers with the third type is more beneficial.
Moreover, for the joint assignment and partition problem, we conjecture
that under Assumption \ref{as:reg-2seg} it is optimal to assign customers with
consecutive service rates to the same queue. We have conducted
extensive numerical experiments which validate such a
conjecture. However, providing a formal proof remains an open challenge, which we leave for future research.

\begin{figure}[ht]
	\centering
	\includegraphics[width=0.8\linewidth]{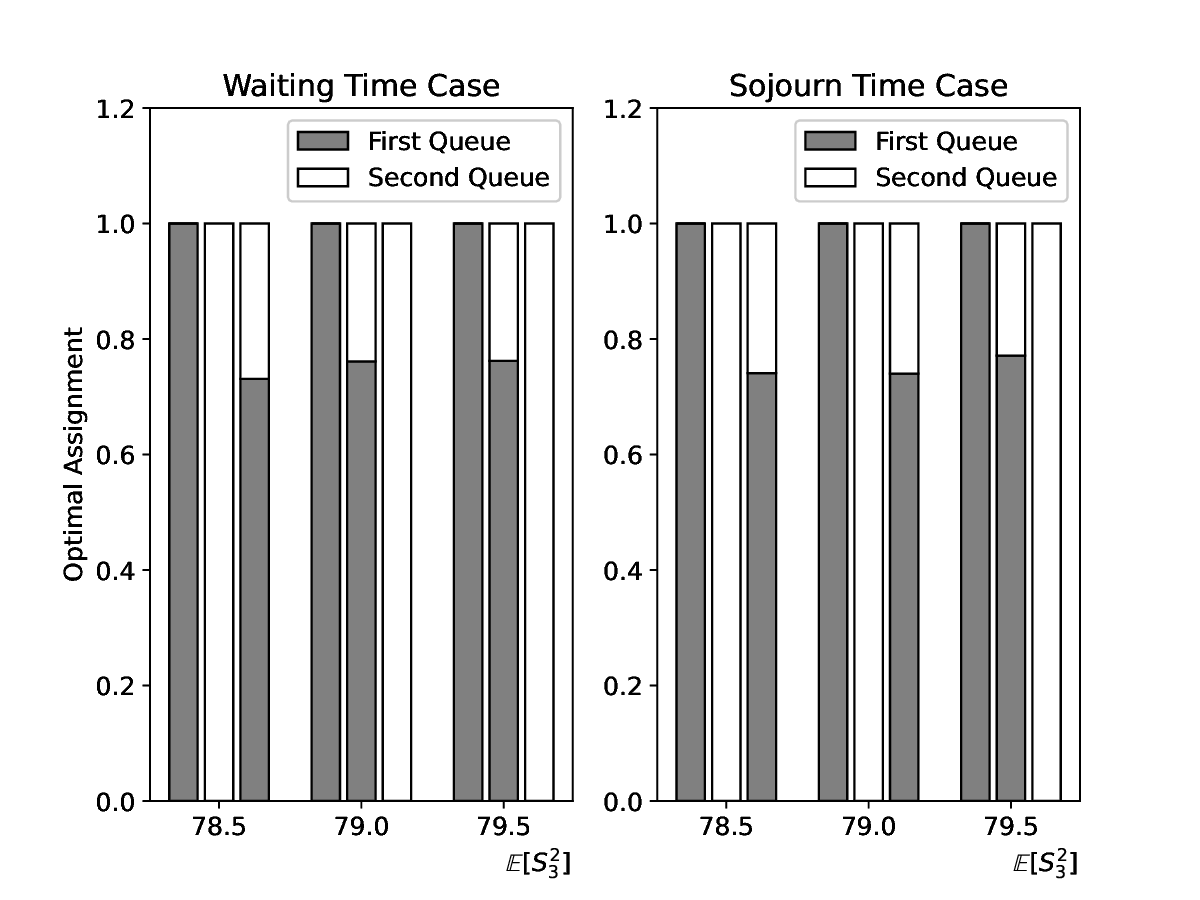}
	\caption{Different grouping strategies for waiting time minimization and sojourn time minimization when there are three different types of customers. \label{fig:w_vs_s}}
\end{figure}
}

\section{Concluding Remarks}
\label{sec:conclusion}
In this paper, we analyzed the optimal partition and assignment strategy for a queueing system with an infinite waiting room and FCFS routing policy. We illustrated that it is possible to improve service efficiency by partitioning the server, and the potential for improvement could be large. We analyzed the optimal structure of the queue partition and assignment and provided efficient algorithms to find the solution. 

We extended the analysis of waiting time minimization to sojourn time optimization. We also extended the analysis to the case where waiting costs for different customer types are non-identical and built a bridge between the queue partition problem and the priority sequencing problem by the famous \(c\mu\) rule. Overall, our work presented a comprehensive analysis of the queue partition problem, which could potentially improve the efficiency of the queueing system.

In practice, one may also want to decide how many queues we should partition into. For this question, we can take a large enough $k$ and solve the partition problem. We conjecture if there are $n$ types of customers, then the optimal partition $k$ is no more than $n$. We have conducted extensive numerical experiments which validate such a conjecture. Unfortunately, we are not able to provide a proof yet. We leave this as a future research direction. However, we note that, in practice, it is impractical to partition a system into too many queues. Therefore, such a result may be less consequential in practice than in theory. Further research could also be conducted to explore the effects of the inexact mean and second-order moment of the service time distribution on the optimal partition decision. A robust or distributionally robust approach may be adopted in those settings. Additionally, algorithms could be developed to learn the service time distribution in an online fashion. These studies could provide valuable insights to enhance the robustness and adaptiveness of the queueing system.

%\THEEndNotes
%\begingroup \parindent 0pt \parskip 4ex
%\def\enotesize{\normalsize} 
%\theendnotes
%\endgroup

% Appendix here
% Options are (1) APPENDIX (with or without general title) or
%             (2) APPENDICES (if it has more than one unrelated sections)
% Outcomment the appropriate case if necessary
%
\newpage
\begin{APPENDIX}{Supplementary Derivations and Proofs}
\noindent \textbf{Proof of Theorem \ref{thm:fix-alpha}.} 
{First we consider \(\alpha = 1\) (\(\alpha = 0\)).} In this case, \(\boldsymbol{x} = \boldsymbol{1}\) (\(\boldsymbol{x} = \boldsymbol{0}\)) is the only feasible solution, and the statement is true.

{When \(0 < \alpha < 1\), we first show that \(\boldsymbol{x} = \boldsymbol 0\) or \(\boldsymbol{x} = \boldsymbol 1\), i.e., assigning all customers to one queue, cannot be optimal.
To proceed, we first introduce some notations. For queue with arrival rate $\lambda$ and service time distribution $S$, we denote the zeroth order load as $C = \lambda\mathbb{E}[S^0]$, the first order load as $B = \lambda\mathbb{E}[S^1]$, and the second order load as $A = \lambda\mathbb{E}[S^2]$. We denote the zeroth, first, second order load for the first (second, resp) sub-queue under a unit server as $C_x, B_x, A_x$ ($C_y, B_y, A_y$, resp.). Therefore, we have
\begin{eqnarray}
&A_x = \sum_{i=1}^n \lambda_i v_i x_i, B_x = \sum_{i=1}^n \frac{\lambda_i}{\mu_i}x_i, C_x = \sum_{i=1}^n \lambda_ix_i, \nonumber \\
&A_y = \sum_{i=1}^n \lambda_i v_i (1-x_i), B_y = \sum_{i=1}^n \frac{\lambda_i}{\mu_i}(1-x_i), C_y = \sum_{i=1}^n \lambda_i(1-x_i).
\end{eqnarray}}
We denote $s = \frac{B_x}{\alpha-B_x}, t = \frac{B_y}{1 - \alpha-B_y}$. We further denote the Lagrangian function \(\mathcal{L}(\boldsymbol{x}, \alpha, \boldsymbol{p}, \boldsymbol{q}) = \frac{A_x \cdot C_x}{\alpha(\alpha - B_x)} + \frac{A_y \cdot C_y}{(1 - \alpha)(1 - \alpha - B_y)} - \sum_i x_i p_i - \sum_i (1 - x_i) q_i\), where \(p_i \ge 0, q_i \ge 0\) (we remove the constant factor \(\frac{1}{2\sum_{i=1}^n \lambda_i}\) for the ease of expression). From the KKT condition \(\frac{\partial\mathcal{L}}{\partial x_i} = 0\), we have
\begin{eqnarray}\label{eq:general-kkt-x}
p_i - q_i
=&& \left(\frac{C_x}{\alpha(\alpha - B_x)} - \frac{C_y}{(1-\alpha)(1-\alpha-B_y)}\right) \lambda_i v_i \nonumber \\
&&+ \left(\frac{A_x}{\alpha(\alpha - B_x)} - \frac{A_y}{(1-\alpha)(1-\alpha-B_y)}\right) \lambda_i \nonumber \\
&&+ \left(\frac{A_xC_x}{\alpha(\alpha - B_x)^2} - \frac{A_yC_y}{(1-\alpha)(1-\alpha-B_y)^2}\right) \frac{\lambda_i}{\mu_i}.
\end{eqnarray}
If \(\boldsymbol{x} = \boldsymbol 0\) is feasible, then for \(\boldsymbol{x} = \boldsymbol 0\),
\[p_i - q_i = -\left(\frac{C_y}{(1-\alpha)(1-\alpha-B_y)} \lambda_i v_i + \frac{A_y}{(1-\alpha)(1-\alpha-B_y)}\lambda_i + \frac{A_yC_y}{(1-\alpha)(1-\alpha-B_y)^2}\frac{\lambda_i}{\mu_i}\right) < 0.
\]
From non-negativity of \(p_i, q_i\) we know \(q_i > 0\). However, complementary slackness condition states that \((1 - x_i)q_i = 0\), thus \(\boldsymbol{x} = \boldsymbol 0\) is not optimal.

Similarly, if \(\boldsymbol{x} = \boldsymbol 1\) is feasible, then for \(\boldsymbol{x} = \boldsymbol 1\),
$$p_i - q_i = \frac{C_x}{\alpha(\alpha - B_x)}\lambda_iv_i + \frac{A_x}{\alpha(\alpha - B_x)} \lambda_i + \frac{A_xC_x}{\alpha(\alpha - B_x)^2} \frac{\lambda_i}{\mu_i} > 0.$$

From non-negativity of \(p_i, q_i\) we know \(p_i > 0\). However, complementary slackness condition states that \(x_i p_i = 0\), thus \(\boldsymbol{x} = \boldsymbol{1}\) is not optimal.

{Now we can restrict our consideration to \(\boldsymbol{x} \neq \boldsymbol{0}, \boldsymbol{1}\). We first formulate a representative function which maps each customer type to a scalar value. In particular, for any feasible assignment \(\boldsymbol{x}\), we define function \(g_{\boldsymbol{x}}\) as
\begin{equation}\label{eq:v-u-genearl}
g_{\boldsymbol{x}}(u) = \left( \frac{\gamma s}{\alpha} - \frac{t}{1-\alpha} \right) \frac{B_y}{A_y} v_{\mathcal{I}}(u) + \left(\frac{\beta s}{\alpha}-\frac{t}{1-\alpha}\right) \frac{B_y}{C_y} + \left(\frac{u^2s^2}{\alpha}-\frac{t^2}{1-\alpha}\right) u.
\end{equation}
It is readily shown that representative function only depends linearly on the first two moment of the customers' service time, i.e., $u$ and $v_{\mathcal{I}}(u)$.

We complete the proof by showing that, under the optimal assignment policy, customer of type $i$ will be allocated to the first server only when the corresponding scalar value $g_{\boldsymbol{x}}(\cdot)$ is non-positive. Consequently, due to linearity, the optimal assignment then corresponds to a polytope partition in a two-dimensional space of the first two moments. In particular, when there are two queues, the two-dimensional space are partitioned by a line into two subspaces.}

Specifically, since \(\boldsymbol{x} \neq \boldsymbol{0}, \boldsymbol{1}\), we have \(A_x, A_y, B_x, B_y, C_x, C_y \neq 0\). Multiplying both sides of \eqref{eq:general-kkt-x} by \(\frac{B_y^2}{A_yC_y}\frac{1}{\lambda_i}\) yields
\begin{eqnarray*}
\frac{B_y^2}{A_yC_y}\frac{p_i - q_i}{\lambda_i}
=&& \left(\frac{\frac{C_xB_y}{C_yB_x}}{\alpha(\alpha - B_x)/B_x} - \frac{1}{(1-\alpha)(1-\alpha-B_y)/B_y}\right) \frac{B_y}{A_y} v_i \nonumber \\
&&+ \left(\frac{\frac{A_xB_y}{A_yB_x}}{\alpha(\alpha - B_x)/B_x} - \frac{1}{(1-\alpha)(1-\alpha-B_y)/B_y}\right) \frac{B_y}{C_y} \nonumber \\
&&+ \left(\frac{\frac{C_xB_y}{C_yB_x}\frac{A_xB_y}{A_yB_x}}{\alpha(\alpha - B_x)^2/B_x^2} - \frac{1}{(1-\alpha)(1-\alpha-B_y)^2/B_y^2}\right) \frac{1}{\mu_i}.
\end{eqnarray*}

Denote \(\gamma = \frac{C_x}{B_x}\frac{B_y}{C_y}, \beta = \frac{A_x}{B_x}\frac{B_y}{A_y}, u = \sqrt{\gamma\beta}\). Recall that \(s = \frac{B_x}{\alpha-B_x}, t = \frac{B_y}{1 - \alpha-B_y}\), the above equation can then be simplified as
\begin{eqnarray}\label{eq:general-kkt-x-sim}
\left(\frac{\gamma s}{\alpha} - \frac{t}{1-\alpha}\right) \frac{B_y}{A_y} v_i
+ \left(\frac{\beta s}{\alpha}-\frac{t}{1-\alpha}\right) \frac{B_y}{C_y}
+ \left(\frac{u^2s^2}{\alpha}-\frac{t^2}{1-\alpha}\right) \frac{1}{\mu_i} = \frac{B_y^2}{A_yC_y}\frac{p_i - q_i}{\lambda_i}.
\end{eqnarray}

Recall definition \eqref{eq:v-u-genearl}, we know that \eqref{eq:general-kkt-x-sim} can be rewritten as \(g_{\boldsymbol{x}}(\frac{1}{\mu_i}) = \frac{B_y^2}{A_yC_y}\frac{p_i - q_i}{\lambda_i}\). By complementary slackness conditions, for any optimal assignment \(\boldsymbol{x}^*\), it holds that
$$
g_{\boldsymbol{x}^*}\left(\frac{1}{\mu_i}\right) > 0 \Rightarrow x^*_i = 0, \quad g_{\boldsymbol{x}^*}\left(\frac{1}{\mu_i}\right) < 0 \Rightarrow x^*_i = 1.
$$

{We then finish the proof by considering two cases.}

\textbf{Case 1.} \(\frac{\gamma s}{\alpha} = \frac{t}{1-\alpha}\). Then \(g_{\boldsymbol{x}^*}(u)\) reduces to a linear function. When \(g_{\boldsymbol{x}^*}(u) \not\equiv 0\), \(g_{\boldsymbol{x}^*}(u)\) intersects with the positive half of \(X\)-axis for at most once. By KKT condition we know there exists \(i_0 \in \{1, 2, \dots, n\}\) such that either \(\begin{cases} x_i^{ * } = 0 & i < i_0 \\ x_i^{ * } = 1 & i > i_0 \\ x_i^{ * } \in [0, 1] & i = i_0. \end{cases}\) or \(\begin{cases} x_i^{ * } = 1 & i < i_0 \\ x_i^{ * } = 0 & i > i_0 \\ x_i^{ * } \in [0, 1] & i = i_0. \end{cases}\). Construct \(U = \{i: i > i_0\}, E = \{i: i = i_0\}, L = \{i: i < i_0\}\). By considering a linear function passing through \((\frac{1}{\mu_{i_0}}, v_0)\) with large enough slope, we have \((U, E, L) \in \mathcal{C}(\mathcal{I})\). When \(g_{\boldsymbol{x}^*}(u) \equiv 0\), we have
\begin{align*}
\frac{\gamma s}{\alpha} = \frac{t}{1-\alpha},
\frac{\beta s}{\alpha} = \frac{t}{1-\alpha},
\frac{u^2s^2}{\alpha} = \frac{t^2}{1-\alpha},
\end{align*}
which implies \(\alpha = \frac{1}{2}\) and \(\frac{\frac{1}{2} - B_x}{\frac{1}{2} - B_y} = \frac{A_x}{A_y} = \frac{C_x}{C_y}\). For any \(i\neq j\), consider another assignment \(\boldsymbol{x}(\epsilon)\) where \(x_i(\epsilon) = x^{ * }_i + \left(\frac{A_x + C_x}{\alpha - B_x}\frac{\lambda_j}{\mu_j} + \lambda_j + \lambda_j v_j\right)\epsilon, x_j(\epsilon) = x^{ * }_j - \left(\frac{A_x + C_x}{\alpha - B_x}\frac{\lambda_i}{\mu_i} + \lambda_i + \lambda_i v_i \right) \epsilon\) and \(x_k(\epsilon)=x^{ * }_k, \forall k \neq i, j\). For \(i = 1, 2, \dots, n\), we denote \(A_i = \lambda_i v_i, B_i = \frac{\lambda_i}{\mu_i}, C_i = \lambda_i\). We define \(l(\epsilon) = f(\boldsymbol{x}(\epsilon), \alpha)\). With some computation, we have
\begin{align*}
l''(\epsilon) = -\frac{1}{\sum_{t=1}^n \lambda_t}
&\left(
\frac{\left(A_x(B_jC_i - B_iC_j) + C_x(A_jB_i - A_iB_j) + (\alpha-B_x)(A_jC_i-A_iC_j)\right)^2}{\alpha\left(\alpha-B_x-(A_jB_i+C_jB_i-A_iB_j-C_iB_j)\epsilon\right)^3} \right. \\+&
\left. \frac{\left(-A_y(B_jC_i - B_iC_j) - C_y(A_jB_i - A_iB_j) - (1-\alpha-B_y)(A_jC_i-A_iC_j)\right)^2}{(1 - \alpha)\left(1-\alpha-B_y+(A_jB_i+C_jB_i-A_iB_j-C_iB_j)\epsilon\right)^3}
\right),
\end{align*}
where the second line uses the fact that
\begin{align*}
x_i(\epsilon) = x^{ * }_i + \left(\frac{A_x + C_x}{\alpha - B_x}\frac{\lambda_j}{\mu_j} + \lambda_j + \lambda_j v_j\right)\epsilon = x^{ * }_i + \left(\frac{A_y + C_y}{1 - \alpha - B_y}\frac{\lambda_j}{\mu_j} + \lambda_j + \lambda_j v_j\right)\epsilon, \\
x_j(\epsilon) = x^{ * }_j - \left(\frac{A_x + C_x}{\alpha - B_x}\frac{\lambda_i}{\mu_i} + \lambda_i + \lambda_i v_i \right) \epsilon = x^{ * }_j - \left(\frac{A_y + C_y}{1 - \alpha - B_y}\frac{\lambda_i}{\mu_i} + \lambda_i + \lambda_i v_i \right) \epsilon.
\end{align*}
Thereby, \(l\) is concave in \(\epsilon\) in the neighborhood of 0. Together with \(l'(0) = 0\), we can claim that there exists an optimal solution with \(\boldsymbol{x}^{ * } = \boldsymbol{1}\) or \(\boldsymbol{x}^{ * } = \boldsymbol{0}\), where contradiction arises.

\textbf{Case 2.} \(\frac{\gamma s}{\alpha} \neq \frac{t}{1-\alpha}\). In this case, consider linear function \(l(u) = -\frac{\left(\frac{\beta s}{\alpha}-\frac{t}{1-\alpha}\right) \frac{B_y}{C_y} + \left(\frac{u^2s^2}{\alpha}-\frac{t^2}{1-\alpha}\right) u}{\left( \frac{\gamma s}{\alpha} - \frac{t}{1-\alpha} \right) \frac{B_y}{A_y}}\) and its induced partition \((U, E, L)\). It is readily shown that \((U, E, L) \in \mathcal{C}(\mathcal{I})\). When \(\frac{\gamma s}{\alpha} > \frac{t}{1-\alpha}\), the KKT condition implies that an optimal solution \(\boldsymbol{x}^*\) satisfies
$$\begin{cases} x_i^{ * } = 0 & i\in U \\ x_i^{ * } = 1 & i\in L \\ x_i^{ * } \in [0, 1] & i\in E. \end{cases}$$
Similarly, when \(\frac{\gamma s}{\alpha} < \frac{t}{1-\alpha}\), an optimal solution \(\boldsymbol{x}^*\) satisfies
$$\begin{cases} x_i^{ * } = 1 & i\in U \\ x_i^{ * } = 0 & i\in L \\ x_i^{ * } \in [0, 1] & i\in E. \end{cases}$$ \(\hfill\square\)

\bigskip
\noindent \textbf{Proof of Theorem \ref{thm:fix-alpha-exp}.}
{First we note that if the feasible set of QAP is nonempty, then an optimal solution exists and is attainable.} To see this, for any given \(\alpha\), as \(\sum_{i=1}^n\frac{\lambda_i}{\mu_i}x_i \to \alpha\)
or \(\sum_{i=1}^n\frac{\lambda_i}{\mu_i}x_i \to \alpha - (1 - \sum_{i=1}^n\frac{\lambda_i}{\mu_i})\), \(f(\boldsymbol{x}, \alpha) \to +\infty\). Also, \(f(\boldsymbol{x}, \alpha)\) is continuous in \(\boldsymbol{x}\). Thereby, a minimum can be obtained in the interior of set \(\mathcal{F}_{\alpha}\). In addition, since the set \(\left\{\boldsymbol{x}: 0 \le x_i \le 1, \forall i\right\}\) is closed, we can conclude that optimal solution \(\boldsymbol{x}^*\) exists.

{Next we prove that for any optimal solution \(\boldsymbol{x}^*\), \(\mathcal{M}^* = \left\{i: 0 < x^*_i < 1\right\}\) is either empty or singleton. We prove it by contradiction and show that if two different types of customers are assigned to both queues, swapping these customers between the two queues will always decrease the expected waiting time.}

Specifically, we consider an optimal assignment \(\boldsymbol{x}^*\) with cardinality \(|\mathcal{M}^*| \ge 2\). Denote two of \(\mathcal{M}^*\)'s elements as \(i\neq j\). Now consider another assignment \(\boldsymbol{x}(\epsilon)\) where \(x_i(\epsilon) = x^*_i + \frac{\mu_i}{\lambda_i}\epsilon, x_j(\epsilon) = x^*_j - \frac{\mu_j}{\lambda_j}\epsilon\),
and \(x_k(\epsilon) = x^*_k, \;\forall k \neq i, j\). We define
\begin{eqnarray}\label{eq:two-x}
l(\epsilon) := f(\boldsymbol{x}(\epsilon), \alpha).
\end{eqnarray}
Therefore, we have \(l''(0) = -\frac{(\mu_i - \mu_j)^2}{\mu_i\mu_j} \frac{1}{\sum_{k=1}^n \lambda_k} \left(\frac{1}{\alpha^2 - \alpha\sum_{k=1}^n \frac{\lambda_k x^*_k}{\mu_k}} + \frac{1}{(1 - \alpha)^2 - (1 - \alpha)\sum_{k=1}^n \frac{\lambda_k(1 - x^*_k)}{\mu_k}} \right) < 0\). Since \(\epsilon = 0\) is an interior point and \(l(\epsilon)\) is strictly concave at \(\epsilon = 0\), there exists \(\delta > 0\) such that for \(\epsilon = \pm\delta\), \(\boldsymbol{x}(\epsilon)\) is feasible and \(l(0) > \min\left\{l(-\delta), l(\delta)\right\}\). Now we find a feasible assignment with a smaller objective function value, which contradicts with the optimality of \(\boldsymbol{x}\). Therefore, \(\mathcal{M}^* = \left\{i: 0 < x^*_i < 1\right\}\) must be either empty or singleton.

{In the following, we prove the rest of the statement, i.e., the ``three continuous segments'' property, by showing that the first two moments of each customer type's service time distribution lie on a (possibly degenerated) quadratic curve. Then by Theorem \ref{thm:fix-alpha}, the partition line then interacts with this quadratic curve for at most twice.}

For simplicity, we denote
\begin{eqnarray}
&A_x = \sum_{i=1}^n \frac{\lambda_i}{\mu_i^2}x_i, B_x = \sum_{i=1}^n \frac{\lambda_i}{\mu_i}x_i, C_x = \sum_{i=1}^n \lambda_ix_i; \nonumber \\
&A_y = \sum_{i=1}^n \frac{\lambda_i}{\mu_i^2}(1-x_i), B_y = \sum_{i=1}^n \frac{\lambda_i}{\mu_i}(1-x_i), C_y = \sum_{i=1}^n \lambda_i(1-x_i); \nonumber \\
&s = \frac{B_x}{\alpha-B_x}, t = \frac{B_y}{1 - \alpha-B_y}. \label{eq:def-st}
\end{eqnarray}
Then we can write \(f(\boldsymbol{x}, \alpha) = \frac{1}{\sum_{i=1}^n \lambda_i}\left(\frac{A_x \cdot C_x}{\alpha(\alpha - B_x)} + \frac{A_y \cdot C_y}{(1 - \alpha)(1 - \alpha - B_y)}\right)\) and
$$
\sum_{i=1}^n \frac{\lambda_i x_i}{\mu_i} < \alpha < \sum_{i=1}^n \frac{\lambda_i x_i}{\mu_i} + 1 - \sum_{i=1}^n \frac{\lambda_i}{\mu_i}
$$
can be reformulated as
$$
B_x < \alpha < 1 - B_y.
$$
First we consider \(\alpha = 1\). In this case, \(\boldsymbol{x} = \boldsymbol{1}\) is the only feasible solution,
and the statement is true. Now we consider \(\frac{1}{2} \le \alpha < 1\). We first show that \(\boldsymbol{x} = \boldsymbol 0\) or \(\boldsymbol{x} = \boldsymbol 1\), i.e., assigning all customers to one queue, cannot be optimal. Denote the Lagrangian function \(\mathcal{L}(\boldsymbol{x}, \alpha, \boldsymbol{p}, \boldsymbol{q}) = \frac{A_x \cdot C_x}{\alpha(\alpha - B_x)} + \frac{A_y \cdot C_y}{(1 - \alpha)(1 - \alpha - B_y)} - \sum_i x_i p_i - \sum_i (1 - x_i) q_i\), where \(p_i \ge 0, q_i \ge 0\) (we remove the constant factor \(\frac{1}{\sum_{i=1}^n \lambda_i}\) for the ease of expression). From the KKT condition \(\frac{\partial\mathcal{L}}{\partial x_i} = 0\), we have
\begin{eqnarray}\label{eq:kkt-x}
p_i - q_i
=&& \left(\frac{C_x}{\alpha(\alpha - B_x)} - \frac{C_y}{(1-\alpha)(1-\alpha-B_y)}\right) \frac{\lambda_i}{\mu_i^2} \nonumber \\
&&+ \left(\frac{A_x}{\alpha(\alpha - B_x)} - \frac{A_y}{(1-\alpha)(1-\alpha-B_y)}\right) \lambda_i \nonumber \\
&&+ \left(\frac{A_xC_x}{\alpha(\alpha - B_x)^2} - \frac{A_yC_y}{(1-\alpha)(1-\alpha-B_y)^2}\right) \frac{\lambda_i}{\mu_i}.
\end{eqnarray}
If \(\boldsymbol{x} = \boldsymbol 0\) is feasible, then for \(\boldsymbol{x} = \boldsymbol 0\),
\[p_i - q_i = -\left(\frac{C_y}{(1-\alpha)(1-\alpha-B_y)}\frac{\lambda_i}{\mu_i^2} + \frac{A_y}{(1-\alpha)(1-\alpha-B_y)}\lambda_i + \frac{A_yC_y}{(1-\alpha)(1-\alpha-B_y)^2}\frac{\lambda_i}{\mu_i}\right) < 0.
\]
From non-negativity of \(p_i, q_i\) we know \(q_i > 0\). However, complementary slackness condition states that \((1 - x_i)q_i = 0\), thus \(\boldsymbol{x} = \boldsymbol 0\) is not optimal.

Similarly, if \(\boldsymbol{x} = \boldsymbol 1\) is feasible, then for \(\boldsymbol{x} = \boldsymbol 1\),
$$p_i - q_i = \frac{C_x}{\alpha(\alpha - B_x)}\frac{\lambda_i}{\mu_i^2} + \frac{A_x}{\alpha(\alpha - B_x)} \lambda_i + \frac{A_xC_x}{\alpha(\alpha - B_x)^2} \frac{\lambda_i}{\mu_i} > 0.$$

From non-negativity of \(p_i, q_i\) we know \(p_i > 0\). However, complementary slackness condition states that \(x_i p_i = 0\), thus \(\boldsymbol{x} = \boldsymbol 1\) is not optimal.

{Now we can restrict our consideration to \(\boldsymbol{x} \neq \boldsymbol{0}, \boldsymbol{1}\).} In this case, \(A_x, A_y, B_x, B_y, C_x, C_y \neq 0\). Multiplying both sides of \eqref{eq:kkt-x} by \(\frac{B_y^2}{A_yC_y}\frac{1}{\lambda_i}\) yields
\begin{eqnarray*}
\frac{B_y^2}{A_yC_y}\frac{p_i - q_i}{\lambda_i}
=&& \left(\frac{\frac{C_xB_y}{C_yB_x}}{\alpha(\alpha - B_x)/B_x} - \frac{1}{(1-\alpha)(1-\alpha-B_y)/B_y}\right) \frac{B_y}{A_y}\frac{1}{\mu_i^2} \nonumber \\
&&+ \left(\frac{\frac{A_xB_y}{A_yB_x}}{\alpha(\alpha - B_x)/B_x} - \frac{1}{(1-\alpha)(1-\alpha-B_y)/B_y}\right) \frac{B_y}{C_y} \nonumber \\
&&+ \left(\frac{\frac{C_xB_y}{C_yB_x}\frac{A_xB_y}{A_yB_x}}{\alpha(\alpha - B_x)^2/B_x^2} - \frac{1}{(1-\alpha)(1-\alpha-B_y)^2/B_y^2}\right) \frac{1}{\mu_i}.
\end{eqnarray*}

Denote \(\gamma = \frac{C_x}{B_x}\frac{B_y}{C_y}, \beta = \frac{A_x}{B_x}\frac{B_y}{A_y}, u = \sqrt{\gamma\beta}\). Recall that \(s = \frac{B_x}{\alpha-B_x}, t = \frac{B_y}{1 - \alpha-B_y}\), the above equation can then be simplified as
\begin{eqnarray}\label{eq:kkt-x-sim}
\left(\frac{\gamma s}{\alpha} - \frac{t}{1-\alpha}\right) \frac{B_y}{A_y}\frac{1}{\mu_i^2}
+ \left(\frac{\beta s}{\alpha}-\frac{t}{1-\alpha}\right) \frac{B_y}{C_y}
+ \left(\frac{u^2s^2}{\alpha}-\frac{t^2}{1-\alpha}\right) \frac{1}{\mu_i} = \frac{B_y^2}{A_yC_y}\frac{p_i - q_i}{\lambda_i}.
\end{eqnarray}

For a feasible assignment \(\boldsymbol{x}\), define a function \(g_{\boldsymbol{x}}\) as
\begin{equation}\label{eq-g-func}
g_{\boldsymbol{x}}(v) = \left( \frac{\gamma s}{\alpha} - \frac{t}{1-\alpha} \right) \frac{B_y}{A_y}v^2 + \left(\frac{\beta s}{\alpha}-\frac{t}{1-\alpha}\right) \frac{B_y}{C_y} + \left(\frac{u^2s^2}{\alpha}-\frac{t^2}{1-\alpha}\right) v.
\end{equation}

Then \eqref{eq:kkt-x-sim} can be rewritten as \(g_{\boldsymbol{x}}(\frac{1}{\mu_i}) = \frac{B_y^2}{A_yC_y}\frac{p_i - q_i}{\lambda_i}\). By complementary slackness conditions, for any optimal assignment \(\boldsymbol{x}^*\), it holds that
$$
g_{\boldsymbol{x}^*}\left(\frac{1}{\mu_i}\right) > 0 \Rightarrow x^*_i = 0, \quad g_{\boldsymbol{x}^*}\left(\frac{1}{\mu_i}\right) < 0 \Rightarrow x^*_i = 1.
$$

Note that \(g_{\boldsymbol{x}^*}(v)\) is a polynomial function with degree being at most two. {We finish the proof by considering two cases.}

In the first case, \(g_{\boldsymbol{x}^*}(v) \equiv 0\). Moreover, if there exist \(i\neq j\) such that \(x^*_i \neq x_j^*\), then we define function \(l(\epsilon)\) as in \eqref{eq:two-x} and denote \(\boldsymbol{x}(\epsilon)\) as the corresponding assignment. By earlier analysis, we have \(l''(0) < 0\). Also since \(g_{\boldsymbol{x}^*}(v) \equiv 0\), we have \(l'(0) = 0\). Moreover, since \(x^*_i \neq x^*_j\), there exists \(\delta \neq 0\) such that when \(\epsilon = \delta\), \(\boldsymbol{x}(\epsilon)\) is feasible and \(l(\delta) < l(0)\), which contradicts with the optimality of \(\boldsymbol{x}^*\). Therefore in this case, \(\boldsymbol{x}\) must all take the same value. Since \(\boldsymbol{x} \neq \boldsymbol{0}\) and \(\boldsymbol{x}\neq\boldsymbol{1}\), we know \(0 < x_i < 1, \;i = 1, \dots, n\), which can only happen when there is only one type of customers. Thus the conclusion is justified. When \(g_{\boldsymbol{x}^*}(v) \not\equiv 0\), it is either a quadratic function or an affine function, and thus intersects with the positive half of \(X\)-axis for at most twice. As a result, we can partition \((0, +\infty)\) into intervals \((0, \underline{\mu}^* ] \cup [\underline{\mu}^*, \bar{\mu}^* ] \cup [\bar{\mu}^* , +\infty)\) for some \(\underline{\mu}^* \le \bar{\mu}^*\). Within each interval, signs of \(g_{\boldsymbol{x}^* }(v)\) are the same, while in adjacent intervals they are different. In other words, there exist \(0 \le l < h \le n\) such that either \(x^*_i = \begin{cases}
    0 &\text{if }l < i < h\\
    1 &\text{if }i > h \text{ or } i < l
\end{cases}\) or \(x^*_i = \begin{cases}
    1 &\text{if }l < i < h\\
    0 &\text{if }i > h \text{ or } i < l
\end{cases}\).

Now we prove that if \(\alpha \ge \frac{1}{2}\), then any optimal solution can be represented by the former case. Again, we prove by contradiction. Suppose there exist \(i < j < k\) such that \(g_{\boldsymbol{x}^*}(\frac{1}{\mu_i}) \ge 0, g_{\boldsymbol{x}^*}(\frac{1}{\mu_j}) \le 0, g_{\boldsymbol{x}^*}(\frac{1}{\mu_k}) \ge 0\), then it implies that
\begin{eqnarray}
\frac{\gamma s}{\alpha} > \frac{t}{1-\alpha} &&\textup{(the quadratic term of }g_{\boldsymbol{x}^*}(v)\textup{ is positive)}; \label{eq:qudra-coeff} \\
\frac{\beta s}{\alpha} > \frac{t}{1-\alpha} &&\textup{(the constant term of }g_{\boldsymbol{x}^*}(v)\textup{ is positive)}; \label{eq:const-coeff} \\
\frac{u^2s^2}{\alpha} < \frac{t^2}{1-\alpha} &&\textup{(the linear term of }g_{\boldsymbol{x}^*}(v)\textup{ is negative)}. \label{eq:linear-coeff}
\end{eqnarray}

Multiplying \eqref{eq:qudra-coeff} with \eqref{eq:const-coeff} and then dividing it by \eqref{eq:linear-coeff}, we have \(\alpha < \frac{1}{2}\), which contradicts with the assumption \(\alpha\ge\frac{1}{2}\). Thus, there exists \(0 \le l < h \le n\) such that \(x^*_i = \begin{cases}0 &\text{if }l < i < h\\ 1&\text{if }i > h\text{ or }i < l \end{cases}\).\(\hfill\square\)

\bigskip
{\begin{theorem}
\label{thm:fix-alpha-convex}
Given \(\lambda_i, i = 1, 2, \dots, n, \,\, \alpha \ge \frac{1}{2}\) and \(\mathcal{I}\) with \(\mu_1 > \cdots > \mu_n, v_1, \dots, v_n\) satisfying $\frac{v_n - v_{n-1}}{1/\mu_n - 1/\mu_{n-1}} \ge \cdots \ge \frac{v_2 - v_1}{1/\mu_2 - 1/\mu_1} \ge \frac{v_1}{1/\mu_1}$. For any optimal solution \(\boldsymbol{x}^*\) to the QAP, there exist \(0 \le l < h \le n\) such that:
\begin{itemize}
\item When \(i > h\) or \(i < l\), \(x^*_i = 1\);
\item When \(l < i < h\), \(x^*_i = 0\);
\item When $ i \in \{l,h\} $, $ x^*_i \in [0,1] $. Moreover, the set \(\mathcal{M}^* = \left\{i: 0 < x^*_i < 1\right\}\) is either empty or singleton.
\end{itemize}
\end{theorem}

\noindent \textbf{Proof of Theorem \ref{thm:fix-alpha-convex}.}
The proof is similar to Theorem \ref{thm:fix-alpha-exp}. We first prove that the set \(\mathcal{M}^* = \left\{i: 0 < x^*_i < 1\right\}\) is either empty or singleton. Recall the $l$ function defined in \eqref{eq:two-x}, we only need to show that $l''(0) < 0$. We notice that
\[l''(0) = (\mu_i - \mu_j)(\mu_iv_i - \mu_jv_j) \frac{1}{2\sum_{k=1}^n \lambda_k} \left(\frac{1}{\alpha^2 - \alpha\sum_{k=1}^n \frac{\lambda_k x^*_k}{\mu_k}} + \frac{1}{(1 - \alpha)^2 - (1 - \alpha)\sum_{k=1}^n \frac{\lambda_k(1 - x^*_k)}{\mu_k}} \right).\]
Since the condition $\frac{v_n - v_{n-1}}{1/\mu_n - 1/\mu_{n-1}} \ge \cdots \ge \frac{v_2 - v_1}{1/\mu_2 - 1/\mu_1} \ge \frac{v_1}{1/mu_1}$ implies $\frac{v_n}{1/\mu_n } \ge \cdots \ge \frac{v_21}{1/\mu_2} \ge \frac{v_1}{1/mu_1}$. We know that $(\mu_i - \mu_j)(\mu_iv_i - \mu_jv_j)$ is negative. Thereby $l''(0) < 0$.

We next prove the rest of the statement. We follow the notation in Theorem \ref{thm:fix-alpha-exp}, except for $A_x = \sum_{i=1}^n\lambda_i\frac{v_i}{2}x_i$ and $A_y = \sum_{i=1}^n\lambda_i\frac{v_i}{2}(1-x_i)$. Recall the definition of $g_{\boldsymbol{x}}$ function in \eqref{eq-g-func} and we show in Theorem \ref{thm:fix-alpha-exp} that $g_{\boldsymbol{x}^*}$ is at most a quadratic function, thus intersecting with the positive half of X-axis for at most twice. Similarly we define
$$
g_{\boldsymbol{x}}(u) = \left( \frac{\gamma s}{\alpha} - \frac{t}{1-\alpha} \right) \frac{B_y}{A_y}\cdot \frac{1}{2} v_{\mathcal{I}}(u) + \left(\frac{\beta s}{\alpha}-\frac{t}{1-\alpha}\right) \frac{B_y}{C_y} + \left(\frac{u^2s^2}{\alpha}-\frac{t^2}{1-\alpha}\right) u,
$$
where $v_{\mathcal{I}}(u)$ is defined in \eqref{eq:v-u-relation}. By condition $\frac{v_n - v_{n-1}}{1/\mu_n - 1/\mu_{n-1}} \ge \cdots \ge \frac{v_2 - v_1}{1/\mu_2 - 1/\mu_1} \ge \frac{v_1}{1/mu_1}$, we know that $v_{\mathcal{I}}(u)$ is a convex function. Therefore, $g_{\boldsymbol{x}^*}$ is either convex or concave. Again it intersects with the positive half of X-axis for at most twice. Following the argument in the proof of Theorem \ref{thm:fix-alpha-exp}, there exists \(0 \le l < h \le n\) such that \(x^*_i = \begin{cases}0 &\text{if }l < i < h\\ 1&\text{if }i > h\text{ or }i < l \end{cases}\).\(\hfill\square\)}

\bigskip
{\begin{lemma}
\label{lemma:2seg}
Given \(\lambda_i, i = 1, 2, \dots, n\) and \(\mathcal{I}\) with \(\mu_1 > \cdots > \mu_n, v_1, \dots, v_n\) satisfying Assumption \ref{as:reg-2seg}. For any optimal solution \((\boldsymbol{x}^*, \alpha^*)\) to the queue partition problem, there exists \(i^*\) such that either \(x^*_i = \begin{cases}
	1 &\text{if }i < i^*\\
	0 &\text{if }i > i^* 
\end{cases}\) or \(x^*_i = \begin{cases}
	0 &\text{if }i < i^*\\
	1 &\text{if }i > i^*
\end{cases}\).
\end{lemma}}

\noindent \textbf{Proof of Lemma \ref{lemma:2seg}.}
{We prove it by showing that otherwise the optimality condition of $\alpha$ and $\boldsymbol{x}$ cannot hold at the same time.

We start with the optimality condition of \(\alpha\)}. Specifically, recall the definition of \(s, t\), we have
$$B_x = \frac{s}{1 + s}\alpha, \quad B_y = \frac{t}{1 + t}(1 - \alpha).$$

Also recall \(\gamma = \frac{C_x}{B_x}\frac{B_y}{C_y}, \beta = \frac{A_x}{B_x}\frac{B_y}{A_y}, u^2= \beta\gamma\). Now we can simplify the optimality condition of \(\alpha\) as
\begin{eqnarray*}
&\frac{A_x C_x}{\alpha(\alpha - B_x)}\left(\frac{1}{\alpha} + \frac{1}{\alpha - B_x}\right) = \frac{A_y C_y}{(1 - \alpha)(1 - \alpha - B_y)}\left(\frac{1}{1 - \alpha} + \frac{1}{1 - \alpha - B_y}\right) \\
\Leftrightarrow& \frac{\frac{C_xB_y}{C_yB_x}\frac{A_xB_y}{A_yB_x}}{\alpha(\alpha - B_x)/B_x^2}(\frac{1}{\alpha} + \frac{1}{\alpha - B_x}) = \frac{1}{(1 - \alpha)(1 - \alpha - B_y)/B_y^2}(\frac{1}{1 - \alpha} + \frac{1}{1 - \alpha - B_y}) \\
\Leftrightarrow& \frac{u^2}{\alpha(\alpha - B_x)/B_x^2} \left(\frac{1}{\alpha} + \frac{1}{\alpha - B_x}\right) = \frac{1}{(1 - \alpha)(1 - \alpha - B_y)/B_y^2} \left(\frac{1}{1 - \alpha} + \frac{1}{1 - \alpha - B_y}\right)\\
\Leftrightarrow& \frac{u^2}{\frac{1}{1+s}/(\frac{s}{1+s})^2} \left(\frac{1}{\alpha} + \frac{1+s}{\alpha}\right) = \frac{1}{\frac{1}{1+t}/(\frac{t}{1+t})^2} \left(\frac{1}{1 - \alpha} + \frac{1 + t}{1 - \alpha}\right)
\end{eqnarray*}
where the second to last line is by the definition of \(u\), and the last line replaces \(B_x, B_y\) with \(B_x = \frac{s}{1 + s}\alpha\) and \(B_y = \frac{t}{1 + t}(1 - \alpha)\).

By rearranging terms, we further have
\begin{eqnarray}\label{eq:general-kkt-alpha}
\frac{\alpha}{1-\alpha} = \frac{u^2s^2 \frac{s+2}{s+1}}{t^2 \frac{t+2}{t+1}}.
\end{eqnarray}

{Next we show that if the argument is wrong, contradiction arises from the optimality condition of $\boldsymbol{x}$.} In particular, by the proof of Theorem \ref{thm:fix-alpha}, if the argument is wrong, then \(v_{\mathcal{I}}(u)\) must intersect with linear function \(l(u) = -\frac{\left(\frac{\beta s}{\alpha}-\frac{t}{1-\alpha}\right) \frac{B_y}{C_y} + \left(\frac{u^2s^2}{\alpha}-\frac{t^2}{1-\alpha}\right) u}{\left( \frac{\gamma s}{\alpha} - \frac{t}{1-\alpha} \right) \frac{B_y}{A_y}}\) for at least twice. Since \(v_{\mathcal{I}}(u)\) satisfies Assumption \ref{as:reg-2seg}, we know that either \(\frac{\gamma s}{\alpha} - \frac{t}{1-\alpha} < 0, \frac{\beta s}{\alpha}-\frac{t}{1-\alpha} < 0\), or \(\frac{\gamma s}{\alpha} - \frac{t}{1-\alpha} > 0, \frac{\beta s}{\alpha}-\frac{t}{1-\alpha} > 0\). Define \(\theta = \frac{s \frac{s+2}{s+1}}{t \frac{t+2}{t+1}}\). Substituting \(\alpha\) by \eqref{eq:general-kkt-alpha}, we have either \(\gamma < u^2\theta, \beta < u^2\theta\), or \(\gamma > u^2\theta, \beta > u^2\theta\). Recall that \(u^2 = \beta\gamma\), we then have
\begin{eqnarray}\label{eq:general-gamma-plus-beta}
\gamma + \beta < u^2\theta + \frac{1}{\theta}.
\end{eqnarray}

We multiply both sides of \eqref{eq:general-kkt-x-sim} with \(\lambda_i\), which yields
\begin{eqnarray}\label{eq:general-kkt-x-two-part}
\left(\frac{\gamma s}{\alpha} - \frac{t}{1-\alpha}\right) \frac{B_y}{A_y}\lambda_i v_i
+ \left(\frac{\beta s}{\alpha}-\frac{t}{1-\alpha}\right) \frac{B_y}{C_y}\lambda_i
+ \left(\frac{u^2s}{\alpha}-\frac{t}{1-\alpha}\right) \frac{\lambda_i}{\mu_i}
= \frac{B_y^2}{A_yC_y}(p_i - q_i).
\end{eqnarray}

By the complementary slackness conditions, we know that \(x_i(p_i - q_i) = -q_i \le 0\). Thus, multiplying \eqref{eq:general-kkt-x-two-part} by \(x_i\) and summing up yields
\begin{eqnarray}\label{eq:general-kkt-x-le-zero}
\left(\frac{\gamma s}{\alpha} - \frac{t}{1 - \alpha}\right) A_x \frac{B_y}{A_y}
+ \left(\frac{\beta s}{\alpha} - \frac{t}{1-\alpha}\right) C_x \frac{B_y}{C_y}
+ \left(\frac{u^2s^2}{\alpha} - \frac{t^2}{1-\alpha}\right) B_x \le 0.
\end{eqnarray}
Here we use the fact that \(A_x = \sum_i \lambda_iv_ix_i,
B_x = \sum_i \frac{\lambda_i}{\mu_i}x_i ,
C_x = \sum_i \lambda_i x_i\). Next, we further have
\begin{eqnarray*}
\eqref{eq:general-kkt-x-le-zero}
&\Rightarrow& \left(\frac{\gamma s}{\alpha} - \frac{t}{1 - \alpha}\right) \beta
+ \left(\frac{\beta s}{\alpha} - \frac{t}{1-\alpha}\right) \gamma
+ \left(\frac{u^2s^2}{\alpha} - \frac{t^2}{1-\alpha}\right) \le 0 \\
&\Rightarrow& 2\beta\gamma + u^2s \le u^2 \theta(\beta+\gamma+t)\\
&\Rightarrow& 2u^2 + u^2s < u^2 \theta(u^2\theta+\frac{1}{\theta}+t) \\
&\Rightarrow& u^2\theta^2 > 1+s-t\theta\\
&\Rightarrow& u^2\theta^2 > \frac{(s+1)^2+(t+1)}{(t+2)(s+1)}
\end{eqnarray*}
where the second line is by \eqref{eq:general-kkt-alpha}, the third line is by \eqref{eq:general-gamma-plus-beta}, and the last line is by the definition of \(\theta\).

Similarly, by the complementary slackness conditions, \((1-x_i)(p_i - q_i) = p_i \ge 0\). Thus, by multiplying \eqref{eq:general-kkt-x-two-part} by \(1 - x_i\) and summing up yields
\begin{eqnarray}\label{eq:general-kkt-x-ge-zero}
\left(\frac{\gamma s}{\alpha} - \frac{t}{1 - \alpha}\right)A_y \frac{B_y}{A_y}+ \left(\frac{\beta s}{\alpha} - \frac{t}{1-\alpha}\right)C_y \frac{B_y}{C_y}+ \left(\frac{u^2s^2}{\alpha} - \frac{t^2}{1-\alpha}\right)B_y \ge 0.
\end{eqnarray}
Here we use the fact that \(A_y = \sum_i \lambda_iv_i(1-x_i),
B_y = \sum_i \frac{\lambda_i}{\mu_i}(1-x_i) ,
C_y = \sum_i \lambda_i(1-x_i)\). And we further have
\begin{eqnarray*}
\eqref{eq:general-kkt-x-ge-zero}
&\Rightarrow& \beta + \gamma + u^2s \ge u^2\theta(t+2)  \\
&\Rightarrow& u^2\theta + \frac{1}{\theta} + u^2s > u^2\theta(t+2)\\
&\Rightarrow& \frac{1}{u^2\theta^2} > t + 1 - \frac{s}{\theta} = \frac{(t+1)^2+(s+1)}{(s+2)(t+1)}
\end{eqnarray*}
where the first line is by \eqref{eq:general-kkt-alpha}, the second line is by \eqref{eq:general-gamma-plus-beta}, and the last line is by definition of \(\theta\).

Together, it implies
$$\frac{(s+1)^2+(t+1)}{(t+2)(s+1)} \cdot \frac{(t+1)^2+(s+1)}{(s+2)(t+1)} < 1,$$
which is equivalent to
$$\left((s+1) + (t+1)\right) (s - t)^2 < 0.$$

Therefore, we reach a contradiction. Thus there can not be \(i < j < k\) such that \(x^*_i > x^*_j\) and \(x^*_k > x^*_j\). It means that for any optimal solution \(\boldsymbol{x}^*\) there exists \(i^*\) such that for \(i \neq i^*\) either \(x^*_i = \begin{cases}
	1 &\text{if }i < i^*\\
	0 &\text{if }i > i^* 
\end{cases}\) or \(x^*_i = \begin{cases}
	0 &\text{if }i < i^*\\
	1 &\text{if }i > i^*
\end{cases}\).\(\hfill\square\)

\bigskip
{\begin{lemma}
\label{lemma:int}
Given \(\lambda_i, i = 1, 2, \dots, n\) and \(\mathcal{I}\) with \(\mu_1 > \cdots > \mu_n, v_1, \dots, v_n\) satisfying Assumption \ref{as:reg-2seg}. Any solution \((1, \dots, 1, x_i, 0, \dots, 0), \alpha\) with \(\alpha \in (0, 1)\) and \(x_i \in (0, 1)\) for some $i$ cannot be optimal.
\end{lemma}}

\noindent \textbf{Proof of Lemma \ref{lemma:int}.}
{We prove it by contradiction and show that if one type of customers are assigned to both queues, either transferring the customers from one queue to the other, or transferring the customers back from the other queue, together with the transfer of serving capacity, will always decrease the expected waiting time. Consequently, it cannot be optimal.}

We assume there exist at least one 1 and one 0 in \(\boldsymbol{x}\), since we can split customers with respect to \(x_i\) to construct artificial customer types. As a result, \(\mu_1 \ge \cdots \ge \mu_i \ge \cdots \ge \mu_n\).

{Before proceeding, we first derive some useful relations from optimality conditions of $\boldsymbol{x}$ and $\alpha$.} Specifically, we denote
\begin{eqnarray}\label{eq:general-beta-gamma-xy}
\beta_x  = \frac{A_x}{B_x}\frac{B_i}{A_i}, \gamma_x = \frac{C_x}{B_x}\frac{B_i}{C_i},
\beta_y  = \frac{A_y}{B_y}\frac{B_i}{A_i}, \gamma_y = \frac{C_y}{B_y}\frac{B_i}{C_i}.
\end{eqnarray}
Since \(\mu_1 \ge \cdots \ge \mu_i \ge \cdots \ge \mu_n\) satisfy Assumption \ref{as:reg-2seg}, we know
\begin{eqnarray}\label{eq:general-beta-gamma-order}
\beta_x \le 1 \le \beta_y,
\gamma_y \le 1 \le \gamma_x.
\end{eqnarray}

By first order optimality condition of \(\alpha\), we have
\begin{eqnarray}\label{eq:general-kkt-alpha-no-frac}
&&\frac{A_x C_x}{\alpha(\alpha - B_x)} \left(\frac{1}{\alpha} + \frac{1}{\alpha - B_x}\right) = \frac{A_y C_y}{(1 - \alpha)(1 - \alpha - B_y)} \left(\frac{1}{1 - \alpha} + \frac{1}{1 - \alpha - B_y}\right) \nonumber\\
&\Leftrightarrow& \frac{A_x / B_x \cdot C_x / B_x}{\alpha \cdot (\alpha - B_x) / B_x}\left(\frac{B_x}{\alpha} + \frac{B_x}{\alpha - B_x}\right) = \frac{A_y / B_y \cdot C_y / B_y}{(1 - \alpha) \cdot (1 - \alpha - B_y) / B_y} \left(\frac{B_y}{1 - \alpha} + \frac{B_y}{1 - \alpha - B_y}\right)\nonumber\\
&\Leftrightarrow& \frac{A_i C_i}{B_i^2}\frac{\beta_x \gamma_x}{\alpha / s} \left(\frac{s}{s+1} + s\right) = \frac{A_i C_i}{B_i^2} \frac{\beta_y \gamma_y}{(1-\alpha) / t} \left(\frac{t}{t+1} + t \right) \nonumber\\
&\Leftrightarrow& \frac{\alpha}{1 - \alpha} = \frac{\beta_x \gamma_x \frac{s^2(s+2)}{s+1}}{\beta_y \gamma_y \frac{t^2(t+2)}{t+1}}
\end{eqnarray}
where the third line is from \eqref{eq:def-st} and \eqref{eq:general-beta-gamma-xy}.

By first order optimality condition of \(x_i\), we have
\begin{eqnarray*}
&&\left(\frac{C_x}{\alpha(\alpha - B_x)} - \frac{C_y}{(1-\alpha)(1-\alpha-B_y)}\right)A_i \\
&+&\left(\frac{A_x}{\alpha(\alpha - B_x)} - \frac{A_y}{(1-\alpha)(1-\alpha-B_y)}\right)C_i \\
&+&\left(\frac{A_xC_x}{\alpha(\alpha - B_x)^2} - \frac{A_yC_y}{(1-\alpha)(1-\alpha-B_y)^2}\right)B_i = 0 \\
\Leftrightarrow &&\left(\frac{C_x / B_x}{\alpha(\alpha - B_x)/B_x} - \frac{C_y/B_y}{(1-\alpha)(1-\alpha-B_y)/B_y}\right)A_i \\
&+&\left(\frac{A_x / B_x}{\alpha(\alpha - B_x) / B_x} - \frac{A_y / B_y}{(1-\alpha)(1-\alpha-B_y) / B_y}\right)C_i\\
&+&\left(\frac{A_x / B_x \cdot C_x / B_x}{\alpha(\alpha - B_x)^2 / B_x^2} - \frac{A_y / B_y \cdot C_y / B_y}{(1-\alpha)(1-\alpha-B_y)^2 / B_y^2}\right)B_i = 0 \\
\Leftrightarrow&& \left(\frac{\gamma_x}{\alpha/s} - \frac{\gamma_y}{(1-\alpha)/t}\right)A_i C_i / B_i \\
&+&\left(\frac{\beta_x}{\alpha / s} - \frac{\beta_y}{(1-\alpha) / t}\right)A_i C_i / B_i \\
&+&\left(\frac{\beta_x \gamma_x}{\alpha/ s^2} - \frac{\beta_y \gamma_y}{(1-\alpha) / t^2}\right)A_i C_i / B_i = 0\\
\Leftrightarrow&&
\frac{\gamma_x}{\alpha/s}
+ \frac{\beta_x}{\alpha / s}
+ \frac{\beta_x \gamma_x}{\alpha/ s^2}
= \frac{\gamma_y}{(1-\alpha)/t}
+ \frac{\beta_y}{(1-\alpha) / t}
+ \frac{\beta_y \gamma_y}{(1-\alpha) / t^2} \\
\Leftrightarrow&&
\frac{\beta_x \gamma_x s}{\alpha}
\left(\frac{1}{\beta_x}
+ \frac{1}{\gamma_x}
+ s\right)
= \frac{\beta_y \gamma_y t}{1 - \alpha}
\left(\frac{1}{\beta_y}
+ \frac{1}{\gamma_y}
+ t\right)
\end{eqnarray*}
where the third to last line is from \eqref{eq:def-st} and \eqref{eq:general-beta-gamma-xy}. Substituting \(\alpha\) by \eqref{eq:general-kkt-alpha-no-frac}, we have
\begin{eqnarray}\label{eq:general-kkt-no-frac}
\frac{ \frac{1}{\beta_x} + \frac{1}{\gamma_x} +  s }{ \frac{s(s+2)}{s+1} }
= \frac{ \frac{1}{\beta_y} + \frac{1}{\gamma_y} + t }{\frac{t(t+2)}{t+1}}.
\end{eqnarray}

{To complete the proof, we construct a way to transfer both $x_i$ customers and serving capacity $\alpha$ such that expected total waiting time must decrease.} Specifically, we denote \(r = \frac{\frac{A_i}{B_i A_y} + \frac{C_i}{B_i C_y} + \frac{1}{1-\alpha-B_y}}{\frac{1}{1-\alpha}+\frac{1}{1-\alpha-B_y}}\). Multiply numerator and denominator with \(B_y\), we have \(r = \frac{ \frac{1}{\beta_y} + \frac{1}{\gamma_y} + t }{\frac{t(t+2)}{t+1}}\). Together with \eqref{eq:general-kkt-no-frac}, we have
$$r = 
\frac{ \frac{1}{\beta_x} + \frac{1}{\gamma_x} +  s }{ \frac{s(s+2)}{s+1} }
= \frac{ \frac{1}{\beta_y} + \frac{1}{\gamma_y} + t }{\frac{t(t+2)}{t+1}}.
$$

Now consider partition \((1,\cdots,1,x_i+\frac{1}{B_i}\epsilon,0,\cdots,0), \alpha + r\epsilon\). Note that as long as \(|\epsilon| > 0\) is small, it is a feasible partition. Multiplied by \(\sum_{i=1}^n \lambda_i\), the corresponding objective function is then denoted as
$$F(\epsilon) = G(\epsilon) + H(\epsilon)$$
where \(G(\epsilon) = \frac{\left(A_x + \frac{A_i}{B_i}\epsilon\right)\left(C_x + \frac{C_i}{B_i}\epsilon\right)}{(\alpha+r\epsilon)\left(\alpha-B_x+(r-1)\epsilon\right)}\) and \(H(\epsilon) = \frac{\left(A_y - \frac{A_i}{B_i}\epsilon\right)\left(C_y - \frac{C_i}{B_i}\epsilon\right)}{(1-\alpha-r\epsilon)(1-\alpha-B_y-(r-1)\epsilon)}\). {We next show that $F'(0) = 0$ and $F''(0) \le 0$.}

Denote \(g(\epsilon) = \log G(\epsilon), h(\epsilon) = \log H(\epsilon)\). Then \(F'(0) = G(0)g'(0) + H(0)h'(0)\). By \(r = \frac{\frac{A_i}{B_i A_y} + \frac{C_i}{B_i C_y} + \frac{1}{1-\alpha-B_y}}{\frac{1}{1-\alpha}+\frac{1}{1-\alpha-B_y}}\), we have \(h'(0) = \frac{-A_i}{B_i A_y} + \frac{-C_i}{B_i C_y} + \frac{r}{1-\alpha} + \frac{r-1}{1-\alpha-B_y} = 0\). By first order optimality condition, we know \(F'(0) = 0\), thus \(g'(0) = h'(0) = 0\). Note that
\begin{eqnarray*}
g(\epsilon) &=&
\log\left(A_x + \frac{A_i}{B_i}\epsilon\right)
+ \log\left(C_x + \frac{C_i}{B_i}\epsilon\right)
\log(\alpha+r\epsilon)
-\log\left(\alpha-B_x+(r-1)\epsilon\right); \\
g'(\epsilon) &=&
\frac{A_i / B_i}{A_x + \frac{A_i}{B_i}\epsilon}
+ \frac{C_i / B_i}{C_x + \frac{C_i}{B_i}\epsilon}
\frac{r}{\alpha+r\epsilon}
\frac{r-1}{\alpha-B_x+(r-1)\epsilon};\\
g''(\epsilon) &=&
\frac{-A_i^2 / B_i^2}{\left(A_x + \frac{A_i}{B_i}\epsilon\right)^2}
+ \frac{-C_i^2 / B_i^2}{\left(C_x + \frac{C_i}{B_i}\epsilon\right)^2}
+\frac{r^2}{(\alpha+r\epsilon)^2}
+\frac{(r-1)^2}{(\alpha-B_x+(r-1)\epsilon)^2}.
\end{eqnarray*}
Thus
\begin{eqnarray*}
g''(0) &=&
\frac{-A_i^2 / B_i^2}{A_x^2}
+ \frac{-C_i^2 / B_i^2}{C_x^2}
+\frac{r^2}{\alpha^2}
+\frac{(r-1)^2}{(\alpha-B_x)^2}\\
\Leftrightarrow B_x^2 g''(0) &=&
-\frac{A_i^2 B_x^2}{A_x^2 B_i^2}
-\frac{C_i^2 B_x^2}{C_x^2 B_i^2}
+\frac{r^2 B_x^2}{\alpha^2}
+\frac{(r-1)^2 B_x^2}{(\alpha-B_x)^2}\\
\Leftrightarrow B_x^2 g''(0) &=& r^2\left(\frac{s}{s+1}\right)^2 + (r-1)^2 s^2 - \left(\frac{1}{\beta_x^2} + \frac{1}{\gamma_x^2}\right)\\
&=& \frac{
-2(s+1)\frac{1}{\beta_x^2}
+\frac{2}{\beta_x}\left(\frac{1}{\gamma_x}(s^2+2s+2)-s^2\right)
+\left(\frac{1}{\gamma_x}+s\right)^2 + \left(\frac{1}{\gamma_x}(s+1)-s\right)^2 - (s+2)^2\frac{1}{\gamma_x}
}{(s+2)^2}
\end{eqnarray*}
where the third line is from \eqref{eq:general-beta-gamma-xy} and \eqref{eq:def-st}. Note that it is a quadratic function of \(\frac{1}{\beta_x}\) whose maximum is obtained at \(\frac{\frac{1}{\gamma_x}(s^2+2s+2)-s^2}{2(s+1)}\). From \eqref{eq:general-beta-gamma-order}, we know \(\gamma_x \ge 1\), thereby \(\frac{\frac{1}{\gamma_x}(s^2+2s+2)-s^2}{2(s+1)} \le 1\). Thereby, such quadratic function is decreasing in \([1, +\infty)\). Since \(\frac{1}{\beta_x} \ge 1\), we know \(B_x^2 g''(0) \le \frac{-2(s+1)}{(s+2)^2} \left(\frac{1}{\gamma_x}-1\right)^2 \le 0\), where equality holds if and only if \(\beta_x = \gamma_x = 1\).

Similarly,
$$
B_y^2 h''(0) = \frac{
-2\left(t+1\right)\frac{1}{\beta_y^2}
+\frac{2}{\beta_y}\left(\frac{1}{\gamma_y}(t^2+2t+2)-t^2\right)
+\left(\frac{1}{\gamma_y}+t\right)^2 + \left(\frac{1}{\gamma_y}(t+1)-t\right)^2 - (t+2)^2\frac{1}{\gamma_y}
}{(t+2)^2}
$$
which is a quadratic function of \(\frac{1}{\beta_y}\) whose maximum is obtained at \(\frac{\frac{1}{\gamma_y}(t^2+2t+2)-t^2}{2(t+1)} \ge 1\). Thereby, such quadratic function is increasing in \([0, 1]\). Since \(\frac{1}{\beta_y} \le 1\), we know \(B_y^2 h''(0) \le \frac{-2(t+1)}{(t+2)^2} \left(\frac{1}{\gamma_y}-1\right)^2 \le 0\), where equality holds if and only if \(\beta_y = \gamma_y = 1\).

Thereby,
$$
F''(0) = G(0) \left(g''(0) + \left(g'(0)\right)^2\right) + H(0)\left(h''(0)+\left(h'(0)\right)^2\right)
= G(0) g''(0) + H(0) h''(0)
\le 0
$$
where equality holds if and only if \(\beta_x = \gamma_x = \beta_y = \gamma_y = 1\).

In the case of \(F''(0) < 0\), due to the openness of feasible domain of \(\epsilon\), \(F'(0) = 0\) and strict concavity of \(F(\epsilon)\) at \(\epsilon = 0\), there exists \(\delta \neq 0\) such that corresponding partition and assignment are feasible and \(F(0) > F(\delta)\). This contradicts with the optimality assumption.

In the case of \(F''(0) = 0\), it is equivalent to the case where there is only one type of customers, and the optimization problem turns into
\begin{eqnarray*}
\inf_{x, \alpha} &\quad \frac{\lambda v}{2}\left(\frac{1}{\frac{\alpha}{x}(\frac{\alpha}{x} - \frac{\lambda}{\mu})} + \frac{1}{\frac{1-\alpha}{1-x}(\frac{1-\alpha}{1-x} - \frac{\lambda}{\mu})}\right)\\
\text{s.t.} &\quad \frac{\alpha}{x} > \frac{\lambda x}{\mu}, \frac{1-\alpha}{1-x} > \frac{\lambda x}{\mu}, 0< x < 1.
\end{eqnarray*}

Note that if \((x, \alpha)\) is a feasible solution, then \((1 - x, 1 - \alpha)\) is also a feasible solution with the same objective function value. Thereby, without loss of generality, assume \(\alpha \le x\). Thus objective function \(\ge \frac{\lambda v}{2}\frac{1}{1 - \frac{\lambda}{\mu}}\). Note that equality holds if and only if \(\alpha = x = 0\) or \(\alpha = x = 1\). Thus single queue is the unique optimal solution.\(\hfill\square\)

\bigskip
\noindent \textbf{Proof of Theorem \ref{thm:opt-alpha}.}
{In Lemma \ref{lemma:2seg}, we prove the ``two continuous segments'' property, i.e., for any optimal solution \((\boldsymbol{x}^*, \alpha^*)\), there exists \(i^*\) such that either \(x^*_i = \begin{cases}
	1 &\text{if }i < i^*\\
	0 &\text{if }i > i^* 
\end{cases}\) or \(x^*_i = \begin{cases}
	0 &\text{if }i < i^*\\
	1 &\text{if }i > i^*
\end{cases}\).

By plugging in Lemma \ref{lemma:int}, we know that for any optimal solution \((\boldsymbol{x}^*, \alpha^*)\), the set \(I^* := \left\{i: 0 < x_i^* < 1\right\}\) is empty. Thereby, either \(x^*_{i^*} = 1\) or \(x^*_{i^*} = 0\). And we can merge it into the \(i < i^*\) case or the \(i > i^*\) case, which completes the proof.} \(\hfill\square\)
\bigskip
\begin{theorem}
\label{thm:np-hard}
The \emph{deterministic assignment problem} for given \(\alpha\):
\begin{eqnarray}
\textup{(DAP)}\quad\inf_{\boldsymbol{x}} &&\quad f(\boldsymbol{x}, \alpha) \nonumber \\
\text{s.t.}&&\quad \boldsymbol{x} \in \mathcal{F}_\alpha; \quad x_i\in\left\{0,1\right\}, \; i = 1, \dots, n,
\end{eqnarray}
is NP-hard.
\end{theorem}

\noindent \textbf{Proof of Theorem \ref{thm:np-hard}.} We prove it by reduction from the \emph{Set Partition} problem. Recall the \emph{Set Partition} problem takes a set \(S\) of numbers and outputs whether there exists a partition \(A\) and \(\bar{A}\) such that \(A \cup \bar{A} = S\) and \(\sum_{w\in A} w= \sum_{w\in\bar{A}} w\).

For set \(S\) in the \emph{Set Partition} problem with elements \(w_1 \le w_2 \le \cdots \le w_n\), we construct an instance of the DAP with parameters \(\lambda_i = w_i\), \(\mu_i = \bar{\mu} = 2\sum_{i=1}^n w_i\), for all \(i\) and \(\alpha = \frac{1}{2}\). To complete the reduction, we next prove
\begin{itemize}
\item If equal set partition exists, then all optimal solutions to the DAP satisfy \(\sum_{i=1}^n \lambda_ix_i = \frac{1}{2}\sum_{i=1}^n \lambda_i\).
\item If equal set partition does not exist, then for any optimal solution to the DAP, \(\sum_{i=1}^n \lambda_ix_i \neq \frac{1}{2}\sum_{i=1}^n \lambda_i\).
\end{itemize}

Note that if the above statements hold, then we can solve the \emph{Set Partition} problem by first solving the constructed DAP and then checking whether the optimal solution satisfies the condition.

Now we consider the DAP with the constructed parameters. Multiply the objective function with \(\sum_{i=1}^n \lambda_i\), the DAP can be reformulated as
\begin{eqnarray*}
\inf_{\boldsymbol{x}} &&\quad \frac{(\sum_{i=1}^n \lambda_i x_i)^2}{\frac{\bar{\mu}}{2} - \sum_{i=1}^n \lambda_i x_i} + \frac{(\sum_{i=1}^n \lambda_i - \sum_{i=1}^n \lambda_ix_i)^2}{\frac{\bar{\mu}}{2} - (\sum_{i=1}^n \lambda_i - \sum_{i=1}^n \lambda_i x_i)}\\
\text{s.t.} &&\quad -\frac{\bar{\mu}}{2} + \sum_{i=1}^n \lambda_i < \sum_{i=1}^n \lambda_i x_i < \frac{\bar{\mu}}{2};\\
&&\quad \boldsymbol{x} \in \left\{0, 1\right\}^n.
\end{eqnarray*}

If we drop the binary constraint and denote \(y = \sum_{i=1}^n \lambda_i x_i\), then it can be relaxed to
\begin{eqnarray*}
\inf_{y} &&\quad h(y) := \frac{y^2}{\frac{\bar{\mu}}{2} - y} + \frac{\left(\sum_{i=1}^n \lambda_i - y\right)^2}{\frac{\bar{\mu}}{2} - \left(\sum_{i=1}^n \lambda_i - y\right)}\\
\text{s.t.} &&\quad -\frac{\bar{\mu}}{2} + \sum_{i=1}^n \lambda_i < y < \frac{\bar{\mu}}{2}.
\end{eqnarray*}

Note that \(h''(y) > 0, y^* = \frac{1}{2}\sum_{i=1}^n \lambda_i = \frac{\bar{\mu}}{4}\) is feasible, and \(h'(y^*) = 0\). 
Thereby, \(y^*\) is the unique optimal solution to the relaxed problem.

Therefore, if there exists subset \(I\) such that \(\sum_{i\in I}\lambda_i = \frac{1}{2}\sum_{i=1}^n \lambda_i\), then \(y^*\) is also feasible for the original problem thus optimal. If there is no such \(I\), then \(\forall I, y_I := \sum_{i\in I}\lambda_i \neq \frac{1}{2}\sum_{i=1}^n \lambda_i\), which completes the proof.\(\hfill\square\)

\bigskip
\noindent \textbf{Proof of Theorem \ref{thm:multi-fixed-alpha}.} Since each element of \(\mathcal{C}(\mathcal{I}, k)\) can be represented by \(\mathcal{L}\) and some \(\mathcal{M}\), we prove that the optimal assignment can be represented by \(\mathcal{L}\) and some \(\mathcal{M}\). Start from the first queue. Following the proof of Theorem \ref{thm:fix-alpha}, we know there exist at most \(k - 1\) of lines to separate customers in this queue from customers in the other \(k - 1\) queues. And if \((\mathbb{E}[S_i], \mathbb{E}[S_i^2])\) is in the interior of the polygon, then \(X_{i1} = 1\). Those in the boundary correspond to \(X_{i1} \in [0, 1]\) and the others \(X_{i1} = 0\). Remove points in the interior of polygon, and continue the process. Thereby the optimal assignment can be induced by these \(\frac{(k-1)k}{2}\) lines which constitute the set \(\mathcal{L}\). \(\hfill\square\)

\bigskip
\noindent \textbf{Proof of Theorem \ref{thm:multi-exp}.} Start from the \(k\)th queue. Following the proof of Theorem \ref{thm:fix-alpha-exp}, we know there exist \(0 \le l_k^1 \le l_k^2 \le n + 1\) such that \(X_{ik} = 1\) if and only if \(l_k^1 < i < l_k^2\). Remove all \(i\)th type of customers such that \(l_k^1 \le i \le l_k^2\). Then consider the \((k-1)\)th queue. We can repeat the process for the rest of the customers until we are considering the first queue. For the first queue, let \(l_1^1 = 0, l_1^2 = n + 1\), obviously, the rest of customers must be assigned to the first queue. Now we have \(2k\) of \(l\)'s, which completes the proof.\(\hfill\square\)

\bigskip
\noindent \textbf{Proof of Theorem \ref{thm:multi-opt-alpha}.} Following the proof of Theorem \ref{thm:opt-alpha}, we know that for any two queues \(j_1, j_2\), an optimal solution \((X^*, \boldsymbol{\alpha}^*)\) splits customers in these two queues by service rate \(\mu_i\), such that there exists \(\bar{\mu}\) and customers whose service rate \(\mu_i \le \bar{\mu}\) are assigned to one queue while customers whose service rate \(\mu_i > \bar{\mu}\) are assigned to the other queue. Define interval \(I_j = [\underline{\mu}_j, \bar{\mu}_j]\) where \(\underline{\mu}_j = \min \left\{\mu_i: X^*_{ik} > 0, k = j\right\}, \bar{\mu}_j = \max \left\{\mu_i: X^*_{ik} > 0, k = j\right\}\). Then, we know for any \(1 \le j_1 < j_2 \le k\), it holds that \(I_{j_1} \cap I_{j_2} = \emptyset\). Thereby, the conclusion is justified.

\bigskip
\noindent \textbf{Proof of Theorem \ref{thm:multi-fixed-alpha-sojourn}.}
The idea is similar to the waiting time case. We only need to prove the \(k = 2\) case, since the general case proof can be extended following the same argument in the proof of Theorem \ref{thm:multi-fixed-alpha}.

First we consider \(\alpha = 1\) (\(\alpha = 0\)). In this case, \(\boldsymbol{x} = \boldsymbol{1}\) (\(\boldsymbol{x} = \boldsymbol{0}\)) is the only feasible solution, and the statement is true.

When \(0 < \alpha < 1\), we first show that \(\boldsymbol{x} = \boldsymbol 0\) or \(\boldsymbol{x} = \boldsymbol 1\), i.e., assigning all customers to one queue, cannot be optimal. We follow the notation of \(A_x, B_x, C_x, A_y, B_y, C_y, s, t\) in the proof of Theorem \ref{thm:fix-alpha} and denote $C, B, A$ as the zeroth, first, second order load of a pooled queue with unit capacity which serves all types of customers. Therefore,
\[A = \sum_{i=1}^n\lambda_iv_i, B = \sum_{i=1}^n \frac{\lambda_i}{\mu_i}, C = \sum_{i=1}^n \lambda_i.\] Then we can write \(\tilde{f}(\boldsymbol{x}, \alpha) = \frac{1}{2\sum_{i=1}^n \lambda_i}\left(\frac{A_x \cdot C_x}{\alpha(\alpha - B_x)} + \frac{2B_x}{\alpha} + \frac{A_y \cdot C_y}{(1 - \alpha)(1 - \alpha - B_y)} + \frac{2B_y}{1-\alpha}\right)\).

We denote the Lagrangian function \(\mathcal{L}(\boldsymbol{x}, \alpha, \boldsymbol{p}, \boldsymbol{q}) = \frac{A_x \cdot C_x}{\alpha(\alpha - B_x)} + \frac{2B_x}{\alpha} + \frac{A_y \cdot C_y}{(1 - \alpha)(1 - \alpha - B_y)} + \frac{2B_y}{1-\alpha} - \sum_i x_i p_i - \sum_i (1 - x_i) q_i\), where \(p_i \ge 0, q_i \ge 0\) (we remove the constant factor \(\frac{1}{2\sum_{i=1}^n \lambda_i}\) for the ease of expression). From the KKT condition \(\frac{\partial\mathcal{L}}{\partial x_i} = 0\), we have
\begin{eqnarray}\label{eq:sojourn-kkt-x}
p_i - q_i
=&& \left(\frac{C_x}{\alpha(\alpha - B_x)} - \frac{C_y}{(1-\alpha)(1-\alpha-B_y)}\right) \lambda_i v_i \nonumber \\
&&+ \left(\frac{A_x}{\alpha(\alpha - B_x)} - \frac{A_y}{(1-\alpha)(1-\alpha-B_y)}\right) \lambda_i \nonumber \\
&&+ \left(\frac{A_xC_x}{\alpha(\alpha - B_x)^2} + \frac{2}{\alpha} - \frac{A_yC_y}{(1-\alpha)(1-\alpha-B_y)^2} - \frac{2}{1-\alpha}\right) \frac{\lambda_i}{\mu_i}.
\end{eqnarray}
If \(\boldsymbol{x} = \boldsymbol 0\) is feasible, then for \(\boldsymbol{x} = \boldsymbol 0\),
\begin{align*}
p_i - q_i &= -\left(\frac{C_y}{(1-\alpha)(1-\alpha-B_y)} \lambda_i v_i + \frac{A_y}{(1-\alpha)(1-\alpha-B_y)}\lambda_i \right.\\
&\qquad\qquad\qquad\quad\left. + \left(\frac{A_yC_y}{(1-\alpha)(1-\alpha-B_y)^2} + \frac{2}{1-\alpha}\right)\frac{\lambda_i}{\mu_i}\right) \\&< 0.
\end{align*}
From non-negativity of \(p_i, q_i\) we know \(q_i > 0\). However, complementary slackness condition states that \((1 - x_i)q_i = 0\), thus \(\boldsymbol{x} = \boldsymbol 0\) is not optimal.

Similarly, if \(\boldsymbol{x} = \boldsymbol 1\) is feasible, then for \(\boldsymbol{x} = \boldsymbol 1\),
$$p_i - q_i = \frac{C_x}{\alpha(\alpha - B_x)}\frac{\lambda_i}{\mu_i^2} + \frac{A_x}{\alpha(\alpha - B_x)} \lambda_i + \left(\frac{A_xC_x}{\alpha(\alpha - B_x)^2} + \frac{2}{\alpha}\right) \frac{\lambda_i}{\mu_i} > 0.$$

From non-negativity of \(p_i, q_i\) we know \(p_i > 0\). However, complementary slackness condition states that \(x_i p_i = 0\), thus \(\boldsymbol{x} = \boldsymbol 1\) is not optimal.

Now we can restrict our consideration to \(\boldsymbol{x} \neq \boldsymbol{0}, \boldsymbol{1}\). We first formulate a representative function which maps each customer type to a scalar value. In particular, for any feasible assignment \(\boldsymbol{x}\), we define function \(g_{\boldsymbol{x}}\) as
\begin{equation}\label{eq:v-u-sojourn}
g_{\boldsymbol{x}}(u) = \left( \frac{\gamma s}{\alpha} - \frac{t}{1-\alpha} \right) \frac{B_y}{A_y} v_{\mathcal{I}}(u) + \left(\frac{\beta s}{\alpha}-\frac{t}{1-\alpha}\right) \frac{B_y}{C_y} + \left(\frac{2B_y^2}{A_yC_y}\left(\frac{1}{\alpha} - \frac{1}{1-\alpha}\right) + \frac{u^2s^2}{\alpha} - \frac{t^2}{1-\alpha} \right) u.
\end{equation}
It is readily shown that representative function only depends linearly on the first two moment of the customers' service time, i.e., $u$ and $v_{\mathcal{I}}(u)$.

We complete the proof by showing that, under the optimal assignment policy, customer of type $i$ will be allocated to the first server only when the corresponding scalar value $g_{\boldsymbol{x}}(\cdot)$ is non-positive. Consequently, due to linearity, the optimal assignment then corresponds to a polytope partition in a two-dimensional space of the first two moments. In particular, when there are two queues, the two-dimensional space are partitioned by a line into two subspaces.

Specifically, since \(\boldsymbol{x} \neq \boldsymbol{0}, \boldsymbol{1}\), we have \(A_x, A_y, B_x, B_y, C_x, C_y \neq 0\). Multiplying both sides of \eqref{eq:sojourn-kkt-x} by \(\frac{B_y^2}{A_yC_y}\frac{1}{\lambda_i}\) yields
\begin{eqnarray*}
\frac{B_y^2}{A_yC_y}\frac{p_i - q_i}{\lambda_i}
=&& \left(\frac{\frac{C_xB_y}{C_yB_x}}{\alpha(\alpha - B_x)/B_x} - \frac{1}{(1-\alpha)(1-\alpha-B_y)/B_y}\right) \frac{B_y}{A_y} v_i \nonumber \\
&&+ \left(\frac{\frac{A_xB_y}{A_yB_x}}{\alpha(\alpha - B_x)/B_x} - \frac{1}{(1-\alpha)(1-\alpha-B_y)/B_y}\right) \frac{B_y}{C_y} \nonumber \\
&&+ \frac{B_y^2}{A_yC_y}\left(\frac{A_xC_x}{\alpha(\alpha - B_x)^2} + \frac{2}{\alpha} - \frac{A_yC_y}{(1-\alpha)(1-\alpha-B_y)^2} - \frac{2}{1-\alpha}\right)\frac{1}{\mu_i}.
\end{eqnarray*}

Denote \(\gamma = \frac{C_x}{B_x}\frac{B_y}{C_y}, \beta = \frac{A_x}{B_x}\frac{B_y}{A_y}\). Recall that \(s = \frac{B_x}{\alpha-B_x}, t = \frac{B_y}{1 - \alpha-B_y}\), the above equation can then be simplified as
\begin{eqnarray}\label{eq:sojourn-kkt-x-sim}
\left(\frac{\gamma s}{\alpha} - \frac{t}{1-\alpha}\right) \frac{B_y}{A_y} v_i
+ \left(\frac{\beta s}{\alpha}-\frac{t}{1-\alpha}\right) \frac{B_y}{C_y}
+ \left(\frac{2B_y^2}{A_yC_y}\left(\frac{1}{\alpha} - \frac{1}{1-\alpha}\right) + \frac{u^2s^2}{\alpha} - \frac{t^2}{1-\alpha} \right) \frac{1}{\mu_i} \nonumber \\
\qquad\qquad = \frac{B_y^2}{A_yC_y}\frac{p_i - q_i}{\lambda_i}.
\end{eqnarray}

Recall definition \eqref{eq:v-u-sojourn}, we know that \eqref{eq:sojourn-kkt-x-sim} can be rewritten as \(g_{\boldsymbol{x}}(\frac{1}{\mu_i}) = \frac{B_y^2}{A_yC_y}\frac{p_i - q_i}{\lambda_i}\). By complementary slackness conditions, for any optimal assignment \(\boldsymbol{x}^*\), it holds that
$$
g_{\boldsymbol{x}^*}\left(\frac{1}{\mu_i}\right) > 0 \Rightarrow x^*_i = 0, \quad g_{\boldsymbol{x}^*}\left(\frac{1}{\mu_i}\right) < 0 \Rightarrow x^*_i = 1.
$$
We then consider two cases. In the first case, \(\frac{\gamma s}{\alpha} \neq \frac{t}{1-\alpha}\). Consider linear function \(l(u) = -\frac{\left(\frac{\beta s}{\alpha}-\frac{t}{1-\alpha}\right) \frac{B_y}{C_y} + \left(\frac{2B_y^2}{A_yC_y}\left(\frac{1}{\alpha} - \frac{1}{1-\alpha}\right)+\frac{u^2s^2}{\alpha}-\frac{t^2}{1-\alpha}\right) u}{\left( \frac{\gamma s}{\alpha} - \frac{t}{1-\alpha} \right) \frac{B_y}{A_y}}\) and its induced partition \((U, E, L)\). It is readily shown that \((U, E, L) \in \mathcal{C}(\mathcal{I})\). When \(\frac{\gamma s}{\alpha} > \frac{t}{1-\alpha}\), the KKT condition implies that an optimal solution \(\boldsymbol{x}^*\) satisfies
$$\begin{cases} x_i^{ * } = 0 & i\in U \\ x_i^{ * } = 1 & i\in L \\ x_i^{ * } \in [0, 1] & i\in E. \end{cases}$$
Similarly, when \(\frac{\gamma s}{\alpha} < \frac{t}{1-\alpha}\), an optimal solution \(\boldsymbol{x}^*\) satisfies
$$\begin{cases} x_i^{ * } = 1 & i\in U \\ x_i^{ * } = 0 & i\in L \\ x_i^{ * } \in [0, 1] & i\in E. \end{cases}$$
Otherwise, \(\frac{\gamma s}{\alpha} = \frac{t}{1-\alpha}\). Then \(g_{\boldsymbol{x}^*}(u)\) reduces to a linear function. When \(g_{\boldsymbol{x}^*}(u) \not\equiv 0\), \(g_{\boldsymbol{x}^*}(u)\) intersects with the positive half of \(X\)-axis for at most once. By KKT condition, we know there exists \(i_0 \in \{1, 2, \dots, n\}\) such that either \(\begin{cases} x_i^{ * } = 0 & i < i_0 \\ x_i^{ * } = 1 & i > i_0 \\ x_i^{ * } \in [0, 1] & i = i_0. \end{cases}\) or \(\begin{cases} x_i^{ * } = 1 & i < i_0 \\ x_i^{ * } = 0 & i > i_0 \\ x_i^{ * } \in [0, 1] & i = i_0. \end{cases}\). Construct \(U = \{i: i > i_0\}, E = \{i: i = i_0\}, L = \{i: i < i_0\}\). By considering a linear function passing through \((\frac{1}{\mu_{i_0}}, v_0)\) with large enough slope, we have \((U, E, L) \in \mathcal{C}(\mathcal{I})\). When \(g_{\boldsymbol{x}^*}(u) \equiv 0\), we have
\begin{align*}
\frac{\gamma s}{\alpha} = \frac{t}{1-\alpha},
\frac{\beta s}{\alpha} = \frac{t}{1-\alpha},
\frac{2B_y^2}{A_yC_y}\left(\frac{1}{\alpha} - \frac{1}{1-\alpha}\right) + \frac{u^2s^2}{\alpha} - \frac{t^2}{1-\alpha} = 0.
\end{align*}
Multiply the first two equations, we have
$$ 
\frac{u^2s^2}{\alpha^2} = \frac{t^2}{(1-\alpha)^2}.
$$
Without loss of generality assume \(\alpha \ge \frac{1}{2}\). Thereby, \(\frac{u^2s^2}{\alpha} \ge \frac{t^2}{1-\alpha}\). And we have \(\frac{A_x}{A_y} = \frac{C_x}{C_y} = \frac{\alpha}{1-\alpha}, \alpha - B_x = 1 - \alpha - B_y, \frac{A_xC_x}{\alpha(\alpha-B_x)}=\frac{2}{1-\alpha}\), which implies \(\frac{AC}{1-B} = \frac{1}{\alpha(1-\alpha)}\). Consider any \(\alpha_0\) in the neighborhood of \(\alpha\), \(\frac{AC}{1-B} \neq \frac{1}{\alpha_0(1-\alpha_0)}\). Thereby, the conclusion holds for \(\alpha_0\). Let \(\alpha_0 \to \alpha\), we can conclude that an optimal solution satisfying the condition exists for \(\alpha\), which completes the proof. \(\hfill\square\)

\end{APPENDIX}

%
%   or
%
% \begin{APPENDICES}
% \section{<Title of Section A>}
% \section{<Title of Section B>}
% etc
% \end{APPENDICES}

% Acknowledgments here
\ACKNOWLEDGMENT{The authors thank the Area Editor, the Associate Editor, and two referees for their diligent reading and thoughtful comments that significantly improved the paper.}

% References here (outcomment the appropriate case)

% CASE 1: BiBTeX used to constantly update the references
%   (while the paper is being written).
\bibliographystyle{informs2014} % outcomment this and next line in Case 1
\bibliography{paper} % if more than one, comma separated

%\bibliographystyle{informs2014} % outcomment this and next line in Case 1
%\bibliography{sample} % if more than one, comma separated

% CASE 2: BiBTeX used to generate mypaper.bbl (to be further fine tuned)
%\input{mypaper.bbl} % outcomment this line in Case 2

%If you don't use BiBTex, you can manually itemize references as shown below.

%\bibliographystyle{nonumber}

%\begin{thebibliography}{3}
%\providecommand{\natexlab}[1]{#1}
%\providecommand{\url}[1]{\texttt{#1}}
%\providecommand{\urlprefix}{URL }
%
%\bibitem[{Smith(2005)}]{smith2005}
%Smith J (2005) Optimal resource allocation in humanitarian logistics.
%  \emph{Journal of Operations Research} 30(2):123--135.
%  
%\bibitem[{Jones(2010)}]{jones2010}
%Jones S (2010) Stochastic programming models for humanitarian logistics.
%  \emph{INFORMS Mathematics of Operations Research} 35(4):567--580.
%
%\bibitem[{Brown(2015)}]{brown2015}
%Brown D (2015) \emph{Introduction to Stochastic Programming} (Springer).
%
%\end{thebibliography}

%%%%%%%%%%%%%%%%%
\end{document}